\documentclass[10pt]{amsart} 
\usepackage{amscd,amssymb}
\usepackage{enumerate} 
\usepackage[all]{xypic} 
\usepackage{multirow}
\usepackage{hhline}
\usepackage{comment}

\theoremstyle{plain} 
\newtheorem{Thm}[equation]{Theorem}
\newtheorem*{Thm*}{Theorem} 
\newtheorem{Cor}[equation]{Corollary}
\newtheorem{Con}[equation]{Conjecture}
\newtheorem{Prop}[equation]{Proposition}
\newtheorem{Lem}[equation]{Lemma} 
\newtheorem{Rmk}[equation]{Remark}

\numberwithin{equation}{section}

\newcommand{\Hom}{\operatorname{Hom}}

\newcommand{\GL}{\operatorname{GL}}
\newcommand{\PGL}{\operatorname{PGL}}

\newcommand{\SL}{\operatorname{SL}}

\newcommand{\SO}{\operatorname{SO}}

\newcommand{\OO}{\operatorname{O}}

\newcommand{\Ind}{\operatorname{Ind}}
\newcommand{\ind}{\operatorname{ind}}

\newcommand{\diag}{\operatorname{diag}}

\newcommand{\GLt}{\widetilde{\operatorname{GL}}}
\newcommand{\SLt}{\widetilde{\operatorname{SL}}}
 
\newcommand{\Bt}{\widetilde{B}}
\newcommand{\Tt}{\widetilde{T}} 
\newcommand{\Pt}{\widetilde{P}}
\newcommand{\Qt}{\widetilde{Q}} 

\newcommand{\MPt}{\widetilde{M_P}} 
\newcommand{\MQt}{\widetilde{M_Q}}
\newcommand{\e}{\operatorname{e}}

\newcommand{\GLtt}{\widetilde{\operatorname{GL}}^{(2)}}

\newcommand{\Bte}{{\widetilde{B}}^{\e}}
\newcommand{\Tte}{\widetilde{T}^{\e}}

\newcommand{\Ptt}{\widetilde{P}^{(2)}}
\newcommand{\Qtt}{\widetilde{Q}^{(2)}}

\newcommand{\MPtt}{\widetilde{M_P}^{(2)}}
\newcommand{\MPttt}{(\widetilde{M_P})^{(2)}}
\newcommand{\MQtt}{\widetilde{M_Q}^{(2)}}
\newcommand{\MQttt}{(\widetilde{M_Q})^{(2)}}

\newcommand{\sigGLt}{{^{\sigma}\widetilde{\operatorname{GL}}}}

\newcommand{\cMPt}{{^c\widetilde{M_P}} }
\newcommand{\cMQt}{{^c\widetilde{M_Q}}}

\newcommand{\cMPtt}{{^c\widetilde{M_P}^{(2)}}}
\newcommand{\cMPttt}{{^c(\widetilde{M_P})^{(2)}}}

\newcommand{\GLtwo}{\GL^{(2)}}

\newcommand{\Zt}{\widetilde{Z}} 
\newcommand{\Mt}{\widetilde{M}}

\newcommand{\timest}{\widetilde{\times}}
\newcommand{\otimest}{\widetilde{\otimes}}

\newcommand{\varphit}{\widetilde{\varphi}}
\newcommand{\pit}{\widetilde{\pi}}

\newcommand{\gt}{\tilde{g}}

\newcommand{\rr}{\mathbf{r}} 
\newcommand{\I}{\mathcal{I}}

\newcommand{\OF}{\mathcal{O}_F}

\newcommand{\G}{\mathbb{G}} 
\newcommand{\Gbt}{\widetilde{\mathbb{G}}} 

\newcommand{\Res}{\operatorname{Res}}

\newcommand{\rec}{\operatorname{rec}}
\renewcommand{\Re}{\operatorname{Re}}

\newcommand{\C}{\mathbb C}
\newcommand{\A}{\mathbb{A}}

\newcommand{\Z}{\mathbb{Z}}
\newcommand{\R}{\mathbb{R}}

\newcommand{\s}{\mathbf{s}}
\newcommand{\sss}{\mathfrak{s}}
\newcommand{\ie}{{\em i.e. }}

\title[symmetric square $L$-function]{The twisted symmetric square
$L$-function of $\GL(r)$} \author{Shuichiro Takeda}
\address{Department of Mathematics, University of Missouri,
  Columbia, MO 65211}

\begin{document}

\maketitle

\begin{abstract} 
In this paper, we consider the (partial) symmetric
square $L$-function $L^S(s,\pi,Sym^2\otimes\chi)$ of an irreducible
cuspidal automorphic representation $\pi$ of $\GL_r(\A)$ twisted by a
Hecke character $\chi$. In particular, we will show that the
$L$-function $L^S(s,\pi,Sym^2\otimes\chi)$ is holomorphic for the
region $\Re(s)>1-\frac{1}{r}$ with the exception that, if
$\chi^r\omega^2=1$, a pole might occur at $s=1$, where $\omega$ is the
central character of $\pi$. Our method of proof is essentially a
(nontrivial) modification of the one by Bump and Ginzburg in
which they considered the case $\chi=1$.
\end{abstract}

%%%%%%%%%%%%%%%%%%%%%%%%%%%%%%%%%%%%%%%%%%%%%%%%%%%%%%%%%%%%%%%%%%%

\section*{\bf Introduction}

%%%%%%%%%%%%%%%%%%%%%%%%%%%%%%%%%%%%%%%%%%%%%%%%%%%%%%%%%%%%%%%%%%%

Let $\pi\cong\otimes'_v\pi_v$ be an irreducible cuspidal automorphic
representation of $\GL_r(\A)$ and $\chi$ a unitary Hecke character on
$\A^\times$, where $\A$ is the ring of adeles over a number field
$F$. By the local Langlands correspondence by Harris-Taylor \cite{HT}
and Henniart \cite{He}, each $\pi_v$ corresponds to an $r$-dimensional
representation $\rec(\pi_v)$ of the Weil-Deligne group $WD_{F_v}$ of
$F_v$. We can also consider the twist of $\rec(\pi_v)$ by $\chi_v$, namely,
\[ \rec(\pi_v)\otimes\chi_v:WD_{F_v}\rightarrow \GL_r(\C),
\] where $\chi_v$ is viewed as a character of $WD_{F_v}$ via local
class field theory. Now for each homomorphism
\[ \rho:\GL_r(\C)\rightarrow\GL_N(\C),
\] one can associate the local $L$-factor
$L_v(s,\pi_v,\rho\circ\rec(\pi_v)\otimes\chi_v)$ of Artin type. Then
one can define the automorphic $L$-function by
\[
L(s,\pi,\rho\otimes\chi):=\prod_vL_v(s,\pi_v,\rho\circ\rec(\pi_v)\otimes\chi_v).
\] In particular in this paper, we consider the case where $\rho$ is
the symmetric square map
\[ Sym^2:\GL_r(\C)\rightarrow\GL_{\frac{1}{2}r(r+1)}(\C),
\] namely we consider the twisted symmetric square $L$-function
$L(s,\pi,Sym^2\otimes\chi)$. By the Langlands-Shahidi method, it can
be shown that the $L$-function $L(s,\pi,Sym^2\otimes\chi)$ admits
meromorphic continuation and a functional equation. 
(See \cite[Theorem 7.7]{Sh90}.)

The Langlands-Shahidi method, however, is unable to determine the
locations of the possible poles of $L(s,\pi,Sym^2\otimes\chi)$. The
main theme of this paper is to determine them to some extent,
though we consider only the partial $L$-function
$L^S(s,\pi,Sym^2\otimes\chi)$. To be more specific, let $S$ be the
finite set of places that contains all the archimedean places and
non-archimedean places where $\pi$ or $\chi$ ramifies. For $v\notin
S$, each $\pi_v$ is parameterized by a set of $r$ complex numbers
$\{\alpha_{v,1},\dots,\alpha_{v,r}\}$ known as the Satake
parameters. Then we have
\[ 
L_v(s,\pi_v,Sym^2\otimes\chi_v)=\prod_{i\leq
j}\frac{1}{(1-\chi_v(\varpi_v)\alpha_{v,i}\alpha_{v,j}q_v^{-s})},
\] 
where $\varpi_v$ is the uniformizer of $F_v$ and $q_v$ is the order
of the residue field. And we set
\[ 
L^S(s,\pi,Sym^2\otimes\chi)=\prod_{v\notin
S}L_v(s,\pi_v,Sym^2\otimes\chi_v).
\] 
As our main theorem (Theorem \ref{T:main}) we will prove

\quad\\
\noindent{\bf Theorem \ref{T:main}.} {\it Let $\pi$ be a cuspidal
automorphic representation of $\GL_r(\A)$ with unitary central character
$\omega$ and $\chi$ a unitary Hecke character. Then for each
archimedean $v$, there exists an integer $N_v\geq 0$ such that the
product
\[ L^S(s,\pi,Sym^2\otimes\chi)\prod_{v|\infty}L_v(rs-r+1,
\chi_v^r\omega_v^2)^{-N_v}
\] is holomorphic everywhere except at $s=0$ and $s=1$. Moreover there
is no pole if $\chi^r\omega^2\neq 1$.}  \quad\\

Here the factor $L_v(rs-r+1, \chi_v^r\omega_v^2)^{-N_v}$ at
each archimedean place is a kind of compensation factor, which
stems from a very subtle issue in the theory of asymptotic expansions
of matrix coefficients of real Lie groups, which will be explained in
detail in the proof of Proposition  \ref{P:poles_Eisenstein}.

Notice that by this theorem the possible poles of
$L^S(s,\pi,Sym^2\otimes\chi)$ other than at $s=0$ and $s=1$ come from
the poles of the archimedean $L$-factors $L_v(rs-r+1,
\chi_v^r\omega_v^2)^{N_v}$, which are gamma functions. Hence in
particular, we have

\quad\\
\noindent{\bf Corollary \ref{C:main}.} {\it The (incomplete) twisted
symmetric square $L$-function $L^S(s,\pi,Sym^2\otimes\chi)$ is
holomorphic everywhere in the region $\Re(s)>1-\frac{1}{2r}$ except at
$s=1$. Moreover there is no pole at $s=1$ if $\chi^r\omega^2\neq 1$.}
\quad\\

The reason we can show the holomorphy only for the region
$\Re(s)>1-\frac{1}{2r}$ is the issue at the archimedean places pointed
out above. However we believe that this can be removed and that we can prove the
following stronger version

\quad\\
\noindent{\bf Conjecture \ref{Con:main}.} {\it The (incomplete) twisted
symmetric square $L$-function $L^S(s,\pi,Sym^2\otimes\chi)$ is
holomorphic everywhere except at
$s=0$ and $s=1$. Moreover there is no pole if $\chi^r\omega^2\neq 1$.}
\quad\\

We will take up this issue in our later work (\cite{Takeda2}).
\quad\\

Let us also note that the above corollary does not tell us that the $L$-function
$L^S(s,\pi,Sym^2\otimes\chi)$ does have a pole at $s=1$ if
$\chi^r\omega^2= 1$. However, based on an observation made by Shahidi,
one can show that if $r$ is odd, then the $L$-function
$L^S(s,\pi,Sym^2\otimes\chi)$ has a pole at $s=1$ if and only if
$\check{\pi}=\pi\otimes\chi$, where $\check{\pi}$ is the
contragredient of $\pi$. (See Corollary \ref{C:main2}.)\\

Our method of proof is by Rankin-Selberg convolution with what we call
the exceptional representation of the metaplectic
double cover $\GLt_r(\A)$ of $\GL_r(\A)$, which is viewed as a natural
generalization of theta series for $r=2$. Indeed for $r=2$, the same
result has been obtained by Gelbart and Jacquet (\cite{GJ}) already in
the late 70's, whose method in turn has its origin in the work by
Shimura (\cite{Shimura}), where he considered the analogous problem in
the classical context of elliptic modular forms. Later Patterson and
Piatetski-Shapiro (\cite{PP}) generalized the method to $r=3$ though
this time they considered only the non-twisted case, \ie
$\chi=1$. Afterwards, Bump and Ginzburg (\cite{BG}) generalized the
method to arbitrary $r$ but again only for $\chi=1$. For the twisted
case, Banks worked out the case $r=3$ in \cite{Banks2}.

Bump and Ginzburg in \cite{BG} used the exceptional representation
constructed by Kazhdan and Patterson (\cite{KP}), which is a
representation of the metaplectic cover $\GLt_r$ of $\GL_r$ both
locally and globally. In order to incorporate character twist into the
work of Bump and Ginzburg, one needs to obtain the twisted version of
the exceptional representation of Kazhdan and Patterson, which we call
the twisted exceptional representation. It turns out that one needs
the twisted exceptional representation of $\GLt_{2q}$, where $2q$ is
such that $r=2q$ or $r=2q+1$.  If $q=1$, the (twisted) exceptional
representation is simply the (twisted) Weil representation of
$\GLt_2$, which is precisely what is used by Gelbart and Jacquet
(\cite{GJ}) for $r=2$ and by Banks (\cite{Banks2}) for $r=3$. For higher
ranks, one needs to construct the twisted exceptional representation. Locally
this is the Langlands quotient of an induced
representation whose inducing representation is essentially $q$ copies
of the (twisted) Weil representation of $\GLt_2$ for the local case,
and globally the residues of the Eisenstein series constructed from the
corresponding global induced representation.
This construction for the non-archimedean local field of odd
residual characteristic is carried out as a main part of the Ph.D
thesis by Banks (\cite{Banks}) supervised by Bump. And part of the
reason that Bump and Ginzburg only considered the non-twisted case is
that the twisted exceptional representation for $q>1$ was not
available at that time. 

In this paper, we first construct the twisted exceptional
representation of $\GLt_{2q}$ both for the local and global cases.
Also one needs the twisted
exceptional representation of the group $\GLtt_{2q}$, which is the
subgroup of $\GLt_{2q}$ consisting of the elements with square
determinant. This exceptional representation is
essentially a (constituent) of restriction of the exceptional
representation of $\GLt_{2q}$. One will need this only for the case
$r=2q$. After those exceptional representations are constructed, we
will prove our main theorem by computing the Rankin-Selberg
integral. For the case $r=2q$, our Rankin-Selberg integral differs
from the one by Bump and Ginzburg even for the non-twisted case. This
is to take care of the issue raised by A. Kable in his Ph.D thesis
(\cite{Kable}). Interested readers should consult his thesis,
especially the appendix, for this issue.\\

Finally, let us mention that the result of this paper will be used in
a work by Asgari and Shahidi (\cite{AS2}) for determination of the
image of the Langlands transfers from the general spin groups to
$\GL_r$ which they obtained in their earlier paper (\cite{AS}).

\quad\\

\begin{center} {\bf Notations}
\end{center} 

Throughout the paper, $F$ will be
either a local or global field of characteristic 0. If $F$ is global,
we denote the ring of adeles by $\A$.  If $F$ is a non-archimedean
local field $F$, we
denote the ring of integers by $\OF$, and the uniformizer by
$\varpi_{F}$ or simply by $\varpi$ when the field is clear from the
context. 

We fix the non-trivial additive character $\psi$ on
$F\backslash\A$ if $F$ is a number field or on $F$ if $F$ is a local
field. Though we often use the same symbol $\psi$ both for the local
and global cases, this will cause no confusion. Whether $F$ is local
or global, for each $a\in F^\times$ we denote by $\psi_a$ the
additive character defined by $\psi_a(x)=\psi(ax)$.
If $F$ is local and
$\chi$ is a character on $F^\times$, by $L(\chi)$ we mean the local Tate
factor for $\chi$. In particular for non-archimedean $F$,
$L(\chi)=(1-\chi(\varpi_F))^{-1}$ (resp. $L(\chi)=1$) if $\chi$ is
unramified (resp. ramified). If $F$ is global and
$\chi=\otimes'_v\chi_v$,  we let $L(\chi)=\prod_v L(\chi_v)$.

For the group $\GL_r$,
we often consider the two cases: $r$ is even and $r$ is odd. For the
former we let $r=2q$ and for the latter $r=2q+1$. If $P$ is a
parabolic subgroup of $\GL_r$, we denote the Levi part by $M_P$ and
the unipotent radical by $N_P$. We always assume that the Levi part
$M_P=\GL_{r_1}\times\cdots\times\GL_{r_k}$
sits in $\GL_r$ diagonally. We often denote each element
$\left(\begin{smallmatrix}g_1&&\\ &\ddots&\\
    &&g_k\end{smallmatrix}\right)\in M_P$ by $(g_1,
\dots,g_k)$ or $\diag(g_i)$ for $g_i\in\GL_{r_i}$ whenever it is
convenient. Also we denote the maximal torus of $M_P$ by $T_P$. We
denote the Borel subgroup by $B$ and we denote $T_B$ simply by $T$. 
Also we let $\delta_P$ be the modular
character of $P$. We let $W$ be the Weyl group of $\GL_r$ and we
choose each element $w\in W$ in such a way that each entry in $w$ is
either $0$ or $1$. We denote the $r\times r$ identity matrix by $I_r$.

For an algebraic group $G$ over
$F$, we sometimes write simply $G$ for the $F$-rational points, when
there is no danger of confusion. Also for a global $F$ we
sometimes denote each element in $G(\A)$ by $\prod_vg_v$ where $g_v\in
G(F_v)$.
If $A$ is a locally compact abelian group, we denote its Pontryagin dual
by $\widehat{A}$.

Let $G$ be any group and $H\subseteq G$ a subgroup. For each $g\in G$
and $h\in H$ we let $^gh=ghg^{-1}$ and $^gH=\{{^gh}: h\in H\}$. If
$\pi$ is a representation of $H$, we define the twist $^g\pi$ of $\pi$
with $g$ to be the representation of $^gH$ given by
$^g\pi(^gh)=\pi(h)$. In particular, if $H$ is normal,
$^g\pi(h)=\pi(g^{-1}hg)$ for $h\in H$. We use the symbol $\Ind$ for
normalized induction and $\ind$ for unnormalized one.

\quad\\

\begin{center} Acknowledgements
\end{center}

The author would like to thank Freydoon Shahidi for suggesting this
problem and constant help throughout. Without his help, this paper
would never have come into existence. Thanks are also due to Anthony
Kable for letting the author know about his thesis work as well as
reading an early draft and giving useful comments. The author would
like to thank Wee Teck Gan for pointing out an error in an early
version and giving him a suggestion on how to fix it. Especially, the
unfolding argument in the $r=2q$ case is based on his
suggestion. He would also like to thank Nolan Wallach for
explaining him about the theory of asymptotic expansions of matrix
coefficients. Finally, he would like to thank the anonymous referee
for carefully and patiently reading the original manuscript and giving him various
suggestions to significantly improve the readability and accuracy of
the paper. The author is partially
supported by NSF grant DMS-1215419.

\quad\\

%%%%%%%%%%%%%%%%%%%%%%%%%%%%%%%%%%%%%%%%%%%%%%%%%%%%%%%%%%%%%%%%%%%

\section{\bf The metaplectic double cover $\GLt_r$ of 
$\GL_r$}\label{S:metaplectic_cover}

%%%%%%%%%%%%%%%%%%%%%%%%%%%%%%%%%%%%%%%%%%%%%%%%%%%%%%%%%%%%%%%%%%%

In this section, we review the theory of the metaplectic double cover $\GLt_r$ of
$\GL_r$ for both local and global cases, which was originally
constructed by Kazhdan and Patterson in \cite{KP}.

%%%%%%%%%%%%%%%%%%%%%%%%%%%%%%%%%%%%%%%%%%%%%%%%%%%%%%%%%%%%%%%%%%%

\subsection{\bf The local metaplectic double cover $\GLt_r$}

%%%%%%%%%%%%%%%%%%%%%%%%%%%%%%%%%%%%%%%%%%%%%%%%%%%%%%%%%%%%%%%%%%%

Let $F$ be a (not necessarily non-archimedean) local field of
characteristic $0$. In this paper, by the metaplectic double cover
$\GLt_r(F)$ of $\GL_r(F)$, we mean the central extension of $\GL_r(F)$ by
$\{\pm1\}$ as constructed in \cite{KP} by Kazhdan and
Patterson. (Kazhdan and Patterson considered more general cover
$\GLt_r^{(c)}(F)$ with a twist by $c\in\{0,1\}$. But we only consider
the non-twisted case, \ie $c=0$.) Later, Banks, Levy, and Sepanski
(\cite{BLS}) gave an explicit description of a 2-cocycle 
\[
\sigma_r:\GL_r(F)\times\GL_r(F)\rightarrow\{\pm 1\}
\]
which defines $\GLt_r(F)$ and shows that their 2-cocycle is
``block-compatible'', by which we mean the
following property of $\sigma_r$:  For the standard
$(r_1,\dots,r_k)$-parabolic $P$ of $\GL_r$, so that its Levi $M_P$ is
of the form
$\GL_{r_1}\times\cdots\times\GL_{r_k}$ which is embedded diagonally into
$\GL_r$, we have
\begin{equation}\label{E:compatibility}
\sigma_r(\begin{pmatrix}g_1&&\\ &\ddots&\\ &&g_k\end{pmatrix},
\begin{pmatrix}g'_1&&\\ &\ddots&\\ &&g'_k\end{pmatrix})
=\prod_{i=1}^k\sigma_{r_i}(g_i,g_i')\prod_{1\leq i<j\leq k}(\det(g_i), \det(g_j'))_F,
 \end{equation}
for all $g_i, g_i'\in\GL_{r_i}(F)$ (\cite[Theorem 11, \S3]{BLS}), where
$(-,-)_F$ is the Hilbert symbol for $F$.
The 2-cocycle of \cite{BLS} generalizes the
well-known cocycle given by Kubota (\cite{Kubota}) for the case
$r=2$. Note that $\GLt_r(F)$ is not the $F$-rational points of an
algebraic group, but this notation seems to be standard. 

We need to recall how this cocycle is constructed. Let
$\G_r=\SL_{r+1}$. Matsumoto in \cite{Matsumoto} constructed the
metaplectic double cover $\Gbt_r$ of $\G_r$. A cocycle
$\sigma_{\G_r}$ defining the cover $\Gbt_r$ is explicitly computed
in \cite{BLS}, and satisfies the block-compatibility in a much
stronger sense (\cite[Theorem 7, \S2]{BLS}). Consider the embedding
\[
l:\GL_r(F)\rightarrow\G_r(F),\quad
g\mapsto \begin{pmatrix}g&\\ &\det(g)^{-1}\end{pmatrix}.
\]
Then the cocycle $\sigma_r$ is defined by
\[
\sigma_r(g, g')=\sigma_{\G_r}(l(g), l(g'))(\det(g), \det(g'))_F.
\]
(See \cite[p.146]{BLS}.)

We define $\sigGLt_r(F)$ to be the group whose underlying set is
\[ 
\sigGLt_r(F) =\GL_r(F)\times\{\pm 1\}=\{(g,\xi):g\in\GL_r(F), \xi\in\{\pm
1\}\},
\] 
and the group law is defined by
\[
(g_1,\xi_1)\cdot (g_2,\xi_2)=(g_1g_2,\sigma_r(g_1,g_2)\xi_1\xi_2).
\]
Since we would like to emphasize the cocycle being used, we write
$\sigGLt_r(F)$ instead of $\GLt_r(F)$.

To use the block-compatible 2-cocycle of \cite{BLS} has obvious
advantages. In particular, it can be explicitly computed and, of course, it is
block-compatible. However it does not allow us to construct the global
metaplectic cover $\GLt_r(\A)$. Namely one
cannot define the adelic block-combatible 2-cocycle simply by taking
the product of the local block-combatible 2-cocycles over all the
places. This can be already observed for the case $r=2$. (See
\cite[p.125]{F}.)

For this reason, we will use a different 2-cocycle $\tau_r$ which
works nicely with the global metaplectic cover $\GLt_r(\A)$. To construct such
$\tau_r$, first assume $F$ is non-archimedean. It is known that an
open compact subgroup $K$
splits in $\GLt_r(F)$, and moreover if the residue characteristic of $F$ is odd,
$K=\GL_r(\OF)$. (See \cite[Proposition 0.1.2]{KP}.) Also for
$k_1,k_2\in K$, we have $(\det(k_1),\det(k_2))_F=1$.
Hence one has a
continuous map $s_r:\GL_r(F)\rightarrow\{\pm1\}$  such that
$\sigma_r(g_1,g_2)s_r(g_1)s_r(g_2)=s_r(g_1g_2)$ for all $g_1,g_2\in
K$. Then define our 2-cocycle $\tau_r$ by
\begin{equation}\label{E:tau_sigma}
\tau_r(g_1,g_2):=\sigma_r(g_1,g_2)s_r(g_1)s_r(g_2)/s_r(g_1g_2)
\end{equation}
for $g_1, g_2\in\GL_r(F)$. If $F$ is archimedean, we set
$\tau_r=\sigma_r$. 

The choice of $s_r$ and hence $\tau_r$ is not unique. However when the
residue characteristic of $F$ is odd, there is a canonical choice with
respect to the splitting of $K$ in the following sense. 
Since the cocycle $\sigma_r$
is the restriction of $\sigma_{\G_r}$ to the image of the embedding
$l$, and it is known that the compact group $\G_r(\OF)$ also splits in
$\Gbt_r(F)$, there is a map $\sss_r:\G_r(F)\rightarrow\{\pm1\}$ such
that the section $\G_r(F)\rightarrow\Gbt_r(F)$ given by $(g,\sss_r(g))$
is a homomorphism on $\G_r(\OF)$. (Here we are assume $\Gbt_r(F)$ is
realized as $\G_r(F)\times\{\pm1\}$ as a set and the group structure
is defined by the cocycle $\sigma_{\G_r}$.) Moreover $\sss_r|_{\G_r(\OF)}$ is
determined up to twists by the elements in $H^1(\G_r(\OF), \{\pm
1\})=\Hom(\G_r(\OF), \{\pm 1\})$. But $\Hom(\G_r(\OF), \{\pm 1\})=1$
since $\G_r(\OF)$ is perfect, and
hence $\sss_r|_{\G_r(\OF)}$ is unique. (See \cite[p. 43]{KP} for this
matter.) 
We choose $s_r$ so that
\begin{equation}\label{E:canonical_section}
s_r=\sss_r|_{l(\GL_r(\OF))}.
\end{equation} 
With this choice, we have the commutative diagram
\begin{equation}\label{E:canonical_diagram}
\xymatrix{
\sigGLt_r(\OF)\ar[r]&\Gbt_r(\OF)\\
K\ar[r]^l \ar[u]^{k\mapsto (k,\;s_r(k))}&\G_r(\OF),\ar[u]_{k\mapsto(k,\;\sss_r(k))}
}
\end{equation}
where the top arrow is $(g,\xi)\mapsto (\l(g),\xi)$ and  all the
arrows can be seen to be
homomorphisms. This choice of $s_r$ is crucial for constructing the
metaplectic tensor product of automorphic representations in Appendix
\ref{S:tensor_product}. Also note that the left vertical arrow in the
above diagram is what is called the canonical lift in \cite{KP} and
denoted by $\kappa^\ast$ there.

Also when $r=2$,  we assume that $\tau_2$ is chosen to
be the cocycle $\beta$ used in \cite[p.125]{F}, which can be shown to be
block-compatible, and equal to the choice we made above when the
residue characteristic of $F$ is odd.

Using $\tau_r$, we realize $\GLt_r(F)$ to be
\[ 
\GLt_r(F)=\GL_r(F)\times\{\pm 1\},
\] 
as a set and the group law is given by
\[ 
(g_1,\xi_1)\cdot(g_2,\xi_2)=(g_1g_2, \tau_r(g_1,g_2)\xi_1\xi_2).
\] 
Note that we have the exact sequence
\[
\xymatrix{
0\ar[r]&\{\pm1\}\ar[r]&\GLt_r(F)\ar[r]^{p_r}&\GL_r(F)\ar[r]&
0
}
\]
given by the obvious maps, where we call $p_r$ the canonical projection.

We define a set theoretic section
\[ 
\kappa:\GL_r(F)\rightarrow\GLt_r(F),\; g\mapsto (g,1).
\] 
Note that $\kappa$ is not a homomorphism. But by our
construction of the
cocycle $\tau_r$, $\kappa|_K$ is a homomorphism if $F$ is
non-archimedean and $K$ is a sufficiently small open compact subgroup,
and moreover if the residue characteristic is odd, one has
$K=\GL_r(\OF)$.

Also we define another set theoretic section
\[
\s:\GL_r(F)\rightarrow\GLt_r(F),\; g\mapsto (g,s_r(g)^{-1})
\]
where $s_r(g)$ is as above. We sometimes write $\s$ for $\s_r$ when
we would like to emphasize the rank of the group. We have the
isomorphism
\[
\GLt_r(F)\rightarrow\sigGLt_r(F),\quad (g,\xi)\mapsto (g,s_r(g)\xi), 
\]
which gives rise to the commutative diagram
\[
\xymatrix{
\GLt_r(F)\ar[rr]&&\sigGLt_r(F)\\
&\GL_r(F)\ar[ul]^{\s}\ar[ur]_{g\mapsto (g,1)}&
}
\]
of set theoretic maps, \ie maps which are not necessarily homomorphisms. Also note
that the elements in the image $\s(\GL_r(F))$ ``multiply via
$\sigma_r$'' in the sense that for $g_1,g_2\in\GL_r(F)$, we have
\begin{equation}\label{E:convenient}
(g_1,s_r(g_1)^{-1}) (g_2,s_r(g_2)^{-1})=(g_1g_2, \sigma_r(g_1,g_2)s_r(g_1g_2)^{-1}).
\end{equation}

For a subgroup
$H\subseteq\GL_r(F)$, whenever the cocycle $\sigma_r$ is trivial on
$H\times H$, the section $\s$ splits 
$H$ by (\ref{E:convenient}). We often denote the image $\s(H)$ by
$H^\ast$ or sometimes simply
by $H$ when it is clear from the context. Particularly important is
that by \cite[Theorem 7 (f), \S3]{BLS}, $\s$ splits $N_B$, the unipotent radical of the Borel
subgroup $B$ of $\GL_r(F)$, and
accordingly we denote $\s(N_B)$ by $N_B^\ast$. (Note that in
\cite{BG} and \cite{KP}, the notation $H^\ast$ seems to be used whenever $H$ splits in
$\GLt_r(F)$ via \emph{any} section. But we avoid this abuse of
notation. For example, if $F$ is non-archimedean of odd residual
characteristic,
$\GL_r(\OF)$ splits via $\kappa$ but not via $\s$, and hence the
notation $\GL_r(\OF)^\ast$ does not make sense in this paper.)

Assume $F$ is non-archimedean of odd residue characteristic. By
\cite[Proposition 0.I.3]{KP} we have
\begin{equation}\label{E:kappa_and_s}
\kappa|_{T\cap K}=\s|_{T\cap K},\quad \kappa|_{W}=\s|_{W},\quad
\kappa|_{N_B\cap K}=\s|_{N_B\cap K},
\end{equation}
where $W$ is the Weyl group and $K=\GL_r(\OF)$. 
In particular, this implies $s_r|_{T\cap K}=s_r|_{W}=s_r|_{N_B\cap
  K}=1$. Also note that $s_r(1)=1$. In particular the section $\s$
splits the Weyl group $W$. If the residue characteristic of $F$ is not odd, however, $\s$ does
not split $W$. Indeed, $\s$ splits $W$ if and only if
$(-1,-1)_F=1$. (See \cite[\S 5]{BLS}.) Yet in either case,
for each element $w\in W$, we
denote $\s(w)$ simply by $w$, when it
 is clear from the context. 

Note that $\GLt_1=\GL_1(F)\times\{\pm 1\}$, where the product is the
direct product, \ie $\sigma_1$ is trivial. (See \cite[Corollary
8, \S3]{BLS}.) Also we define $\widetilde{F^\times}$ to be
$\widetilde{F^\times}=F^\times\times\{\pm 1\}$ as a set but the
product is given by $(a_1,\xi_1)\cdot(a_2,\xi_2)=(a_1a_2,
(a_1,a_2)_F\xi_1\xi_2)$. (It is known that $\widetilde{F^\times}$ is
isomorphic to $\GLt_1$ if and only if $(-1,-1)_F=1$. It is our
understanding that this is due to J. Klose (\cite[p.42]{KP}), though we
do not know where his proof is written. See
\cite{Adams} for a proof for a more general statement.)

For each subgroup $H(F)\subseteq\GL_r(F)$, we denote the preimage
$p_r^{-1}(H(F))$ of $H(F)$ via the canonical projection $p_r$ by
$\widetilde{H}(F)$  or sometimes simply by
$\widetilde{H}$ when the base field is clear from the context. We call
it the ``metaplectic preimage'' of $H(F)$. 

If $P$ is a parabolic subgroup of $\GL_r$ whose Levi is
$M_P=\GL_{r_1}\times\cdots\times\GL_{r_k}$, we often write
\[ 
\MPt=\GLt_{r_1}\timest\cdots\timest\GLt_{r_k}
\] 
for the metaplectic preimage of $M_P$. One can check
\[
\Pt=\MPt N_P^\ast
\]
because each element in
$\widetilde{N_P}$ is written in the form $(1,\xi)n^\ast$ for $n^\ast\in
N_P^\ast$ and $\xi\in\{\pm1\}$, and $(1, \xi)\in\MPt$. Moreover, one
can check by using \cite[Theorem 7 (f), \S3]{BLS} that $N_P^\ast$ is
normalized by
$\MPt$. Also we have $\MPt\cap N_P^\ast=\{(1,1)\}$. Hence if $\pi$
is a representation of $\MPt$, one can consider the induced
representation $\Ind_{\MPt N_P^\ast}^{\GLt_r}\pi$ as usual by letting
$N_P^\ast$ act trivially. This is the reason we prefer to write
$\Pt=\MPt N_P^\ast$ rather than $\Pt=\MPt\widetilde{N_P}$.

Next let
\[ 
\GL_r^{(2)}=\{g\in\GL_r:\det g\in (F^\times)^2\}, 
\] 
and $\GLtt_r$ its metaplectic preimage. Also we define
\[ 
M_P^{(2)}=\{(g_1,\dots,g_k)\in M_P: \det g_i\in (F^\times)^2\}
\] 
and often denote its preimage by
\[ 
\MPtt=\GLtt_{r_1}\timest\cdots\timest\GLtt_{r_k}.
\] 
We write $P^{(2)}=M_P^{(2)}N_P$ and denote its preimage by
$\Ptt$. Then we have
\[
\Ptt=\MPtt N_P^\ast.
\]
As explained for $\Pt$, it is preferable to write $\Ptt$ in this way for
forming induced representations.

\quad

Let us mention the following important fact. Let $Z\subseteq\GL_r$ be
the center of $\GL_r$. Then $\Zt$, which is abelian, is not the center of
$\GLt_r$ in general. And it is the center only when $r=2q+1$ or
$F=\C$. If $r=2q$ and $F\neq \C$, the preimage of
\[
    Z^{0}:=\{aI_r:a\in (F^\times)^2\}\subset\GL_r
\]
is the center of $\GLt_r$, where $I_r$ is the identity matrix. From
(\ref{E:compatibility}), one can compute
\[
\sigma_r(a_1I_r, a_2I_r)=\prod_{1\leq i<j\leq r}(a_1,a_2)_F=(a_1,a_2)_F^{\frac{1}{2}r(r-1)}.
\]
Hence for either $r=2q$
or $r=2q+1$, $\Zt$ is isomorphic to
$\widetilde{F^\times}$ if $q$ is odd, and isomorphic to $\GLt_1$ if
$q$ is even. Also note that for $r=2q$ we have $\Zt\subset\GLtt_r$ and
it is the center of $\GLtt_r$.

\quad

Let $\pi$ be an admissible representation of a subgroup
$\widetilde{H}\subseteq \GLt_r$. We say $\pi$ is ``genuine'' if each
element $(1,\xi)\in\widetilde{H}$ acts as multiplication by $\xi$, so
if $\pi$ is genuine, it does not descend to a representation of $H$ via the canonical projection
$\widetilde{H}\rightarrow H$. On the other hand, if $\pi$ is a
representation of $H$, one can always view it as a (non-genuine)
representation of $\widetilde{H}$ by pulling back $\pi$ via the
canonical projection $\widetilde{H}\rightarrow H$, which we denote by the
same symbole $\pi$. In particular, for a parabolic subgroup $P$, we
view the modular character $\delta_P$ as a character on $\Pt$ in this way.

%%%%%%%%%%%%%%%%%%%%%%%%%%%%%%%%%%%%%%%%%%%%%%%%%%%%%%%%%%%%%%%%%%%

\subsection{\bf The global metaplectic double cover $\GLt_r$}\label{S:group}

%%%%%%%%%%%%%%%%%%%%%%%%%%%%%%%%%%%%%%%%%%%%%%%%%%%%%%%%%%%%%%%%%%%

In this subsection we consider the global metaplectic group. So we let
$F$ be a number field and $\A$ the ring of adeles. We shall define the
2-fold metaplectic cover $\GLt_r(\A)$ of $\GL_r(\A)$. (Just like the
local case, we write $\GLt_r(\A)$ even though it is not the adelic
points of an algebraic group.) The construction of $\GLt_r(\A)$ has
been done in various places such as \cite{KP, FK}. 

First define the adelic 2-cocycle $\tau_r$ by
\[
    \tau_r(g_1,g_2):=\prod_v\tau_{r,v}({g_1}_v,{g_2}_v),
\]
for $g_1, g_2\in\GL_r(\A)$, where $\tau_{r,v}$ is the local
cocycle defined in the previous subsection and ${g_i}_v$ is the
$v$-component of $g_i$ as usual. By definition of $\tau_{r,v}$, we
have $\tau_{r,v}({g_1}_v,{g_2}_v)=1$ for almost all $v$, and hence the
product is well-defined. 

We define $\GLt_r(\A)$ to be the group whose underlying set is
$\GL_r(\A)\times\{\pm 1\}$ and the group structure is defined as in the
local case, \ie
\[
    (g_1,\xi_1)\cdot(g_2,\xi_2)=(g_1g_2, \tau_r(g_1,g_2)\xi_1\xi_2),
\]
for $g_i\in\GL_r(\A)$, and $\xi_i\in\{\pm 1\}$.  Just as the local case, we have
\[
\xymatrix{
0\ar[r]&\{\pm1\}\ar[r]&\GLt_r(\A)\ar[r]^{p_r}&\GL_r(\A)\ar[r]&0,
}
\]
where we call $p_r$ the canonical projection. Define a set theoretic section
$\kappa:\GL_r(\A)\rightarrow\GLt_r(\A)$ by
$g\mapsto(g,1)$. 

It is well-known that $\GL_r(F)$ splits in $\GLt_r(\A)$. However the
splitting is not via $\kappa$. In what follows, we will write the
splitting $\GL_r(F)\rightarrow \GLt_r(\A)$ explicitly.

Let us start with
\begin{Prop}\label{P:product_formula}
For $g_1, g_2\in\GL_r(F)$, we have $\sigma_{r,v}(g_1,g_2)=1$ for almost
all $v$, and further
\[
\prod_v\sigma_{r,v}(g_1,g_2)=1.
\]
\end{Prop}
\begin{proof}
From the explicit description of the cocycle $\sigma_{r,v}(g_1,g_2)$
given at the end of $\S 4$ of \cite{BLS}, one can see that
$\sigma_{r,v}(g_1,g_2)$ is written as a product of Hilbert
symbols of the form $(t_1,t_2)_{F_v}$ for $t_i\in F^\times$. This proves
the first part of the proposition. The second part follows from the
product formula for the global Hilbert symbol.
\end{proof}

This ``product formula'' of the block-compatible 2-cocycle implies
\begin{Prop}
If $g\in\GL_r(F)$, then for almost all $v$, we have $s_{r,v}(g)=1$, where
$s_{r,v}$ is the map $s_{r,v}:\GL(F_v)\rightarrow\{\pm 1\}$ defining the
section $\s:\GL(F_v)\rightarrow\GLt_r(F_v)$.
\end{Prop}
\begin{proof}
By the Bruhat decomposition we have $g=bwb'$ for some $b, b'\in B(F)$ and
$w\in W$. Then for each place $v$
\begin{align*}
s_{r,v}(g)
&=s_{r,v}(bwb')\\
&=\sigma_{r,v}(b, wb')s_{r,v}(b)s_{r,v}(wb')/\tau_{r,v}(b,wb')\quad\text{by
  (\ref{E:tau_sigma})}\\
&=\sigma_{r,v}(b, wb')s_{r,v}(b) \sigma_{r,v}(w,
b')s_{r,v}(w)s_{r,v}(b')/\tau_{r,v}(w,b')\tau_{r,v}(b,wb') \quad\text{again by
  (\ref{E:tau_sigma})}.
\end{align*}
By the previous proposition, $\sigma_{r,v}(b, wb')=\sigma_{r,v}(w,
b')=1$ for almost all $v$.
By (\ref{E:kappa_and_s}) we know $s_{r,v}(b)=s_{r,v}(w)=s_{r,v}(b')=1$ for almost
all $v$. Finally by definition of $\tau_{r,v}$,
$\tau_{r,v}(w,b')=\tau_{r,v}(b,wb')=1$ for almost all $v$.
\end{proof}

This proposition implies that the expression 
\[
s_r(g):=\prod_vs_{r,v}(g)
\]
makes sense for all $g\in\GL_r(F)$, and one can define the map
\[
\s:\GL_r(F)\rightarrow\GLt_r(\A),\quad g\mapsto (g, s_r(g)^{-1}).
\]
Moreover, this is a homomorphism because of Proposition
\ref{P:product_formula} and (\ref{E:convenient}).

Unfortunately, however, the expression $\prod_vs_{r,v}(g_v)$ does
not make sense for every $g=\prod_vg_v\in\GL_r(\A)$ because one does not know
whether $s_{r,v}(g_v)=1$ for almsot all $v$. But whenever the product
$\prod_vs_{r,v}(g_v)$ makes sense we denote the element $(g,
\prod_vs_{r,v}(g_v)^{-1})$ by $\s(g)$. This defines a partial global section
$\s:\GL_r(\A)\rightarrow\GLt_r(\A)$. For example, if $g\in B(\A)$,
$\s(g)$ is defined thanks to (\ref{E:kappa_and_s}). 
(See the last paragraph of \cite[p.150]{BG} as well.) Also we know that
$\s$ splits $N_B(\A)$ thanks to \cite[Theorem 7(f), \S3]{BLS}.

Analogously to the local case, if the partial global section $\s$ is
defined on a subgroup $H\subseteq\GL_r(\A)$ and $\s|_H$ is a
homomorphism, we denote the image $\s(H)$ by $H^\ast$ or simply by $H$
when there is no danger of confusion. This applies
to, for example, $H=\GL_r(F)$ or $N_B(\A)$. But let us emphasize that we
reserve this notation only for the subgroup split by $\s$. 

Moreover we have
\begin{Lem}\label{L:cocycle_generated}
For $g\in\GL_r(F)$ and $n\in N_B(\A)$, both $\s(gn)$ and $\s(ng)$ are
defined and moreover $\s(gn)=\s(g)\s(n)$ and $\s(ng)=\s(n)\s(g)$.
\end{Lem}
\begin{proof}
To show $\s(gn)$ is defined, it suffices to show $s_r(gn)$ is defined. 
We know both $s_r(g)$ and $s_r(n)$ are defined. Moreover for all places $v$,
we have $\sigma_{r,v}(g_v, n_v)=1$ by \cite[Theorem 7(f), \S3]{BLS}. Hence
for all $v$,
$s_{r,v}(gn_v)=s_{r,v}(g)s_{r,v}(n_v)/\tau_{r,v}(g,n_v)$. For
almost all $v$, the right hand side is $1$. Hence the global $s_r(gn)$
is defined. Also this equality shows that $\s(gn)=\s(g)\s(n)$. The
same argument works for $ng$.
\end{proof}

We define the groups like $\GLtt_r(\A)$, $\MPtt(\A)$, $\Ptt(\A)$, etc
completely analogously to the local case. Also $\widetilde{\A^\times}$ is a group whose
underlying set is $\A^\times\times\{\pm1\}$ and the group structure is given by
the global Hilbert symbol analogously to the local case. Also just like the local
case, the preimage $\Zt(\A)$ of the center $Z(\A)$ is the center of
$\GLt_r(\A)$ only if $r=2q+1$. If $r=2q$, then the center of
$\GLt_r(\A)$ is $\Zt^0(\A)$, and $\Zt(\A)$ is the center of only
$\GLtt_r(\A)$.

\quad

Let $\pi$ be a representation of $\widetilde{H}(\A)\subseteq \GLt_r(\A)$. Just like the local
case, we call $\pi$ genuine if each element $(1,\xi)\in\widetilde{H}(\A)$ acts as
multiplication by $\xi$, so it does not descend to a representation
of $H(\A)$ via the projection $\widetilde{H}(\A)\rightarrow H(\A)$. 
On the other hand, any representation of
$H(\A)$ is viewed as a representation of $\widetilde{H}(\A)$ by
pulling it back by $p_r$, which we also denote by $\pi$. In
particular, this applies to the modular character $\delta_P$ for each
parabolic $P(\A)$.

\quad

We can also describe $\GLt_r(\A)$ as a quotient of a restricted direct
product of the groups $\GLt_r(F_v)$ as follows. Consider the
restricted direct product $\prod_v'\GLt_r(F_v)$ with respect to the
groups $\kappa(K_v)=\kappa(\GL_r(\mathcal{O}_{F_v}))$ for all $v$ with
$v\nmid 2$ and $v\nmid\infty$. If we denote each element in this
restricted direct product by $\Pi_v(g_v,\xi_v)$ so that $g_v\in
K_v$ and $\xi_v=1$ for almost all $v$, we have the
surjection
\begin{equation}\label{E:surjection}
    \rho:{\prod_v}'\GLt_r(F_v)\rightarrow\GLt_r(\A),\quad
    \Pi_v(g_v,\xi_v)\mapsto (\Pi_vg_v, \Pi_v\xi_v).
\end{equation}
This is a group homomorphism by our definition of $\GLt_r(F_v)$ and
$\GLt_r(\A)$. Of course
\[
    {\prod_v}'\GLt_r(F_v)/\ker\rho\cong \GLt_r(\A),
\]
where $\ker\rho$ consists of the elements of the form
$\Pi_v(1,\xi_v)$ with $\xi_v=-1$ at an even number of $v$.

Suppose we are given a collection of irreducible admissible
representations $\pi_v$ of $\GLt_r(F_v)$ such that $\pi_v$ is
$\kappa(K_v)$-spherical for almost all $v$. Then we can form an
irreducible admissible representation of $\prod_v'\GLt_r(F_v)$ by
taking a restricted tensor product $\otimes_v'\pi_v$ as usual. Suppose
further that $\ker\rho$ acts trivially on $\otimes_v'\pi_v$, which is
always the case if each $\pi_v$ is genuine. Then it
descends to an irreducible admissible representation of $\GLt_r(\A)$,
which we denote by $\otimest'_v\pi_v$, and call it the ``metaplectic
restricted tensor product''. Let us emphasize that the space for
$\otimest'_v\pi_v$ is the same as that for
$\otimes_v'\pi_v$. Conversely, if $\pi$ is an irreducible admissible
representation of $\GLt_r(\A)$, it is written as
$\otimest'_v\pi_v$ where $\pi_v$ is an irreducible admissible
representation of $\GLt_r(F_v)$, and for almost all $v$, $\pi_v$ is
$\kappa(K_v)$-spherical. (To see it, view $\pi$ as a representation of
the restricted product $\prod_v'\GLt_r(F_v)$ by pulling it back by
$\rho$ and apply the usual
tensor product theorem for the restricted product, which gives
$\otimes_v'\pi_v$, and it descends to $\otimest_v'\pi_v$.) Note
that though the restricted tensor product
(metaplectic or not) is far from canonical, each local component
$\pi_v$ is uniquely determined up to equivalence.\\

%%%%%%%%%%%%%%%%%%%%%%%%%%%%%%%%%%%%%%%%%%%%%%%%%%%%%%%%%%%%%%%%%%%

\section{\bf The exceptional representations of
$\GLt_r$}\label{S:exceptional}

%%%%%%%%%%%%%%%%%%%%%%%%%%%%%%%%%%%%%%%%%%%%%%%%%%%%%%%%%%%%%%%%%%%

In this section, we first review the theory of the (non-twisted) exceptional
representation of $\GLt_r$ of Kazhdan-Patterson (\cite{KP}), and after
that we construct 
the twisted version of it. Throughout the section we write
\[
r=\begin{cases}2q\\2q+1\end{cases}
\]
depending on the parity of $r$.

%%%%%%%%%%%%%%%%%%%%%%%%%%%%%%%%%%%%%%%%%%%%%%%%%%%%%%%%%%%%%%%%%%%

\subsection{\bf The non-twisted exceptional representation of
$\GLt_r$}\label{S:non-twisted}

%%%%%%%%%%%%%%%%%%%%%%%%%%%%%%%%%%%%%%%%%%%%%%%%%%%%%%%%%%%%%%%%%%%

Let us consider the non-twisted exceptional
representation of $\GLt_r$ developed by Kazhdan and Patterson in
\cite{KP}. We treat both the $r=2q$ and
$2q+1$ cases at the same time. Also in this subsection, all
the groups are over the local field $F$ (non-archimedean or
archimedean) or the adeles $\A$, and most of the time we consider the local and global
case at the same time.

Roughly speaking, this exceptional representation is the
Langlands quotient of a certain induced representation of $\GLt_r$
induced from the metaplectic preimage $\Bt$ of the Borel subgroup $B$,
which we will define now.

First for the maximal torus $T\subseteq B$, we let 
\[ 
T^{\e}=\{\begin{pmatrix}t_1&&\\&\ddots&\\&&t_r\end{pmatrix}\in T:
    t_1t_2^{-1}, t_3t_4^{-1},\dots,t_{2q-1}t_{2q}^{-1}\text{ are squares}\}.
\] 
The metaplectic preimage $\Tte$ of $T^{\e}$ is a maximal
abelian subgroup of $\Tt$. Also we denote $T^{\e}N_B$ by $B^{\e}$.

To define the exceptional representation, we need to recall the notion
of the Weil index attached to each (local or global) additive
character $\psi$, which was first defined by Weil in his important
paper (\cite{Weil}). A good reference (for the local case) is
\cite[Appendix]{Rao}. First consider the local case. For the
additive character $\psi$ on $F$, the map $F\rightarrow \C^\times$
defined by $x\mapsto \psi(x^2)$ is what Weil called a character of
second degree. Weil attached to any character of second degree $f$ an
eight root of unity $\gamma(f)$, which is called the Weil index of
$f$. In particular, we denote by $\gamma(\psi)$ the Weil index of
$x\mapsto \psi(x^2)$, which we call the Weil index of $\psi$. Of
course, we can also define $\gamma(\psi_a)$ for each $a\in F$
analogously. We let
\[
\mu_\psi(a):=\frac{\gamma(\psi_a)}{\gamma(\psi)}.
\]
Various properties of $\mu_\psi$ as well as those of $\gamma(\psi)$
are reviewed in \cite[Appendix]{Rao}. In particular, one has
\begin{equation}\label{E:Weil_index}
\mu_\psi(ab)=\mu_\psi(a)\mu_\psi(b)(a,b)_F.
\end{equation}
This property implies that the map
$\widetilde{F^\times}\rightarrow \C^\times$ defined by $(a,\xi)\mapsto
\xi\mu_\psi(a)$ is a homomorphism. Let us also mention that 
\begin{equation}\label{E:Weil_index2}
\mu_{\psi_a}=\mu_{\psi_b}\text{ if and only if } a\equiv b\mod (F^\times)^2.
\end{equation}
Next assume $F$ is global and $\psi$ is an additive
character on $\A$. We define $\mu_\psi:=\prod_v\mu_{\psi_v}$. By
\cite[Proposition A.11]{Rao}, $\mu_{\psi_v}=1$ on $\mathcal{O}_{F_v}$
for almost all $v$, and hence the product is well-defined. As in the
local case $\mu_\psi$ defines a character on $\widetilde{\A^\times}$.

The non-twisted exceptional representation of $\GLt_r$ is the unique
irreducible quotient of the induced representation $\Ind_{\Tte
N_B^\ast}^{\GLt_r}\omega_\chi^\psi\otimes\delta_B^{1/4}$, where
$\omega_\chi^\psi$ is the character on $\Tte$ defined as follows: 
Let $\chi$ be a unitary character of $F^\times$ if $F$ is local, and a
unitary Hecke character of $\A^\times$ if $F$ is global. Define a
character ${\omega_\chi^{\psi}}$ on $\Tte$ by
\begin{equation}\label{E:exceptional_character}
    {\omega_{\chi}^{\psi}}((1,\xi)\s(t))=\xi\chi(\det
    t)\mu_{\psi}(t_1)\mu_\psi(t_3)\mu_\psi(t_5)\cdots \mu_\psi(t_{2q-1}).
\end{equation}
 Here even when
$F$ is global, the section $\s$ is defined on $T(\A)$ and the
expression $\s(t)$ makes sense. Note that if $t=\diag(t_i), t'=\diag(t_i')\in
T^{\e}$, then one can see from (\ref{E:compatibility}) together with
basic properties of the Hilbert symbol that
\[
    \sigma_r(t,t')=(t_1,t'_1)(t_3,t'_3)(t_5,t'_5)\cdots
    (t_{2q-1},t'_{2q-1}).
\]
Then (\ref{E:Weil_index}) implies that
${\omega_\chi^\psi}$ is indeed a character on $\Tte$.

It is shown in \cite{KP} that
\begin{Prop} 
The induced representation $\Ind_{\Tte
N_B^\ast}^{\GLt_r}\omega_\chi^{\psi}\otimes\delta_B^{1/4}$ has a unique
irreducible quotient, which we denote by $\theta_\chi^{\psi}$. For the local
case, it is the image of the intertwining integral
\[ \Ind_{\Tte N_B^\ast}^{\GLt_r}\omega_\chi^{\psi}\otimes\delta_B^{1/4}
\rightarrow \Ind_{(^{w_0}\Tte)
N_B^\ast}^{\GLt_r}\;^{w_0}(\omega_\chi^{\psi}\otimes\delta_B^{1/4}),
\] 
where $w_0$ is the longest Weyl group element. (See the notation
section for the notations for the superscript $^{w_0}$. Also note that
$w_0$ is actually $\s(w_0)$ when viewed as an element in $\GLt_r$.) For the global case,
it is generated by the residues of the Eisenstein series for this
induced space, and $\theta_\chi^{\psi}$ is a square integrable automorphic
representation of $\GLt_r(\A)$. Moreover for the global $\theta_\chi^{\psi}$,
one has the decomposition
$\theta_\chi^{\psi}=\otimest'_v\theta_{\chi_v}^{\psi_v}$. 
\end{Prop}

We call the representation $\theta_\chi^{\psi}$ the non-twisted exceptional
representation of $\GLt_r$ with the \emph{determinantal character} $\chi$.

\begin{Rmk}\label{R:independence}
Assume $F$ is local. Define
\[
\Omega_\chi^{\psi}:=\Ind_{\Tte}^{\Tt}\omega_\chi^{\psi}.
\]
This is irreducible (\cite[p.55]{KP}). Also if $r$ is even, this is
independent of $\psi$. This is because each element in $\Tte\backslash\Tt$ is
represented by $\s(a_1,\dots, a_{2q})$ with $a_i\in (F^\times)^2\backslash
F^\times$ and by direct computation one can
check that the twists of $\omega_\chi^{\psi}$ by $\s(a_1,\dots, a_{2q})$
are all distinct by using (\ref{E:Weil_index2}).

By inducing in stages, one can see that
\[
    \Ind_{\Tte N_B^\ast}^{\GLt_r}{\omega_\chi^{\psi}}\otimes\delta_B^{1/4}
    =\Ind_{\Bt}^{\GLt_r}{\Omega_\chi^{\psi}}\otimes\delta_B^{1/4},
\]
which implies $\theta_\chi^{\psi}$ is independent of $\psi$ if $r$ is
even. 
\end{Rmk}

One of the important properties of the exceptional representation is
that the constant term is again an exceptional representation, which
can be called the ``periodicity'' of Jacquet module for
the non-archimedean case and the periodicity of constant terms for the
global case. Namely, locally we have
\begin{Prop}[Local Periodicity]\label{P:local_periodicity1}
Assume $F$ is non-archimedean. Let ${(\theta_\chi^{\psi})}_{N_B}$ be he
Jacquet module of $\theta_\chi^{\psi}$ along the parabolic $\Bt$. Then 
\[
{(\theta_\chi^{\psi})}_{N_B}={^{w_0}(\Omega_\chi^{\psi})}\otimes\delta_B^{1/4}
=\Omega_\chi^{\psi}\otimes\delta_B^{1/4},
\]
where $w_0$ is the longest element in the Weyl group.
\end{Prop}
\begin{proof}
The first equality is \cite[Theorem I.2.9(e)]{KP} with the notations adjusted to
ours. The second equality follows because the metaplectic tensor
products behaved in the expected way under conjugation by a Weyl group
element as proven in \cite{Takeda1}.
\end{proof}

Globally, we have
\begin{Prop}[Global Periodicity]\label{P:global_periodicity1}
Assume $F$ is a number field. Let ${(\theta_\chi^{\psi})}_{N_B}$ be the
space generated by the constant terms of the automorphic forms in
$\theta_\chi^{\psi}$ along the Borel $\Bt(\A)$. Then as a representation
of $\Tt(\A)$, we have
\[
{(\theta_\chi^{\psi})}_{N_B}={^{w_0}(\Omega_\chi^{\psi})}\otimes\delta_B^{1/4}
=\Omega_\chi^{\psi}\otimes\delta_B^{1/4},
\]
where $w_0$ is the longest element in the Weyl group.
\end{Prop}
\begin{proof}
This is not proven in \cite{KP}. But it can be proven by using the
theory of Eisenstein series developed in \cite{MW}. We will give
the detailed argument later for the twisted
case, and one may simply mimic the argument there.
\end{proof}

Finally let us mention that (locally or globally) if $r=2q+1$, under $\theta_\chi^{\psi}$
the center $\Zt$ acts by the character
\begin{equation}\label{E:central_character1}
(1,\xi)\s(z)\mapsto \xi\chi(a)^{2q+1}\mu_\psi(a)^q,\quad
z=\begin{pmatrix}a&&\\ &\ddots&\\ &&a\end{pmatrix}\in\GL_{2q+1}.
\end{equation}
As we see in Section \ref{S:group}, if $q$ is odd,
$\Zt\cong\widetilde{F^\times}$ or $\widetilde{\A^\times}$, and hence certainly $z\mapsto
\mu_\psi(a)^q$ is a character on $\Zt$. If $q$ is even, $\Zt\cong\GLt_1$
(trivial extension) but by (\ref{E:Weil_index}) one can see that the
map $z\mapsto \mu_\psi(a)^q$ is also a character.

%%%%%%%%%%%%%%%%%%%%%%%%%%%%%%%%%%%%%%%%%%%%%%%%%%%%%%%%%%%%%%%%%%%

\subsection{\bf The Weil representation of $\GLt_2$}\label{S:Weil_rep}

%%%%%%%%%%%%%%%%%%%%%%%%%%%%%%%%%%%%%%%%%%%%%%%%%%%%%%%%%%%%%%%%%%%

To construct the twisted exceptional representation of $\GLt_r$, one
needs the Weil representation of $\GLt_2$ both for the local and
global cases. In this subsection, we review the basics of the theory
of the Weil representation of $\GLt_2$. The definitive references for
this are \cite{G} and \cite{GPS}.\\

\noindent\textbf{Local case:}

Let us consider the local case, and hence $F$ will be a (not
necessarily non-archimedean) local field of characteristic
$0$. Everything stated below without any
specific reference is found in \cite[\S 2]{GPS} for the
non-archimedean case and in \cite[\S 4]{G} for the archimedean case. Let
$S(F)$ be the space of Schwartz-Bruhat functions on $F$, \ie smooth
functions with compact support if $F$ is non-archimedean, and
functions with all the derivatives rapidly decreasing if $F$ is
archimedean. Let $\rr^\psi$ denote the representation of $\SLt_2(F)$
on $S(F)$ such that
\begin{align*} \rr^\psi(\s\begin{pmatrix}0&1\\-1&0\end{pmatrix})f(x)&=\gamma(\psi)\hat{f}(x)\\
\rr^\psi(\s\begin{pmatrix}1&b\\0&1\end{pmatrix})f(x)&=\psi(bx^2)f(x),\qquad b\in F\\
\rr^\psi(\s\begin{pmatrix}a&0\\0&a^{-1}\end{pmatrix})f(x)&=|a|^{1/2}\mu_\psi(a)f(ax),\qquad a\in
F^\times\\ \rr^\psi(1,\xi)f(x)&=\xi f(x),
\end{align*} 
where $\hat{f}(x)=\int f(y)\psi(2xy)\,dy$ with the Haar
measure $dy$ normalized in such a way that
$\hat{\hat{f}}(x)=f(-x)$. Also $\gamma(\psi)$ is the Weil
index of $\psi$, and $\mu_\psi(a)=\gamma(\psi_a)/\gamma(\psi)$. 
It is well-known that $\rr^\psi$ is reducible and
written as $\rr^\psi=\rr^\psi_+\oplus\rr^\psi_-$, where $\rr^\psi_+$
(resp. $\rr^\psi_-$) is an irreducible representation realized in the
subspace of even functions (resp. odd functions) in $S(F)$. 

Let $a\in F^\times$. For each $g\in\SLt_2(F)$ let us write
\[
g^a=\s\begin{pmatrix}1&\\ &a\end{pmatrix}^{-1}
g\; \s\begin{pmatrix}1&\\ &a\end{pmatrix}.
\]
\begin{Lem}\label{L:Weil_twist}
Let $\epsilon\in\{\pm\}$ be fixed. For each $a\in F^\times$, let
$^a\rr_\epsilon^\psi$ be the representation of $\SLt_2(F)$ defined by
$^a\rr_\epsilon^\psi(g)=\rr_\epsilon^\psi({g^a})$ for all
$g\in\SLt_2(F)$. Then
\[
^a\rr_\epsilon^\psi=\rr_\epsilon^{\psi_a}.
\]
\end{Lem}
\begin{proof}
This is \cite[Proposition 2.27]{G}.
\end{proof}

Let $\chi$ be a unitary character on $F^\times$. If $\chi(-1)=1$
(resp. $\chi(-1)=-1$), one can extend $\rr^\psi_+$
(resp. $\rr^\psi_-$) to a representation $\rr^\psi_\chi$ of
$\GLtt_2(F)$ by letting
\[ \rr^\psi_\chi(\s\begin{pmatrix}1&0\\0&a^2\end{pmatrix})f(x)=\chi(a)|a|^{-1/2}f(a^{-1}x).
\]
This is indeed a well-defined irreducible
representation of $\GLtt_2(F)$ and we call it the Weil representation of
$\GLtt_2(F)$ associated with $\chi$. We denote by $S_\chi(F)$ the subspace of $S(F)$
in which $\rr_\chi^\psi$ is realized, which is the
space of even functions if $\chi(-1)=1$ and odd functions if
$\chi(-1)=-1$. Note that
\begin{equation}\label{E:central_character}
\rr^\psi_\chi(\s\begin{pmatrix}a&0\\0&a\end{pmatrix})f(x)=\chi(a)\mu_\psi(a)f(x).
\end{equation}
Lemma \ref{L:Weil_twist} implies
\begin{Lem}\label{L:Weil_twist2}
For $a\in F^\times$, let $^a\rr_\chi^\psi$ be the representation of $\GLtt_2(F)$ obtained by
conjugating $\rr_\chi^{\psi}$ by $\s\left(\begin{smallmatrix}1&\\
  &a\end{smallmatrix}\right)$. Then
\[
^a\rr_\chi^\psi=\rr_\chi^{\psi_a}.
\]
\end{Lem}
Also note
\begin{Lem}\label{L:Weil_equivalent}
$\rr^{\psi_a}_\chi$ and $\rr^{\psi_b}_\chi$ are equivalent if and
only if $a\equiv b\mod (F^\times)^2$. 
\end{Lem}
\begin{proof}
See \cite[(1.3)]{GPS}.
\end{proof}

The Weil representation $\rr_\chi$ of $\GLt_2(F)$ is defined by
\[ 
\rr_\chi=\Ind_{\GLtt_2(F)}^{\GLt_2(F)}\rr_\chi^\psi.
\] 
Then $\rr_\chi$ is irreducible and
independent of the choice of $\psi$, and hence our notation. By Mackey
theory together with $^a\rr_\chi^\psi=\rr_\chi^{\psi_a}$, we have
\begin{equation}\label{E:restriction}
\rr_\chi|_{\GLtt_2(F)}=\bigoplus_{\alpha\in\Sigma}\rr_\chi^{\psi_\alpha},
\end{equation}
where $\Sigma$ is a set of representatives of
$(F^\times)^2\backslash F^\times$, because
$\GLtt_2(F)\backslash\GLt_2(F)=\Sigma$.

If $\chi(-1)=1$, one can check that $\rr_\chi$ is the exceptional
representation of Kazhdan-Patterson for $r=2$ with the determinantal
character $\chi^{1/2}$. (See \cite[Proposition 2.3.3]{GPS} for the
non-archimedean case, and \cite[\S 6]{GPS} for the archimedean case.)
Namely, we have the embedding
\begin{equation}\label{E:Weil_embedding}
\rr_\chi\hookrightarrow\Ind_{\Bt}^{\GLt_2}{^s(\Omega_{\chi^{1/2}}^\psi\otimes\delta_B^{1/4})}
=\Ind_{(^s\Tte)N_B^\ast}^{\GLt_2}{^s(\omega_{\chi^{1/2}}^\psi\otimes\delta_B^{1/4})},
\end{equation}
where $s$ is the Weyl group element
$s=\left(\begin{smallmatrix}&1\\1&\end{smallmatrix}\right)$.
Similarly we have the embedding
\begin{equation}\label{E:Weil_embedding2}
\rr_\chi^\psi\hookrightarrow
\Ind_{(^s\Tte)N_B^\ast}^{\GLtt_2}{^s(\omega_{\chi^{1/2}}^\psi\otimes\delta_B^{1/4})}.
\end{equation}
Let us mention that one can choose any $\chi^{1/2}$ because in general
for any quadratic
character $\epsilon$ and character $\eta$, we have
$\omega_{\epsilon\eta}^\psi=\omega_{\eta}^\psi$ for a character of
$\Tte\subseteq\GLt_r$ as long as $r$ is even.

If $\chi(-1)=-1$, then $\rr_\chi$ is supercuspidal for the
non-archimedean case (\cite[Proposition 3.3.3]{GPS}), is a discrete
series representation of lowest weight $3/2$ for the real case
(\cite[\S 6]{GPS}) and is identified with a certain induced
representation for the complex case (\cite[\S 6]{GPS}).

\begin{Prop}\label{P:Whittaker_local_Weil} 
The Weil representation
$\rr_\chi^\psi$ of $\GLtt_2(F)$ is $\psi_a$-generic if and only if
$a=b^2$. Also in this case, the $\psi_{b^2}$-Whittaker functional on
$S_\chi(F)$ is (a scalar multiple of) the functional given by $f\mapsto
f(b)$.
\end{Prop}
\begin{proof}
This seems to be folkloric, though the author does not know any reference for
it. So we will give a brief proof here. 
First of all, since $\rr_\chi^\psi$ is extended from the
representation $\rr^\psi_{\pm}$ of $\SLt_2(F)$, it suffices to show
the corresponding statement for $\rr^\psi_{\pm}$. From the explicit
description of the action of $\SLt_2(F)$, it is immediate that the
functional given by $f\mapsto f(b)$ is a $\psi_{b^2}$-Whittaker
functional. This shows one direction.

The non-obvious part is the converse. One way to prove this is to
invoke the theory of Waldspurger developed in \cite{Wald1, Wald2}, according
to which an irreducible admissible representation $\pit$ of
$\SLt_2(F)$ has a
non-zero theta lift with respect to $\psi_a$ to $\PGL_2(F)$ if and
only if $\pit$ has a $\psi_{a}$-Whittaker functional. But from the
explicit theta correspondences obtained in \cite[Theorem 1]{Wald2}, one
can see that this is possible only when $\rr_{\pm}^\psi$ is isomorphic
to $\rr_{\pm}^{\psi_a}$, which implies $a\in
(F^\times)^2$. (Apparently to use the theory of Waldspurger is
overkill and too indirect. One can directly prove it by using
a theory of distributions. But in the interest of space, we only give
this indirect proof here.)
\end{proof}

This proposition together with (\ref{E:restriction}) implies that the
Weil representation $\rr_\chi$ of $\GLt_2(F)$ is $\psi_a$-generic for
any $a$.\\

\noindent\textbf{Global case:}

Next we consider the global Weil representation. So we let $F$ be a
number field, $\A$ the ring of adeles and $\chi$ a unitary Hecke
character on $\A^\times$. We define the global Weil representation
$\rr_\chi$ of $\GLt_2(\A)$ as the restricted tensor product of the
local Weil representations, \ie
\[ 
\rr_\chi=\otimest'\rr_{\chi_v}.
\] 
It is shown in \cite[\S8]{GPS} that $\rr_\chi$ is a square
integrable automorphic representation of $\GLt_2(\A)$, and moreover it
is cuspidal if and only if $\chi^{1/2}$ does not exist. Also one can
see that if $\chi^{1/2}$ exists, then just like the local case,
$\rr_\chi$ is the exceptional representation of Kazhdan-Patterson for
$r=2$, namely $\rr_\chi=\theta_{\chi^{1/2}}$. (Again as in the local
case, it is independent of the choice of $\chi^{1/2}$.)

We also define the global Weil representation $\rr_\chi^\psi$ of
$\GLtt_2(\A)$ by
\[ 
\rr_\chi^\psi=\otimest'{\rr_{\chi_v}^{\psi_v}}.
\] 
Then  $\rr_\chi^\psi$ can be realized in the subspace
$S_\chi(\A)=\otimes'S_{\chi_v}(F_v)$ of the space $S(\A)$ of
Schwartz-Bruhat functions on $\A$ with the action given by the same
formulas as the local case.

The two representations $\rr_\chi$ and $\rr_\chi^\psi$ are related by
\begin{Prop}\label{P:Weil_restriction1}
Let $\rr_\chi^{(2)}$ be the representation of $\GLtt_2(\A)$ whose
space is $\{f|_{\GLtt_2(\A)}:f\in \rr_{\chi}\}$,
namely the space of restrictions to $\GLtt_2(\A)$ of automorphic forms
in $\rr_\chi$. Then as a representation of $\GLtt_2(\A)$, we have
\[ 
\rr_\chi^{(2)}=\bigoplus_{a\in (F^\times)^2\backslash
F^\times}\rr_\chi^{\psi_a}.
\] 
\end{Prop}
\begin{proof} 
As in the proof of \cite[Proposition 8.1.1]{GPS}, a
typical element in the space of $\rr_\chi$ is written as
$\Phi=(\Phi_a)_{a\in\Sigma}$, where the indexing set is
$\Sigma=(\A^\times)^2\backslash\A^\times$, and each $\Phi_a$ is
in $S(\A)$, on which $\GLtt_2(\A)$ acts as $\otimest'
\rr_{\chi_v}^{\psi_{a_v}}$. Then the function $\varphi_\Phi$ on
$\GLt_2(\A)$ defined by
\[ \varphi_\Phi(g)=\sum_{a\in(F^\times)^2\backslash
F^\times}\sum_{\xi\in F}(\rr_\chi(g)\Phi)_a(\xi)
\] gives an automorphic realization of $\rr_\chi$. Here note that the
natural map $(F^\times)^2\backslash
F^\times\rightarrow(\A^\times)^2\backslash\A^\times$ is an injection
by the Hasse-Minkowski theorem.

Notice that the representation $\otimest'
\rr_{\chi_v}^{\psi_{a_v}}$ is an automorphic representation of
$\GLtt_2(\A)$ if and only if $a\in(F^\times)^2\backslash
F^\times$. Then one can see that for $a\in(F^\times)^2\backslash
F^\times$, the function on $\GLtt_2(\A)$ given by $g\mapsto
\sum_{\xi\in F}(\rr_\chi(g)\Phi)_a(\xi)$ is an automorphic form
on $\GLtt_2(\A)$, which is in the space of
$\rr_\chi^{\psi_a}$. Hence $\rr_\chi^{\psi_a}$ is a
constituent of $\rr_\chi^{(2)}$. Since we know
$\rr_\chi$ is square integrable and so $\rr_\chi^{(2)}$ is in
the space of square integrable automorphic forms on $\GLtt_2(\A)$, it
is completely reducible. The proposition follows.
\end{proof}

The following is the global analogue of Proposition
\ref{P:Whittaker_local_Weil}

\begin{Prop}\label{P:Whittaker_global_Weil} 
The Weil representation
$\rr_\chi^\psi$ of $\GLtt_2(\A)$ is $\psi_a$-generic if and only if
$a=b^2$ for $b\in F^\times$.
\end{Prop}
\begin{proof} This is implied by the local case, or one may directly
compute the $\psi_a$-Whittaker coefficient for the automorphic
realization $\sum_{\xi\in F}(\rr_\chi(g)\Phi)_1(\xi)$ of
$\rr_\chi^\psi$ as in the proof of the above proposition.
\end{proof}

%%%%%%%%%%%%%%%%%%%%%%%%%%%%%%%%%%%%%%%%%%%%%%%%%%%%%%%%%%%%%%%%%%%

\subsection{\bf The Weil representation of $\MPt$}\label{S:Weil_rep2}

%%%%%%%%%%%%%%%%%%%%%%%%%%%%%%%%%%%%%%%%%%%%%%%%%%%%%%%%%%%%%%%%%%%

In this subsection, we assume $r=2q$ and $P$ is the
$(2,\dots,2)$-parabolic, so that
\[
M_P=\underbrace{\GL_2\times\cdots\times\GL_2}_{q \text{ times}}.
\]
Recall from Section \ref{S:group} that  we write
$\MPt=\GLt_2\timest\cdots\timest\GLt_2$ and
$\MPtt=\GLtt_2\timest\cdots\timest\GLtt_2$. 
Let $R=F$ if $F$ is local and $R=\A$ if $F$ is global. Then we let
\[
(M_p)^{(2)}:=M_p\cap\GL_{2q}^{(2)}=\{(g_1,\dots,g_q)\in
M_P:\prod\det(g_i)\in (R^\times)^2\}.
\]
We let $\MPttt$ be the metaplectic preimage of
$(M_p)^{(2)}$. Let us note the inclusions
\[
\MPtt\unlhd\MPttt\unlhd\MPt.
\]
Also note that $\MPtt\unlhd\MPt$. Then we have
\[
\MPtt(R)\backslash\MPt(R)
=\underbrace{(R^\times)^2\backslash R^\times\times\cdots\times
  (R^\times)^2\backslash R^\times}_{q \text{ times}}
\]
and
\[
\MPttt(R)\backslash\MPt(R)= (R^\times)^2\backslash R^\times.
\]

In this subsection, we extend the theory of the Weil
representation as discussed in the previous subsection to the groups
$\MPtt, \MPttt,$ and $\MPt$. Naively, the Weil representations of
those groups are simply the tensor products of $q$ copies of the Weil representation for
$\GLt_2$ or $\GLtt_2$. 

To construct a representation of $\Mt$ or $\MPtt$ out of
representations of $\GLt_2$, it is convenient to
consider the groups $\cMPt$ and $\cMPtt$ constructed by the
block-compatible cocycle $\tau_P$ in Appendix
\ref{S:tensor_product}. Since in this subsection we often use the results and notations
from Appendix \ref{S:tensor_product}, the reader is advised to read
Appendix \ref{S:tensor_product} before moving on.\\

\noindent\textbf{Local case:}

Let us consider the local case, so $F$ is a local field and $\chi$ is
a unitary character on $F^\times$. We would like to work with $q$ different additive
characters. For this purpose, we let
\[
\bar{a}=(a_1,\dots,a_q)\in
\underbrace{F^\times\times\cdots\times F^\times}_{\text{$q$ copies}}
\]
be a $q$-tuple of elements of $F^\times$.

For each $i\in\{1,\dots,q\}$, let $\rr_\chi^{\psi_{a_i}}$ be the Weil
representation of $\GLtt_2$. We define the Weil representation
$\pi_\chi^{\psi^{\bar{a}}}$ of $\cMPtt$ with respect to $\chi, \psi$ and $\bar{a}$
by the metaplectic tensor product
\[
\pi_\chi^{\psi^{\bar{a}}}:=\rr_\chi^{\psi_{a_1}}\otimest\cdots \otimest\rr_\chi^{\psi_{a_q}}.
\]
In particular, the space of
$\pi_\chi^{\psi^{\bar{a}}}$ is the usual tensor product of the spaces
of $\rr_\chi^{\psi_{a_i}}$, \ie $S_\chi(F)\otimes\cdots\otimes S_\chi(F)$ ($q$-times).
If $\bar{a}=(1,\dots,1)$, we simply write
$\pi_\chi^\psi$ for $\pi_\chi^{\psi^{\bar{a}}}$. 

\begin{Lem}\label{L:tensor_Weil}
Let $\bar{a}=(a_1,\dots,a_q), \bar{b}=(b_1,\dots,b_q)\in
F^\times\times\cdots\times F^\times$. Then
$\pi_\chi^{\psi^{\bar{a}}}\cong\pi_\chi^{\psi^{\bar{b}}}$ if and only
if $a_i\equiv b_i\mod (F^\times)^2$ for each $i$.
\end{Lem}
\begin{proof}
The if-part follows from Lemma \ref{L:Weil_equivalent}.
For the converse, recall from Appendix
\ref{S:tensor_product} that the metaplectic tensor product is defined
in terms of the tensor product representation
$\rr_\chi^{\psi_{a_1}}\otimes\cdots \otimes\rr_\chi^{\psi_{a_q}}$ of
the group $\GLtt_2\times\cdots\times\GLtt_2$ (direct product). But if
$\rr_\chi^{\psi_{a_1}}\otimes\cdots \otimes\rr_\chi^{\psi_{a_q}}\cong
\rr_\chi^{\psi_{a_1}}\otimes\cdots \otimes\rr_\chi^{\psi_{a_q}}$, then
$\rr_\chi^{\psi_{a_i}}\cong \rr_\chi^{\psi_{b_i}}$ for each $i$,
which implies $a_i\equiv b_i\mod (F^\times)^2$ by Lemma \ref{L:Weil_equivalent}.
\end{proof}

Let $\pi$ be a representation of $\cMPtt$. For each $m\in{\cMPt}$,
recall from the notation section that $^m\pi$ is the representation of
$\cMPtt$ twisted by $m$. The set of the elements of the form
\[
m=((\begin{pmatrix}1&\\ &a_1\end{pmatrix},\dots,\begin{pmatrix}1&\\
  &a_q\end{pmatrix}) 1)\in\cMPt,
\]
where each $a_i$ is chosen modulo $(F^\times)^2$, is a complete set of the
representatives of $\cMPtt\backslash\cMPt$. For each such $m$, we have
\[
^m\rr_\chi^\psi=\rr_\chi^{\psi_{a_1}}\otimest\cdots\otimest\rr_\chi^{\psi_{a_k}}
\]
because for each $i$ we have
$^{a_i}\rr_\chi^\psi=\rr_\chi^{\psi_{a_i}}$ by Lemma
\ref{L:Weil_twist2}. By Lemma \ref{L:tensor_Weil},
$^m\rr_\chi^\psi\cong\rr_\chi^\psi$ if
and only if $m\in\cMPtt$. Thus Mackey's irreducibility criterion is
satisfied and hence the induced representation
\[
\Pi_\chi:=\Ind_{\cMPtt}^{\cMPt}\pi_\chi^{\psi^{\bar{a}}}
\]
is irreducible. This is independent of the choice of $\bar{a}$ and
$\psi$. Indeed, this is the metaplectic tensor product of $q$ copies
of $\rr_\chi$ in the sense of \cite{Mezo}.

For our purposes, we woud like to consider the representation
\[
\varpi_\chi^{\psi^{\bar{a}}}:=\Ind_{\cMPtt}^{\cMPttt}\pi_\chi^{\psi^{\bar{a}}},
\]
where $\cMPttt$ is the subgroup of $\cMPt$ whose underlying set is
$(M_P)^{(2)}\times\{\pm 1\}$  and the group law is defined
via the block-compatible cocycle $\tau_P$ as defined in Appendix
\ref{S:tensor_product}. This induced representation is irreducible
because $\Pi_\chi$ is, but is dependent on $\psi$ and $\bar{a}$.
Also note that by inducing in states, we have
\[
\Pi_\chi=\Ind_{\cMPtt}^{\cMPt}\pi_\chi^{\psi^{\bar{a}}}
=\Ind_{\cMPttt}^{\cMPt}\varpi_\chi^{\psi^{\bar{a}}}.
\]

Now the set
\[
    \{((\begin{pmatrix}1&\\&a\end{pmatrix},\begin{pmatrix}1&\\&1\end{pmatrix}
\dots,\begin{pmatrix}1&\\&1\end{pmatrix}),1)\in
    \cMPt:a\in (F^\times)^2\backslash F^\times\}
\]
is a complete set of the representatives of $\cMPttt\backslash
\cMPt$. For an element of the form $\bar{a}=(a,1,\dots,1)$, we denote
$\varpi_\chi^{\psi^{\bar{a}}}$ simply by
$\varpi_\chi^{\psi^{a}}$. By Mackey theory
\[
    \Pi_\chi|_{(\cMPt)^{(2)}}=\bigoplus_{a\in(F^\times)^2\backslash
F^\times}\varpi_\chi^{\psi^{a}}.
\]

Also the set
\[
    \{((\begin{pmatrix}1&\\&a_1\end{pmatrix},\begin{pmatrix}1&\\&a_2\end{pmatrix}
\dots,\begin{pmatrix}1&\\&a_q\end{pmatrix}),1)\in
    (\cMPt)^{(2)}:a_i\in (F^\times)^2\backslash F^\times,\;a_1\cdots
    a_q\in (F^\times)^{2}\}
\]
is a complete set of the representatives of $\cMPtt\backslash
(\cMPt)^{(2)}$, and again by Mackey theory, one sees that
\[
    \varpi_\chi^{\psi^a}|_{\cMPtt}=\bigoplus_{(a_1,\dots,a_q)}
    \rr_\chi^{\psi_{aa_1}}\otimest\rr_\chi^{\psi_{a_2}}\otimest\cdots\otimest\rr_\chi^{\psi_{a_q}},
\]
where the sum is over the elements of the form $(a_1,\dots,a_q)\in (F^\times)^2\backslash
F^\times\times\cdots\times(F^\times)^2\backslash F^\times$ and
$a_1\cdots a_q\in (F^\times)^2$. 

Also from the above decomposition of $\varpi_\chi^{\psi^a}|_{\cMPt}$ we have
\begin{Lem}\label{L:space_of_decomposition}
The induced representation
$\varpi_{\chi}^{\psi^a}$ is realized in the space
$\bigoplus_{\bar{a}}S_\chi(F^q)$, where $\bar{a}=(a_1,\dots,a_q)$ runs through the
elements in $((F^\times)^2\backslash F^\times)^q$ with $a_1\cdots
a_q\in (F^\times)^2$ and each summand
$S_\chi(F^q)$ realizes the representation
$\rr_{\chi}^{\psi_{aa_1}}\otimest\cdots\otimest
\rr_{\chi}^{\psi_{a_q}}$.
\end{Lem}
Further by the above decomposition of $\Pi_\chi|_{\cMPttt}$, we have
\begin{Lem}\label{L:space_of_decomposition2}
The representation
$\Pi_\chi$ is realized in the space
$\bigoplus_{a\in(F^\times)^2\backslash F^\times}\bigoplus_{\bar{a}}S_\chi(F^q)$, where
$\bar{a}=(a_1,\dots,a_q)$ runs through the elements in
$((F^\times)^2\backslash F^\times)^q$ with $a_1\cdots a_q\in(F^\times)^2$ and each summand
$S_\chi(F^q)$ realizes the representation
$\rr_{\chi}^{\psi_{aa_1}}\otimest\cdots\otimest
\rr_{\chi}^{\psi_{a_q}}$.
\end{Lem}

Finally for the local case, let us mention the genericity of
$\varpi_\chi^\psi$. Recall from Proposition
\ref{P:Whittaker_local_Weil} that the Weil representation
$\rr_\chi^{\psi_a}$ is $\psi_b$-generic if and only if $b\equiv a \mod
(F^\times)^2$. Hence if we define the additive character
$\psi_{(a_1,\dots,a_q)}$ on the unipotent part $N_B\cap M_P$ of $M_P$
by
\begin{equation}\label{E:additive_a}
    \psi_{(a_1,\dots,a_q)}(n)=\psi(a_1x_1+\cdots+a_qx_q),
\end{equation}
where
\[
    n=\left(\begin{array}{ccccc}
    1&\multicolumn{1}{c|}{x_1}&&&\\
    &\multicolumn{1}{c|}{1}&&&\\ \cline{1-2}
    &&\ddots&&\\ \cline{4-5}
    &&&\multicolumn{1}{|c}{1}&\multicolumn{1}{c}{x_q}\\
    &&&\multicolumn{1}{|c}{}&\multicolumn{1}{c}{1}
    \end{array}\right),
\]
we have
\begin{Lem}
The Weil representation $\rr_\chi^{\psi_{a_1}}\otimest\cdots\otimest\rr_\chi^{\psi_{a_q}}$
is $\psi_{(b_1,\dots,b_q)}$-generic if and only if
$b_i\equiv a_i\mod (F^\times)^2$ for each $i$.
\end{Lem}
  
Then we have
\begin{Prop}\label{P:tensor_product_Weil_local}
The representation $\varpi_\chi^\psi$ is
$\psi_{(b_1,\dots,b_q)}$-generic if and only if $b_1\cdots b_q\in
(F^\times)^2$.
\end{Prop}
\begin{proof}
Assume $\varpi_\chi^\psi$ is $\psi_{(b_1,\dots,b_q)}$-generic. Then
some $\rr_\chi^{\psi_{a_1}}\otimest\cdots\otimest\rr_\chi^{\psi_{a_q}}$ in
the decomposition of $\varpi_\chi^\psi|_{\MPtt}$ is
$\psi_{(b_1,\dots,b_q)}$-generic. Hence by the above lemma, we have
$b_i\equiv a_i \mod (F^\times)^2$. But since  $a_1\cdots a_q\in
(F^\times)^2$, we also have $b_1\cdots b_q\in (F^\times)^2$.

Conversely assume $b_1\cdots b_q\in (F^\times)^2$. Then in the decomposition
of $\varpi_\chi^\psi|_{\MPtt}$, there is a constituent 
$\rr_\chi^{\psi_{b_1}}\otimest\cdots\otimest\rr_\chi^{\psi_{b_q}}$
which is $\psi_{(b_1,\dots,b_q)}$-generic. Hence $\varpi_\chi^\psi$ is
$\psi_{(b_1,\dots,b_q)}$-generic.
\end{proof}

Let us mention how the Weil representation $\Pi_\chi$ is related to
the non-twisted exceptional representation of
Kazhdan-Patterson when $\chi^{1/2}$ exists, which we fix once an for
all. Recall from (\ref{E:Weil_embedding2}) that the Weil
representation $\rr_\chi^\psi$ embeds into the induced representation
$\Ind_{\Bte_2}^{\GLt_2}{^{s}\omega}_{\chi^{1/2}}^\psi\otimes\delta_{B_2}^{-1/4}$,
where by $B_2$ we mean the Borel subgroup of $\GL_2$ with the maximal torus
$T_2$ and $\omega_{\chi^{1/2}}^\psi$ is the character on the group
$\Tte_2$ as defined in (\ref{E:exceptional_character}). Hence we have
the embedding
\[
\pi_\chi^\psi\hookrightarrow(
\Ind_{\Bte_2}^{\GLtt_2}{^{s}\omega}_{\chi^{1/2}}^\psi
\otimes\delta_{B_2}^{-1/4})\otimest\cdots\otimest
(\Ind_{\Bte_2}^{\GLtt_2}{^{s}\omega}_{\chi^{1/2}}^\psi\otimest\delta_{B_2}^{-1/4}).
\]

Now let $B_{2,\dots,2}=B_2\times\cdots\times B_2$
(resp. $B_{2,\dots,2}^{\e}=B_2^{\e}\times\cdots\times B_2^{\e}$ ) be the product of
$q$ copies of $B_2$ (resp. $B^{\e}$), and view them as subgroups of $M_P$. Also let
$\Bt_{2,\dots,2}$ (resp. $\Bte_{2,\dots,2}$) be the metaplectic
preimage of $B_{2,\dots,2}$ (resp. $B^{\e}_{2,\dots,2}$) in
$\cMPt$, so in particular we assume that their group structures are
given by the cocycle $\tau_P$. Also for the maximal
torus $T_2\subseteq B_2$, we have $T_2^{\e}=T_2^{(2)}$, and for the
maximal torus $T$ of $B_{2,\dots,2}$, which is the same as the maximal
torus of $\GL_{2q}$, we have
$T^{\e}=T_2^{(2)}\times\cdots\times T_2^{(2)}$ ($q$ times). We view both
$\Tte$ and $\Tt$ as subgroups of $\Bt_{2,\dots,2}$, so in particular
the group structures are given by the cocycle $\tau_P$. (Actually one
can check that the restriction of  the cocycle $\tau_P$ to $T$ is the
same as $\sigma_r$.)

With those notations, one can check 
\[
(\Ind_{\Bte_2}^{\GLtt_2}{^{s}\omega}_{\chi^{1/2}}^\psi
\otimes\delta_{B_2}^{-1/4})\otimest\cdots\otimest
(\Ind_{\Bte_2}^{\GLtt_2}{^{s}\omega}_{\chi^{1/2}}^\psi\otimest\delta_{B_2}^{-1/4})
=\Ind_{\Bt_{2,\dots,2}^{\e}}^{\cMPtt}{^{w_1}\bar{\omega}}_{\chi^{1/2}}^\psi
\otimes\delta_{B_2\times\cdots\times B_2}^{-1/4},
\]
where $\bar{\omega}_{\chi^{1/2}}^\psi$ is the
character on $\Tte$ associated with $\chi^{1/2}$ as defined by
(\ref{E:exceptional_character})
and $w_1$ is the Weyl group element of the form
\[
w_1=\begin{pmatrix}s&&\\ &\ddots&\\ &&s\end{pmatrix}
\]
where $s=\left(\begin{smallmatrix}&1\\1&\end{smallmatrix}\right)$.

To sum up, we have the embedding
\begin{equation}\label{E:Weil_embedding3}
\pi_\chi^\psi\hookrightarrow \Ind_{\Bt_{2,\dots,2}^{\e}}^{\cMPtt}{^{w_1}\omega}_{\chi^{1/2}}
\otimes\delta_{B_2\times\cdots\times B_2}^{-1/4},
\end{equation}
and, by inducing both sides to $\cMPt$, we have the embedding
\begin{equation}\label{E:Weil_embedding4}
\Pi_\chi\hookrightarrow \Ind_{\Bt_{2,\dots,2}^{\e}}^{\cMPt}{^{w_1}\omega}_{\chi^{1/2}}
\otimes\delta_{B_2\times\cdots\times B_2}^{-1/4}.
\end{equation}

\quad\\

\noindent\textbf{Global case:}

Next we consider the global case, so $F$ is global, $\chi$ is a
unitary Hecke character and $\psi$ is our fixed additive character on
$F\backslash\A$. As we did in the local case, it is convenient to
consider the groups $\cMPt(\A), \cMPtt(\A)$ and $\cMPttt(\A)$ instead
of $\MPt(\A), \MPtt(\A)$ and $\MPttt(\A)$ for the sake of constructing
metaplectic tensor products. (See Appendix \ref{S:tensor_product}.)

As we did in the local case we would like to consider the tensor
product of the Weil representations with respect to possibly
different additive characters. Namely, we let
\[
\bar{a}=(a_1,\dots,a_q)\in \underbrace{F^\times\times\cdots\times
  F^\times}_{q \text{ times}},
\]
and we define
\[
    \pi_\chi^{\psi^{\bar{a}}}=\rr_\chi^{\psi_{a_1}}\otimest\cdots\otimest\rr_\chi^{\psi_{a_q}}
\]
to be the metaplectic tensor product representation of $\cMPtt(\A)$, where
$\rr_\chi^{\psi_{a_i}}$ is the global Weil representation of
$\GLtt(\A)$ with respect to the additive character $\psi_{a_i}$. As in
Appendix \ref{S:tensor_product}, the space of this metaplectic tensor
product is the same as that of the usual tensor product
$\rr_\chi^{\psi_{a_1}}\otimest\cdots\otimest\rr_\chi^{\psi_{a_q}}$,
and moreover since each $\rr_\chi^{\psi_{a_i}}$ is automorphic, so is
the metaplectic tensor product by Proposition
\ref{P:tensor_product_automorphy}. 
Note that
\[
    \pi_\chi^{\psi^{\bar{a}}}=\otimest_v'\pi_{\chi_v}^{\psi^{\bar{a}}_v},
\]
where at each $v$ we view $\bar{a}=(a_1,\dots,a_q)$ naturally as in
$((F_v^\times)^2\backslash F_v^\times)^q$.

Next we let
\[
    \varpi_\chi^{\psi^{\bar{a}}}=\otimest_v'\varpi_{\chi_v}^{\psi^{\bar{a}}_v}.
\]
To see its automorphy, recall from Section \ref{S:Weil_rep} that each
$\rr_{\chi_v}^{\psi_v}$ is realized in the subspace
$S_{\chi_v}(F_v)$ and accordingly $\rr_\chi^\psi$ is realized in a subspace
$S_\chi(\A)$ of the space of Schwartz functions on $\A$. Hence the
representation $\pi_\chi^{\psi^{\bar{a}}}$ is realized in a subspace
$S_\chi(\A^q) =S_\chi(\A)\otimes\cdots\otimes S_\chi(\A)$
of the space of Schwartz functions on $\A^q$. (Once again, the space
of the metaplectic tensor product is the same as that of the usual
tensor product.) Now let
\[
    \Sigma_v=\{\bar{a}=(a_1,\dots,a_q)\in((F_v^\times)^2\backslash F_v^\times)^{q}:a_1\cdots a_q\in(F_v^\times)^2\}.
\]
From Lemma \ref{L:space_of_decomposition} the representation
$\varpi_{\chi_v}^{\psi_v}$ is realized in the space
$\bigoplus_{\bar{a}\in\Sigma_v}S_{\chi_v}({F_v}^q)$, where each
$S_{\chi_v}({F_v}^q)$ realizes the representation
$\rr_{\chi_v}^{\psi_{v\;a_1}}\otimest\cdots\otimest
\rr_{\chi_v}^{\psi_{v\;a_q}}$ for each
$(a_1,\dots,a_q)\in\Sigma_v$. Then as we have seen for the Weil
representation of $\GLtt_2$
in Section \ref{S:Weil_rep}, the global representation
$\varpi_\chi^{\psi^a}$ is realized in the space of elements of the
form $\Phi=(\Phi_{\bar{a}})_{\bar{a}\in\Sigma_\A}$, where the indexing
set $\Sigma_\A$ is given by
\[
    \Sigma_\A=\{\bar{a}=(a_1,\dots,a_q)\in((\A^\times)^2\backslash
    \A^\times)^{q}:a_1\cdots a_q\in(\A^\times)^2\}.
\]
Now the representation $\varpi_\chi^{\psi^a}$ has an automorphic
realization similarly to the Weil representation of $\GLt_2(\A)$. Namely for each element
$\Phi=(\Phi_{\bar{a}})$, we put
\[
    \varphi_\Phi(g)=\sum_{\bar{a}\in \Sigma_F}\sum_{\xi\in
      F}(\varpi_\chi^{\psi^a}(g)\Phi_{\bar{a}})(\xi),
\]
where $g\in (\MPt(\A))^{(2)}$ and
\[
    \Sigma_F=\{\bar{a}=(a_1,\dots,a_q)\in((F^\times)^2\backslash
    F^\times)^{q}:a_1\cdots a_q\in(F^\times)^2\}.
\]
Then as in \cite[Proposition 8.1.1]{GPS}, one sees that the map
$\Phi\mapsto\varphi_\Phi$ defines an embedding of
$\varpi_\chi^{\psi^a}$ into the space of automorphic forms on
$(\MPt(\A))^{(2)}$.

Once we obtain this automorphic realization of
$\varpi_\chi^{\psi^a}$, the following global analogue of Proposition
\ref{P:tensor_product_Weil_local} follows just as Proposition
\ref{P:Whittaker_global_Weil}.
\begin{Prop}\label{P:tensor_product_Weil_global}
Let $(b_1,\dots,b_q)\in((F^\times)^2\backslash F^\times)^q$. Then
$\varpi_\chi^\psi$ is $\psi_{(b_1,\dots,b_q)}$-generic if and only if
$b_1\cdots b_q\in (F^\times)^2$, where the additive character
$\psi_{(b_1,\dots,b_q)}$ is defined analogously to the local case.
\end{Prop}
\begin{proof}
One can prove it in the same way as Proposition
\ref{P:Whittaker_global_Weil}. 
\end{proof}

Essentially this says that many of the
Whittaker-Fourier coefficients for the forms in $\varpi_\chi^\psi$  vanish. 
This proposition will play a crucial role in our computation for
unfolding of our Rankin-Selberg integral for the case $r=2q$.

Finally, we define the global Weil representation $\Pi_\chi$ of
$\cMPt(\A)$ by
\[
\Pi_\chi:=\otimest'\Pi_{\chi_v},
\]
where each $\Pi_{\chi_v}$ is the local Weil representation of
$\cMPt(F_v)$ as defined previously. One can prove the automorphy of
$\Pi_\chi$ in the same way as the automorphy of
$\varpi_\chi^\psi$. (Let us mention that this is precisely the
metaplectic tensor product of $q$ copies of the Weil representation in
the sense defined in \cite{Takeda1}.)

Analogously to Proposition \ref{P:Weil_restriction1}, we have
\begin{Prop}\label{P:Weil_restriction2}
Let $\Pi_\chi^{(2)}$ be the representation of $\cMPttt(\A)$ whose
space is $\{f|_{\cMPttt(\A)}:f\in\Pi_\chi\}$, namely the space of
restrictions to $\cMPttt(\A)$ of automorphic forms in $\Pi_\chi$. Then
as a representation of $\cMPttt(\A)$, we have
\[
\Pi_\chi^{(2)}=\bigoplus_{a\in(F^\times)^2\backslash F^\times}\varpi_\chi^{\psi^a}.
\]
\end{Prop}
\begin{proof}
The proof is essentially identical to Proposition \ref{P:Weil_restriction1}.
\end{proof}

%%%%%%%%%%%%%%%%%%%%%%%%%%%%%%%%%%%%%%%%%%%%%%%%%%%%%%%%%%%%%%%%%%%

\subsection{\bf The twisted exceptional representation of $\GLt_{2q}$}

%%%%%%%%%%%%%%%%%%%%%%%%%%%%%%%%%%%%%%%%%%%%%%%%%%%%%%%%%%%%%%%%%%%

We construct the twisted exceptional representation of $\GLt_r$
when $r=2q$ for both the local and global cases. But for a
non-archimedean local field of odd residual characteristic, this is
one of the main achievements of the Ph.D thesis by Banks (\cite{Banks}). The
basic idea for the local case is that just like the non-twisted
exceptional representation of Kazhdan-Patterson the twisted one is
constructed as a quotient of the induced representation $\Ind_{\MPt
  N^\ast_P}^{\GLt_r}\Pi_\chi\otimes\delta_P^{1/4}$, where $\Pi_\chi$
is the Weil representation of $\MPt$ constructed in the previous subsection and extended
trivially on $N^\ast_P$. For
this purpose, Banks explicitly computed the local coefficients for
intertwining operators on this induced representation and showed that
it has a unique irreducible quotient, which is
the image of an intertwining operator. This quotient is precisely the
twisted exceptional representation. But for technical reasons,
Banks treated only the case of odd residual characteristic. 

However, thanks to the recent work by Ban and Jantzen (\cite{BJ}) that
proves the Langlands quotient theorem for metaplectic covers over the
$p$-adic field, the construction of the twisted exceptional
representation for the non-archimedean case is very simple. (But let us mention
that the approach taken by \cite{Banks} gives more
information about the induced representation such as the point of
reducibility.) Also let us note that over the archimedean field the
Langlands quotient theorem has been already available for groups like
$\GLt_{2q}$ (See \cite[Chapter IV]{BW}. Note that the groups
$\GLt_{2q}(\R)$ and $\GLt_{2q}(\C)$ are real reductive groups in the
sense of \cite[0.3.1]{BW}, which are also called real reductive groups
in the Harrish-Chandra class in \cite[p.289]{Wallach} to which the
general theory of \cite{BW} applies.) Hence the construction of the twisted exceptional
representation for the archimedean case is very simple as
well. Indeed, Kazhdan-Patterson constructed the
non-twisted exceptional representation over the archimedean field by
the Langlands quotient theorem as well. The
twisted case can be treated in the same way. The global case is a
standard argument in the Langlands theory of Eisenstein series, which
is also the method employed by Kazhdan-Patterson for the non-twisted
case.

\quad

Throughout this subsection, $r=2q$ and $P$ is the
$(2,\dots,2)$-parabolic whose Levi $M_P$ is
$\GL_2\times\cdots\times\GL_2$
($q$-times). Also we need to view the
group $\cMPt$ as the subgroup $\MPt$ of $\GLt_r$ via the embedding
$\varphit_P:\cMPt\rightarrow \GLt_r$ (See
Appendix \ref{S:tensor_product}.) In other words when we treat the
group $\MPt$ by itself, we always
mean $\cMPt$ and when we would like to view it as a subgroup of
$\GLt_r$ we consider it as the image of the embedding $\varphit_P$. 
Accordingly, we view the Weil
representation $\Pi_\chi$ constructed in the previous subsection as a
representation of $\MPt$ via $\varphit_P$, namely
$\Pi_\chi\circ\varphit_P^{-1}$. But we simply write $\Pi_\chi$ for
$\Pi_\chi\circ\varphit_P^{-1}$ since this does not produce any
confusion. The same applies to the representations $\varpi_\chi^\psi$
and $\pi_\chi^\psi$.

Let us set up some general notations. For each standard parabolic subgroup $Q$ of
$\GL_{2q}$, we let $T_Q$ be the maximal
torus in $Q$, and $\Phi_Q$ the set of roots of $\GL_{r}$ relative to
$T_Q$. The choice of $Q$ determines the positive roots in $\Phi_Q$. We have
the natural inclusion $\Phi_Q(\C)\subset\Phi_B(\C)$ via the inclusion
$M_B\hookrightarrow M_Q$. We write
$\rho_Q$ for half the sum of positive roots in $\Phi_Q$. Assume the
Levi part $M_Q$ of $Q$ is of the form $\GL_{r_1}\times\cdots\times\GL_{r_k}$
Each root
$\beta\in\Phi_Q$ is identified with a pair of integers $\beta=(i,j)$
for $1\leq i,j\leq k$ with $i\neq j$, and $\beta=(i,j)$ is
positive if $i<j$. To be precise, let us denote each element in $M_Q$
by $\diag(g_l)$ where each $g_l\in\GL_{r_l}$ for $l=1,\dots,k$. Then
for $\beta=(i,j)$, we have $\beta(\diag(g_l))=\det(g_i)\det(g_j)^{-1}$.

Now assume $Q=P$, \ie the $(2,\dots,2)$-parabolic subgroup of
$\GL_{2q}$. Define $W_P$ to be the set of block matrices
\[
    W_P=\{(\delta_{wi,j}I_2):w\in S_q\},
\]
where $S_q$ is the symmetric group on $q$ letters and $\delta_{i,j}$
is the Kronecker delta function. Then $W_P$, which is isomorphic to
$S_q$, is a subgroup of the Weyl group $W_B$ of $\GL_{2q}$. For each
elements $\diag(h_k)\in M_P$ and $w\in W_P$, we have
$w\diag(h_k)w^{-1}=\diag(h_{wk})$.  We often view each element $w\in
W_P$ as the element $\s(w)$ in $\GLt_{2q}$. For each root
$\beta=(i,j)\in\Phi_P$, we let $\beta^\vee$
be the corresponding coroot, so that for each
$t\in F^\times$ we have $\beta^\vee(t)=\diag(g_l)$ where $g_l=I_{r_l}$ for
$l\neq i,j$, $g_i=tI_2$ and $g_j=t^{-1}I_2$, \ie
\[
    \beta^\vee(t)=\begin{pmatrix}\ddots&&&&
\\&tI_2&&&\\&&\ddots&&\\&&&t^{-1}I_2&\\&&&&\ddots\end{pmatrix},
\]
where $tI_2$ and $t^{-1}I_2$ are in the $i^\text{th}$ and
$j^\text{th}$ entries respectively and all the other diagonal entries are
$I_2$. The space $\Phi_P(\C):=\Phi_P\otimes_{\Z}\C$ is identified with
$\C^{q-1}$ by choosing a basis to be the set of the simple roots \ie the roots of
the form $(i, i+1)$. For each $\nu\in\Phi_P(\C)$ and a representation
$\Pi$ of $\MPt$ (locally or globally), we define the representation
\[
    {\Pi^{\nu}}:=\Pi\otimes\exp(\nu,H_P(\;))
\]
of $\Pt$ where $\Pi$ is extended trivially to the unipotent
part, and $H_P$ is the Harish-Chandra map as usual (or
strictly speaking the Harish-Chandra map composed with the canonical
projection $\GLt_{2q}\rightarrow\GL_{2q}$). If
$\nu=\rho_P/2\in\Phi_P(\C)$, then
$\Pi^\nu=\Pi\otimes\delta_P^{1/4}$, where $\delta_P$ is the modular
character of $P$. We often write
\[
\delta_P^\nu:=\exp(\nu,H_P(\;)).
\]
Note that $\delta_P$ can be computed as
\begin{equation}\label{E:modular}
\delta_P(g_1,\dots,g_q)=|\det(g_1)|^{2(q-1)}|\det(g_2)|^{2(q-3)}\cdots |\det(g_q)|^{-2(q-1)}
\end{equation}
for the element $(g_1,\dots,g_q)\in M_P$.

By following the notation of \cite[p.62]{KP}, for each irreducible representation
$\Pi$ of $\MPt(\A)$ (resp. $\MPt(F)$) if $F$ is global (resp. local) and for each root
$\beta\in\Phi(\C)$, we define the character on $\A^\times$
(resp. $F^\times$) by 
\begin{equation}\label{E:Pi_character1}
(\Pi)_\beta(t)=\Pi((\beta^\vee(t^2),1)
\end{equation}
for $t\in\A^\times$ (resp. $\in F^\times$). Note that the map
$t\mapsto(\beta^\vee(t^2), 1)$ is indeed a homomorphism from
$\A^\times$ to $\MPt(\A)$ (resp. $F^\times$ to $\MPt(F)$). 
Namely this is just the central
character of $\Pi$ evaluated at $(\beta^\vee(t^2),1)$. In particular, by
considering
$\pi_\chi^\psi=\rr_\chi^\psi\otimest\cdots\otimest\rr_\chi^\psi$ and
the central character of $\rr_\chi^\psi$ is given by
(\ref{E:central_character}), one can see that
\begin{equation}\label{E:Pi_character2}
(\Pi_\chi^\nu)_\beta(t)=\delta_P^\nu(\beta^\vee(t^2)),
\end{equation}
and if $\nu=\rho_P/2$, so $\delta_P^\nu=\delta_P^{1/4}$, and
$\beta=(i,j)$, then by using (\ref{E:modular})
\begin{equation}\label{E:Pi_character3}
\delta^{1/4}_P(\beta^\vee(t^2))=|t|^{4(j-i)}.
\end{equation}

Let us recall the notion of intertwining integrals. First assume
$F$ is local and $\Pi$ is an irreducible admissible representation of
$\MPt$. For $w_1,w_2\in W_P$, we define the intertwining integral
\begin{equation}\label{E:local_intertwining}
    A({^{w_1}\nu}, {^{w_1}\Pi},w_2):\Ind_{{\Pt}}^{\GLt_{2q}}{^{w_1}}
(\Pi^\nu)\rightarrow
\Ind_{{\Pt}}^{\GLt_{2q}}{^{w_2w_1}}(\Pi^\nu)
\end{equation}
by
\[
    A({^{w_1}\nu}, {^{w_1}\Pi},w_2)f(g)=\int_{N_P^\ast}f({w_2}^{-1}ng)\;dn
\]
for
$f\in\Ind_{\Pt}^{\GLt_{2q}}{^{w_1}}(\Pi_\chi^{\nu})$.

Next assume $F$ is global and $\Pi$ an irreducible automorphic
representation of $\MPt(\A)$. By following \cite{MW} we view the
induced representation $\Ind_{\Pt(\A)}^{\GLt_{2q}}\Pi$ as a space of
automorphic forms on $N_P(\A)^\ast M_P(F)^\ast\backslash\GLt_{2q}(\A)$. Then the
global intertwining integral $M({^{w_1}\nu}, {^{w_1}\Pi},w_2)$ is
defined in the completely analogous way as the local case.

\quad

\noindent\textbf{Local case:}

Let us consider the local case. But as we mentioned at the beginning
of this subsection, the construction of the twisted exceptional
representation is quite simple thanks to the Langlands quotient
theorem. But first we should mention
\begin{Lem}\label{L:Ind_stages}
Let $\pi$ be an irreducible admissible representation of $\MPtt$ such
that $\Pi:=\Ind_{\MPtt}^{\MPt}\pi$ is irreducible, so
$\varpi:=\Ind_{\MPtt}^{\MPttt}\pi$ is irreducible as well. Then
\[
\Ind_{\MPt N_P^\ast}^{\GLt_r}\Pi=\Ind_{\MPttt N_P^\ast}^{\GLt_r}\varpi
=\Ind_{\MPtt N_P^\ast}^{\GLt_r}\pi.
\]
\end{Lem}
\begin{proof}
The proof is straightforward. See \cite[Proposition 4.1]{Banks} as well.
\end{proof}

By this lemma, together with the fact that $\MPtt$ is better
behaved in the sense that each $\GLtt_2$-factor in the Levi $\MPtt$
commutes with each other, it is easier to work with
$\Ind_{\MPtt N_P^\ast}^{\GLt_r}\pi$ than $\Ind_{\MPt N_P^\ast}^{\GLt_r}\Pi$.

With this said, the local twisted exceptional representation is
constructed as follows:
\begin{Prop} 
The induced representation
$\Ind_{\Ptt}^{\GLt_r}\pi_\chi^\psi\otimes\delta_P^{1/4}$ has a unique
irreducible quotient, which we denote by $\vartheta_\chi$. It is
the image of the intertwining integral
\[ 
\Ind_{\Ptt}^{\GLt_r}\pi_\chi^\psi\otimes\delta_P^{1/4}\rightarrow
\Ind_{\Ptt}^{\GLt_r}\;^{w_0}(\pi_\chi^\psi\otimes\delta_P^{1/4}),
\] 
where $w_0$ is the longest Weyl group element relative to $P$. (Recall
from Section \ref{S:group} that $\Ptt=\MPtt N_P^\ast$.)
\end{Prop}
\begin{proof}
Let us first note that if $\chi^{1/2}$ exists, from the embedding
(\ref{E:Weil_embedding3}), one can see that the situation boils down
to the non-twisted case of Kazhdan-Patterson. Hence we assume
that $\chi^{1/2}$ does not exist. 

Let us consider the non-archimedean case.
 As we noted in Section \ref{S:Weil_rep} the Weil
representation $\rr_\chi^\psi$ is supercuspidal, and hence in
particular tempered. Thus $\pi_\chi^\psi$ and so $\Pi_\chi$ are tempered.
Then the Langlands quotient theorem for
metaplectic covers (\cite[Theorem 4.1]{BJ}) applies to this situation,
and implies that the induced representation
$\Ind_{\Pt}^{\GLt_r}\Pi_\chi\otimes\delta_P^{1/4}=\Ind_{\Ptt}^{\GLt_r}
\pi_\chi^\psi\otimes\delta_P^{1/4}$
has a unique
irreducible quotient, which is also obtained as a unique irreducible
subrepresentation of
$\Ind_{\Ptt}^{\GLt_r}\;^{w_0}(\pi_\chi^\psi\otimes\delta_P^{1/4})$. One
needs to show that this irreducible quotient
is indeed obtained as the image of the intertwining
integral. (Unlike the usual Langlands quotient theorem, it is not
shown in \cite{BJ} that the Langlands quotient is indeed obtained as
the image of the intertwining integral.) But this can be easily proven
for the case at hand
because in exactly in the same way as the proof of \cite[Proposition
I.2.2]{KP} one can show that 
\[
      \dim\Hom_{\GLt_{2q}}(\Ind_{\Ptt}^{\GLt_r}\pi_\chi^\psi\otimes\delta_P^{1/4},\;
      \Ind_{\Ptt}^{\GLt_r}\;^{w_0}(\pi_\chi^\psi\otimes\delta_P^{1/4}))\leq 1
\]
by standard computations of the Jacquet modules of the induced
representations. (Also see \cite[Corollay 6.7]{Banks}.)

The archimedean case follows from the Langlands quotient theorem,
which is available for the groups like $\GLt_r$ as we mentioned at the
beginning of this subsection.
\end{proof}

We call the representation $\vartheta_\chi$ the ``twisted exceptional
representation'' of $\GLt_{2q}$. By Lemma \ref{L:Ind_stages},
$\vartheta_\chi$ is also the quotient of
$\Ind_{\Pt}^{\GLt_r}\Pi_\chi\otimes\delta_P^{1/4}$. Since $\Pi_\chi$ is independent of the
choice of $\psi$, so is $\vartheta_\chi$, and hence our notation. If
$\chi^{1/2}$ exists, then the twisted and non-twisted ones are related
as $\vartheta_\chi=\theta_{\chi^{1/2}}^\psi$. (As we mentioned in Remark
\ref{R:independence}, if $r$ is even, $\theta_{\chi^{1/2}}^\psi$ is
independent of the choice of $\psi$.)

Finally, we have the analogue of Proposition \ref{P:local_periodicity1}.
\begin{Prop}[Local Periodicity]\label{P:local_periodicity2}
Assume $F$ is non-archimedean. Let ${(\vartheta_\chi)}_{N_P}$ be he
Jacquet module of $\vartheta_\chi^\psi$ along the parabolic $\Pt$. Then 
\[
{(\vartheta_\chi)}_{N_p}={^{w_0}(\Pi_\chi)}\otimes\delta_P^{1/4}
=\Pi_\chi\otimes\delta_P^{1/4}.
\]
\end{Prop}
\begin{proof}
This proof is completely analogous to the proof the non-twisted case
(\cite[Theorem I.2.9(e)]{KP}) and left to the reader.\\
\end{proof}

\noindent\textbf{Global case:}

We construct the global twisted exceptional representation of
$\GLt_r(\A)$, so $F$ is a number field, $\chi$ is a unitary Hecke
character and $\psi$ is our fixed additive character on
$F\backslash\A$. The construction is analogous to the local case in
that the exceptional representation is obtained as a unique
irreducible quotient of the global induced space
$\Ind_{\Pt(\A)}^{\GLt_r(\A)}\Pi_\chi\otimes\delta_P^{1/4}$, where $\Pi_\chi$
is the global Weil representation of $\MPt(\A)$. (Strictly speaking $\Pi_\chi$
is the pullback of the Weil representation of $\cMPt(\A)$ via the map
$\varphit_P^{-1}:\MPt(\A)\rightarrow\cMPt(\A)$, which is also
automorphic. See Corollary \ref{C:diagram}.) Moreover the exceptional
representation is generated by the residues of certain Eisenstein
series to be defined below.

Let us start with the definition of the Eisenstein series. Although
this might be already quite familiar to experts, let us repeat some of
the essential points of the theory of Eisenstein series. The
best reference (probably the only one for metaplectic groups) for the theory of Eisenstein
series is \cite{MW}. For a (cuspidal) automorphic representation $\Pi$ of
$\MPt(\A)$, the induced representation
$\Ind_{\Pt(\A)}^{\GLt_{2q}(\A)}{\Pi^{\nu}}$ is realized in a space of
automorphic forms on $N_P(\A)M_P(F)\backslash\GLt_{2q}(\A)$. (Here of
course we are viewing $N_P(\A)$ and $M_P(F)$ as subgroups of $\GLt_{2q}(\A)$
via the splitting $\s$ and writing simply $N_P(\A)$ and $M_P(F)$ for
$N_P(\A)^\ast$ and $M_P(F)^\ast$, respectively.) To be
precise, we have the Iwasawa decomposition 
\[
\GLt_{2q}(\A)=N(\A)\MPt(\A)\widetilde{K},
\]
where $K\subseteq\GL_{2q}(\A)$ is the usual maximal compact subgroup of
$\GL_{2q}(\A)$, namely $K=\prod' K_v$, where $K_v$ is $\GL_{2q}({\OF}_v)$ at
non-archimedean $v$, $\OO(2q)$ for real $v$ and
$\operatorname{U}(2q)$ for complex $v$. Then
\[
N_P(\A)M_P(F)\backslash\GLt_{2q}(\A)
=(M_P(F)\backslash\MPt(\A))\cdot\widetilde{K}.
\]
Hence for each automorphic form $\phi$ on
$N_P(\A)M_P(F)\backslash\GLt_{2q}(\A)$ and each $k\in\widetilde{K}$,
the function $\phi_k$ on $M_P(F)\backslash\MPt(\A) $ defined by $\phi_k(m)=
\phi(mk)$ is an automorphic form on $\MPt(\A)$.  Each $f^\nu\in
\Ind_{\Pt(\A)}^{\GLt_{2q}(\A)}{\Pi^{\nu}}$ is of the form 
\begin{equation}\label{E:section}
f^\nu=\phi\otimes\exp(\nu+\rho_P, H_P(\;)),
\end{equation}
where $\phi:N_P(\A)M_P(F)\backslash\GLt_{2q}(\A)\rightarrow \C$ is
such that for each $k\in \widetilde{K}$, the function $\phi_k$ is in
the space of $\Pi$. Also note that our induction is normalized
so that we have the shift by $\rho_P$. Also note 
\begin{equation}\label{E:section_restriction}
f^\nu|_{\MPt(\A)}\in\Pi^{\nu+\rho_P},
\end{equation}
\ie the restriction of $f^\nu$ to $\MPt(\A)$ is in the space of $\Pi^{\nu+\rho_P}$.
For each $f^\nu$, we define
the Eisenstein series by
\[
    E(g,{\Pi},f^\nu)=
    \sum_{\gamma\in P(F)\backslash\GL_{2q}(F)}f^\nu(\gamma g),
\]
where $g\in \GLt_{2q}(\A)$. It converges absolutely when $\nu$ is in
a sufficiently positive part of the Weyl chamber, and admits
meromorphic continuation. This is an automorphic form on
$\GLt_{2q}(\A)$ whenever it is holomorphic. If the inducing
representation $\Pi$ is cuspidal and $\nu$ is in the positive chamber,
the poles of $E(g,{\Pi},f^\nu)$ are
at most simple, and when it has a
(simple) pole, the residue is an automorphic form on $\GLt_{2q}(\A)$
and the space generated by the residues is a space of a square
integrable automorphic representation of $\GLt_{2q}$. 

The twisted exceptional representation to be constructed is
generated by the residues of the Eisenstein series
$E(g,{\Pi_\chi},f^\nu)$ associated with the induced representation
$\Ind_{\Pt(\A)}^{\GL_r(\A)}\Pi_\chi^\nu$ at
$\nu=\rho_P/2$. To see it, one needs to know this Eisenstein series
indeed has  a simple pole at $\nu=\rho_P/2$. But to
study poles of the Eisenstein series, one should look at the global
intertwining operator
\[
M(\nu, \Pi_\chi, w):\Ind_{\Pt(\A)}^{\GLt_r(\A)}\Pi_\chi\otimes\delta_P^{\nu}\rightarrow
\Ind_{\Pt(\A)}^{\GLt_r(\A)}\;^{w}(\Pi_\chi\otimes\delta_P^{\nu}),
\]
where $w\in W_P$. (See \cite[Proposition IV.1.11]{MW}.) We will show
that the global intertwining operator $M(\nu, \Pi_\chi, w)$ (and hence
the Eisenstein series $E(g,{\Pi_\chi},f^\nu)$) has a pole at
$\nu=\rho_P/2$ if and only if $w=w_0$. 

The computation of poles of the global intertwining operator
essentially boils down the computation of the ``normalizing factor'' for
the corresponding local intertwining operator 
\[
A(\nu, {\Pi_{\chi}}_v, w):\Ind_{\Pt(F_v)}^{\GLt_r(F_v)}{\Pi_{\chi}}_v\otimes\delta_P^{\nu}\rightarrow
\Ind_{\Pt(F_v)}^{\GLt_r(F_v)}\;^{w}({\Pi_{\chi}}_v\otimes\delta_P^{\nu})
\]
at the unramified place $v$. Namely
\begin{Lem}\label{L:spherical_section}
Assume $v$ is a place where all the data defining
$\Ind_{\Pt(F_v)}^{\GLt_r(F_v)}{\Pi_{\chi}}_v\otimes\delta_P^{\nu}$ are
unramified. Let $f_v^\nu\in
\Ind_{\Pt(F_v)}^{\GLt_r(F_v)}{\Pi_{\chi}}_v\otimes\delta_P^{\nu}$ be
such that $f_v^\nu(1)=1$. Then
\[
    A(\nu,{\Pi_\chi}_v,w_0)f_v^\nu(1)=\prod_{\substack{\beta>0\\w\beta<0}}
    \frac{L(|\;|_v^{-1}((\Pi_{\chi_v}^{\nu})_\beta)^{1/2})
     \cdot L(((\Pi_{\chi_v}^{\nu})_\beta)^{1/2})}
    {L(|\;|_v((\Pi_{\chi_v}^{\nu})_\beta)^{1/2})
     \cdot L(|\;|_v^2((\Pi_{\chi_v}^{\nu})_\beta)^{1/2})},
\]
where recall the character $(\Pi_{\chi_v}^{\nu})_\beta$ is defined in
(\ref{E:Pi_character1}), and $L$ is the local Tate factor as defined
in the notation section.
\end{Lem}
\begin{proof}
Since all the data are unramified, we have $\chi_v(-1)=1$. Then
the embedding (\ref{E:Weil_embedding4}) gives us the commutative diagram
\[
\xymatrix{
\Ind_{\Bte_{2,\dots,2}(F_v)}^{\GLt_r(F_v)}\,
^{w_1}\omega_{\chi^{1/2}}^{\psi\;\nu}
\ar[rrr]^{A({^{w_1}\nu},{^{w_1}\omega_{\chi^{1/2}}^\psi},w)}
&&&
\Ind_{\Bte_{2,\dots,2}(F_v)}^{\GLt_r(F_v)}\, ^{ww_1}(\omega_{\chi^{1/2}}^{\psi\;\nu})\\
\Ind_{\Pt(F_v)}^{\GLt_r(F_v)}{{\Pi_\chi^\nu}_v}\ar[rrr]^{A(\nu,{\Pi_\chi^\nu}_v,w)}\ar@{^{(}->}[u]
&&&
\Ind_{\Pt(F_v)}^{\GLt_r(F_v)}\, ^{w}({\Pi_\chi^\nu}_v)\ar@{^{(}->}[u]
}
\]
where $w_1$ is the Weyl group element in $W_P$ of the form
$(\left(\begin{smallmatrix}&1\\1&\end{smallmatrix}\right),
\dots,\left(\begin{smallmatrix}&1\\1&\end{smallmatrix}\right))$
and the top arrow is the intertwining operator for the
corresponding induced representations, which is studied by Kazhdan and
Patterson in \cite{KP}. Hence
\[
A(\nu,{\Pi_\chi}_v,w){f_v}^\nu(1)
=A({^{w_1}\nu},{^{w_1}\omega_{\chi^{1/2}}^\psi},w){f_v}^\nu(1).
\]
But the right hand side is computed by Kazhdan and Patterson in 
\cite[Proposition I.2.4]{KP}. Then one can see that this formula by
Kazhdan-Patterson is rewritten as in the lemma.
\end{proof}

\begin{Rmk}
The inverse of the product of the Tate factors appearing in the above
lemma can be used as a normalizing factor of the corresponding intertwining
operator.  A similar expression can be obtained for all places
$v$, which give more refined results on the twisted exceptional
representation. The author has carried out this computation, which will appear elsewhere.
\end{Rmk}

Next we need
\begin{Lem}\label{L:local_holomorphy}
Just for the sake of this lemma, let us assume $F$ is local. Then the local intertwining
operator
\[
A(\nu, {\Pi_{\chi}}, w):\Ind_{\Pt(F)}^{\GLt_r(F)}{\Pi_{\chi}}\otimes\delta_P^{\nu}\rightarrow
\Ind_{\Pt(F)}^{\GLt_r(F)}\;^{w}({\Pi_{\chi}}\otimes\delta_P^{\nu})
\]
is holomorphic for $\nu$ in the positive chamber, in particular at
$\nu=\rho_P/2$.
\end{Lem}
\begin{proof}
This is a general fact for the region of the convergence for
intertwining integrals when the inducing representation is tempered,
at least for non-metaplectic groups. For archimedean $F$, it is indeed
known (\cite[Lemma 4.2, p.84]{BW}) even for the metaplectic case.
(A proof for the $p$-adic non-metaplectic case is also available in
\cite[Proposition 2.6, p.217]{BW}.) But at this moment, it is not known
for the non-archimedean metaplectic case, or at least to the best of
our knowledge, a proof is not written anywhere, though the author has been
informed by D. Ban that this might be included in \cite{BJ} for a
future revision. But here, to be absolutely rigorous, we will give an
alternate indirect proof, which works at least for the case of our
interest. 

The idea is to use a global argument. Namely one can always
globalize the character $\chi$ to a Hecke character $\hat{\chi}$ in such a way that
at two places $v_1\neq v_2$,
$\hat{\chi}_{v_1}=\hat{\chi}_{v_2}=\chi$ and there is at least one
place $v$ at which $\hat{\chi}_v(-1)=-1$, so that $\Pi_{\hat{\chi}}$ is
cuspidal. Now if the local
intertwining operator $A(\nu, {\Pi_{\chi}}, w)$ has a pole at positive
$\nu$, the
global operator $M(\nu, {\Pi_{\hat{\chi}}}, w)=\otimes_v' A(\nu,
{\Pi_{\hat{\chi}_v}}, w)$ must have at least a double pole. But in the
positive chamber the global intertwining operator has at most a simple
pole by \cite[Proposition IV.1.11]{MW}. Hence $A(\nu, {\Pi_{\chi}},
w)$ cannot have a pole.
\end{proof}

\begin{Prop}\label{P:pole_global_intertwining}
The global intertwining operator $M(\nu, \Pi_\chi, w)$ has a pole at
$\nu=\rho_P/2$ if and only if $w=w_0$.
\end{Prop}
\begin{proof}
Assume $M(\nu, \Pi_\chi, w)$ has a pole at $\nu=\rho_P/2$. 
Let $f^\nu=\otimes' f_v^\nu\in
\Ind_{\Pt(\A)}^{\GLt_{2q}(\A)}\Pi_\chi^\nu$, where for
almost all $v$ $f_v^\nu(1)=1$. By Lemma \ref{L:spherical_section},  we can write
\begin{align*}
M(\nu,{\Pi_\chi},w)f^\nu=&\prod_{\substack{\beta>0\\w\beta<0}}
\left(\frac{L(\|\;\|_\A^{-1}(({\Pi_\chi^\nu})_\beta)^{1/2})\cdot
L((({\Pi_\chi^\nu})_\beta)^{1/2})}
{L(\|\;\|_\A(({\Pi_\chi^\nu})_\beta)^{1/2})\cdot
L(\|\;\|_\A^2(({\Pi_\chi^\nu})_\beta)^{1/2})}\right)\\
&\qquad\underset{v}{\otimes'}\prod_{\substack{\beta>0\\ w\beta<0}}
\frac{L(|\;|_v(({\Pi_\chi^\nu}_v)_\beta)^{1/2})\cdot
L(|\;|_v^2(({\Pi_\chi^\nu}_v)_\beta)^{1/2})}
{L(|\;|_v^{-1}(({\Pi_\chi^\nu}_v)_\beta)^{1/2})\cdot
L((({\Pi_\chi^\nu}_v)_\beta)^{1/2})}A(\nu,{\Pi_\chi}_v,w)f_v^\nu.
\end{align*}
Note that Lemma \ref{L:spherical_section} guarantees that the
restricted tensor product in the right hand side makes sense. 

At $\nu=\rho_P/2$, by (\ref{E:Pi_character3}) the local part in the
above decomposition is written as
\[
\frac{L(|\;|_v|\;|_v^{2(j-i)})\cdot
L(|\;|_v^2|\;|_v^{2(j-i)})}
{L(|\;|_v^{-1}|\;|_v^{2(j-i)})\cdot
L(|\;|_v^{2(j-i)})}A(\nu,{\Pi_\chi}_v,w_0)f_v^\nu
\]
for each $\beta=(i,j)$. By Lemma \ref{L:local_holomorphy} together
with the fact that all the local Tate factors appearing here have no pole,
we conclude that, if  $M(\nu,{\Pi_\chi},w)f^\nu$ has a pole at
$\nu=\rho_P/2$, it comes from the global factor. But for each
$\beta=(i,j)$, (\ref{E:Pi_character3}) gives
\begin{equation}\label{E:simple_pole}
\frac{L(\|\;\|_\A^{-1}(({\Pi_\chi^\nu})_\beta)^{1/2})\cdot
L((({\Pi_\chi^\nu})_\beta)^{1/2})}
{L(\|\;\|_\A(({\Pi_\chi^\nu})_\beta)^{1/2})\cdot
L(\|\;\|_\A^2(({\Pi_\chi^\nu})_\beta)^{1/2})}
=
\frac{L(\|\;\|_\A^{-1}\|\;\|_\A^{2(j-i)})\cdot
L(\|\;\|_\A^{2(j-i)})}
{L(\|\;\|_\A\|\;\|_\A^{2(j-i)})\cdot
L(\|\;\|_\A^2\|\;\|_\A^{2(j-i)})},
\end{equation}
at $\nu=\rho_P/2$. This has a pole if $j=i+1$, \ie $\beta$ is
simple. Now in order for the product over all $\beta>0$ with
$w\beta<0$ to have a pole at $\nu=\rho_P/2$, it must be the case that
(\ref{E:simple_pole}) has a pole for all simple $\beta>0$. This is
because by a pole, we mean a pole of a meromorphic function for
$(q-1)$ variables on $\Phi(\C)=\C^{q-1}$, which are indexed by the
simple roots, and hence to have a pole at $\nu=\rho_P^{1/2}$,
it much have a pole
at all simple $\beta>0$.
But if $w\neq w_0$, there is a simple $\beta>0$ such that
$w\beta>0$. Hence we must have $w=w_0$.

Conversely assume $w=w_0$. By reversing the reasoning, one can see
that  $M(\nu, \Pi_\chi, w_0)$ has a pole at $\nu=\rho_P/2$. (Let us note
that one can always choose the local $f_v^\nu$ so that $A(\nu,
\Pi_{\chi_v}, w_0)f_v^\nu\neq 0$, and hence the pole of
(\ref{E:simple_pole}) is not cancelled by the local factors.)
\end{proof}

Now we are ready to construct the twisted exceptional representation
as follows:
\begin{Thm}\label{T:global_twisted_exceptional_even}
At $\nu=\rho_P/2\in\Phi_P(\C)$, the global induced space
$\Ind_{\Pt(\A)}^{\GLt_{2q}(\A)}\Pi_\chi^\nu$ has a unique
irreducible quotient. Moreover, this irreducible quotient is
(equivalent to) a square integrable automorphic representation
realized in the space of automorphic forms on $\GLt_{2q}(\A)$, which are generated by
the residues of the Eisenstein for series $E(-,\Pi_\chi,f^\nu)$
at $\nu=\rho_P/2$. Let us denote this irreducible quotient by
$\vartheta_\chi$. Then
\[
    \vartheta_\chi=\otimest'\vartheta_{\chi_v},
\]
where $\vartheta_{\chi_v}$ is the local twisted exceptional
representation. We call $\vartheta_\chi$ the global twisted
exceptional representation.
\end{Thm}
\begin{proof}
First of all, let us note that if $\chi^{1/2}$ exists, then just as in
the local case this theorem is subsumed under the Kazhdan-Patterson
construction of the global exceptional representation which is
discussed in Part II of \cite{KP}. Also in this case, one can see that
our $\vartheta_\chi$ is the exceptional representation of
Kazhdan-Patterson with the determinantal character $\chi^{1/2}$, \ie
$\vartheta_\chi=\theta_{\chi^{1/2}}^\psi$. The way to reduce this case to
\cite{KP} is completely analogous to the local case, and the
detail is left to the reader.

Hence we consider the case where $\chi^{1/2}$ does not exist, and so
by \cite[Proposition 8.1.1]{GPS} we know that $\Pi_\chi$ is
cuspidal. However, even for this case the argument is essentially the
same as \cite{KP}, which is a reworking of the Langlands theory of
Eisenstein series for metaplectic groups. Of course, thanks to
\cite{MW}, this theory has been completely worked out. 

By (\cite[Proposition IV.1.11]{MW}), our Eisenstein series $E(g,{\pi_\chi^\psi},f^\nu)$ has
meromorphic continuation and has a pole when the global intertwining
operator $M(\nu,\Pi_\chi,w_0)$ has a pole, which is simple. By the
above proposition, this happens at $\nu=\rho_P/2$. Then the residues make up
the residual spectrum. Note that our inducing representation is
cuspidal. Thus the residues are square integrable automorphic
forms. (See \cite[Theorem (iii), V.3.13]{MW})

Let us write
\[
    E_{-1}(-, {\Pi_\chi},f)=
     \underset{\nu=\rho_P/2}{\Res}E(-, {\Pi_\chi},f^\nu).
\]
The map $f^\nu\mapsto  E_{-1}(-, {\pi_\chi^\psi},f)$ defines a
$\GLt_{2q}(\A)$ intertwining operator
\[
    E_{-1}:\Ind_{\Pt(\A)}^{\GLt_{2q}(\A)}({\Pi_\chi^{\rho_P/2}})\rightarrow 
    \mathcal{A}_2(\GLt_{2q}(\A)),
\]
where $\mathcal{A}_2(\GLt_{2q}(\A))$ is the space of square integrable
automorphic forms on $\GLt_{2q}(\A)$. Also let
\[
    M_{-1}({\Pi_\chi},w_0,f)=\underset{\nu=\rho_P/2}{\Res}
   M(\nu,{\Pi_\chi},w_0)f^\nu
\]
be the residue of the intertwining operator. The map
$f^\nu\mapsto  M_{-1}({\Pi_\chi},w_0,f)$ defines a
$\GLt_{2q}(\A)$ intertwining operator
\[
    M_{-1}:\Ind_{\Pt(\A)}^{\GLt_{2q}(\A)}(\Pi_\chi^{\rho_P/2})
   \rightarrow\Ind_{\Pt(\A)}^{\GLt_{2q}(\A)}(\,^{w_0}(\Pi_\chi^{\rho_P/2})).
\]
That the global induced space $\Ind_{\Pt(\A)}^{\GLt_{2q}(\A)}(\Pi_\chi^{\rho_P/2}) $ has
a unique irreducible quotient follows from the corresponding statement
for the local induced representations. The image of $M_{-1}$ is the
unique irreducible quotient, which we denote by
$\vartheta_\chi$. By decomposing the intertwining operator into local
constituents, we see
$\vartheta_\chi\cong\otimest'\vartheta_{\chi_v}$. 

It remains to show that $\vartheta_\chi$ is generated by the residues
of the Eisenstein series. This follows from the inner product
formula of pseudo-Eisenstein series (\cite[Theorem II.2.1]{MW}), which implies
(up to a suitable normalization of inner products)
that
\[
    \langle E_{-1}(-,{\Pi_\chi},f_1), E_{-1}(-,{\Pi_\chi},f_2)\rangle
    =\langle f_1, M_{-1}({\Pi_\chi},w_0,f_2)\rangle,
\]
where the inner product on the left hand side is the usual inner product on
$\mathcal{A}_2(\GLt_{2q}(\A))$ and the one on the right hand side is the usual pairing
on $\Ind_{\Pt(\A)}^{\GLt_{2q}(\A)}({\Pi_\chi^{\rho_P/2}})\times
\Ind_{\Pt(\A)}^{\GLt_{2q}(\A)}({^{w_0}({\Pi_\chi^{\rho_P/2}}}))$. (For
the derivation of this inner product, see the proof of \cite[Theorem
II.1.4]{KP}, which is based on the argument by Langlands in \cite{Langlands}.)
This inner product formula implies
\[
    \ker M_{-1}\subseteq \ker E_{-1},
\]
and so the image of $E_{-1}$ is equivalent to a quotient of the image
of $M_{-1}$ \ie
$\vartheta_\chi$. But since $\vartheta_\chi$ is the unique irreducible
quotient, the image of $E_{-1}$ is indeed isomorphic to
$\vartheta_\chi$. This completes the proof of the theorem.
\end{proof}

Finally in this section, let us give a proof of the global
periodicity of $\vartheta_\chi$, which is the twisted analogue of
Proposition \ref{P:global_periodicity1}
\begin{Prop}[Global Periodicity]\label{P:global_periodicity2}
Let ${(\vartheta_\chi)}_{N_P}$ be the space generated by the constant terms of the
automorphic forms in $\vartheta_\chi$ along the parabolic
$\Pt$. Then 
\[
{(\vartheta_\chi)}_{N_P}={^{w_0}(\Pi_\chi)\otimes\delta_P^{1/4}}
=\Pi_\chi\otimes\delta_P^{1/4}.
\] 
\end{Prop}
\begin{proof}
This follows from Proposition \ref{P:pole_global_intertwining} and a
well-known computation of contant terms of Eisenstein series. Namely by
\cite[Proposition II.1.7]{MW}, the constant term $E_P(-, \Pi_\chi,
f^\nu)$ can be computed as
\[
E_P(-, \Pi_\chi, f^\nu)=\sum_{w\in W_P}M(\nu, \Pi_\chi, w)f^\nu(-),
\]
where both sides are viewed as automorphic forms on
$\MPt(\A)$. Proposition \ref{P:pole_global_intertwining} implies
\[
\underset{\nu=\rho_P/2}{\Res}E_P(-, \Pi_\chi, f^\nu)
=\underset{\nu=\rho_P/2}{\Res}M(\nu, \Pi_\chi, w_0)f^\nu(-).
\]
But the space generated by $\left(\underset{\nu=\rho_P/2}{\Res}M(\nu,
\Pi_\chi, w_0)f^\nu\right)\Big|_{\MPt(\A)}$ is
equal to ${^{w_0}(\Pi_\chi^{\rho_P/2}})^{\rho_P}$ because the residue
$\underset{\nu=\rho_P/2}{\Res}M(\nu, \Pi_\chi, w_0)f^\nu$ is in the space of
$\Ind_{\MPt(\A)}^{\GLt_{2q}}{^{w_0}(\Pi_\chi^{\rho_P/2})}$. (Recall how
$f^\nu$ is defined in (\ref{E:section}) as well as
(\ref{E:section_restriction})). Finally the residue of the
constant term of the Eisenstein series is the same as the constant
term of the residue of the Eisenstein series. (To see this, note that
the constant term is obtained by an integral over the compact set 
$N_P(F)\backslash N_P(\A)$,
and the residue is obtained by an integral over a
closed path around the singularity, which is agan an integration over
a compact set. Two integrations over compact sets can be interchanged.)
This completes the proof.
\end{proof}

%%%%%%%%%%%%%%%%%%%%%%%%%%%%%%%%%%%%%%%%%%%%%%%%%%%%%%%%%%%%%%%%%%%

\subsection{\bf The exceptional representation of
  $\GLtt_{2q}$}\label{S:restriction}

%%%%%%%%%%%%%%%%%%%%%%%%%%%%%%%%%%%%%%%%%%%%%%%%%%%%%%%%%%%%%%%%%%%

For our purposes, especially for taking care of the issue raised by
Kable in his thesis (\cite{Kable}) for the case $r=2q$, we need to
construct the exceptional representation of $\GLtt_{2q}$ both for the
local and global cases. This exceptional representation is simply a
constituent of the restriction of the twisted exceptional
representation $\vartheta_\chi$ of $\GLt_{2q}$ constructed
in the previous subsection. The important property of those
representations (especially the global one) is the vanishing of
many of the Whittaker-Fourier coefficients, which is essentially a
generalization of the analogous property of the Weil representations
of $\GLtt_2$ as stated in  Proposition \ref{P:Whittaker_local_Weil} and
\ref{P:Whittaker_global_Weil}. The use of this property of the
exceptional representations is one of the key points for our unfolding
argument for the Rankin-Selberg integral to be considered in the next section.\\

\noindent\textbf{Local case:}

For the local case, let us first note the following.
We have the intertwining operator
\[
A(\rho_P/2, \Pi_\chi,
w_0):\Ind_{\MPt N_P^\ast}^{\GLt_{2q}}\Pi_\chi\otimes\delta_P^{1/4}\rightarrow
\Ind_{{\MPt N_P^\ast}}^{\GLt_{2q}}{^{w_0}(\Pi_\chi\otimes\delta_P^{1/4})}
\]
given by the integral as in (\ref{E:local_intertwining}). Also we have the identity
\[
    (\Ind_{\MPt N_P^\ast}^{\GLt_{2q}}\Pi_\chi\otimes\delta_P^{1/4})|_{\GLtt_{2q}}
    =\Ind_{(\MPt)^{(2)}N_P^\ast}^{\GLtt_{2q}}(\Pi_\chi\otimes\delta_P^{1/4})|_{(\MPt)^{(2)}}
    =\bigoplus_{a\in\Sigma}\Ind_{(\MPt)^{(2)}N_P^\ast}^{\GLtt_{2q}}\varpi_\chi^{\psi^a}\otimes
\delta_P^{1/4},
\]
where for the first space the restriction of representation actually
coincides with the restriction of functions in the induced space, and
$\Sigma=(F^\times)^2\backslash F^\times$. Hence by composing $A(\rho_P/2, \Pi_\chi,
w_0)$ with restriction to $\GLtt_{2q}$, we obtain a $\GLtt_{2q}$
intertwining map
\[
\bigoplus_{a\in\Sigma}\Ind_{(\MPt)^{(2)}N_P^\ast}^{\GLtt_{2q}}\varpi_\chi^{\psi^a}\otimes
\delta_P^{1/4}\rightarrow
\bigoplus_{a\in\Sigma}\Ind_{(\MPt)^{(2)}N_P^\ast}^{\GLtt_{2q}}{^{w_0}(\varpi_\chi^{\psi^a}\otimes
\delta_P^{1/4})}.
\]
Moreover one can check by direct computation that the
``$a$ component'' for each $a\in\Sigma$ on the left hand
side maps into the component for the same $a$. Namely we have
\[
A(\rho_P/2, \Pi_\chi, w_0)
=\bigoplus_{a\in\Sigma}A(\rho_P/2, \varpi_\chi^{\psi^a}, w_0).
\]

With this said, one can prove

\begin{Prop}\label{P:local_twisted_exceptional_even2}
For the local twisted exceptional representation $\vartheta_\chi$ of $\GLt_{2q}$, we have
\[
    \vartheta_\chi|_{\GLtt_{2q}}=\bigoplus_{a\in\Sigma}\vartheta_\chi^{\psi^a},
\]
where $\Sigma=(F^\times)^2\backslash F^\times$ and
$\vartheta_\chi^{\psi^a}$ is a unique irreducible quotient of
$\Ind_{(\MPt)^{(2)}N_P}^{\GLtt_{2q}}\varpi_\chi^{\psi^a}\otimes\delta_P^{1/4}$. Moreover
$\vartheta_\chi^{\psi^a}$ is the image of the intertwining integral
$A(\rho_P/2, \varpi_\chi^{\psi^a}, w_0)$.
\end{Prop}
\begin{proof}
We have the commutative diagram
\[
\xymatrix{
 (\Ind_{\MPt  N_P^\ast}^{\GLt_{2q}}\Pi_\chi\otimes\delta_P^{1/4})|_{\GLtt_{2q}}
\ar@{=}[r]\ar[d]_{A(\rho_P/2, \Pi_\chi, w_0)}&
\underset{a\in\Sigma}{\bigoplus}\Ind_{(\MPt)^{(2)}N_P^\ast}^{\GLtt_{2q}}\varpi_\chi^{\psi^a}\otimes
\delta_P^{1/4}\ar[d]^{\underset{a\in\Sigma}{\bigoplus} A(\rho_P/2, \varpi_\chi^{\psi^a}, w_0)}\\
\vartheta_\chi|_{\GLtt_{2q}}\ar@{^{(}->}[r]&\underset{a\in\Sigma}{\bigoplus}
\Ind_{(\MPt)^{(2)}N_P^\ast}^{\GLtt_{2q}}{^{w_0}(\varpi_\chi^{\psi^a}\otimes\delta_P^{1/4})}.
}
\]
Hence each
irreducible constituent of $\vartheta_\chi|_{\GLtt_{2q}}$ is an irreducible quotient of
$\Ind_{(\MPt)^{(2)}N_P^\ast}^{\GLtt_{2q}}\varpi_\chi^{\psi^a}\otimes\delta_P^{1/4}$
for some $a\in\Sigma$, which is the image of the intertwining integral
$A(\rho_P/2, \varpi_\chi^{\psi^a}, w_0)$. Moreover, one and only one
of the constituents 
of $\vartheta_\chi|_{\GLtt_{2q}}$ appears as a quotient of each
$\Ind_{(\MPt)^{(2)}N_P^\ast}^{\GLtt_{2q}}\varpi_\chi^{\psi^a}\otimes\delta_P^{1/4}$
because the number of irreducible constituents of
$\vartheta_\chi|_{\GLtt_{2q}}$ is at most the size of $\Sigma$. This
shows that for each $a\in\Sigma$, there is a unique subrepresentation $\sigma^a$ of
$\Ind_{(\MPt)^{(2)}N_P^\ast}^{\GLtt_{2q}}{^{w_0}(\varpi_\chi^{\psi^a}\otimes\delta_P^{1/4})}$
that is obtained as the image of the intertwining integral
$A(\rho_P/2, \varpi_\chi^{\psi^a}, w_0)$ such that
$\vartheta_\chi|_{\GLtt_{2q}}=\underset{a\in\Sigma}{\bigoplus}\sigma^a$.

To show the uniqueness, assume there exists $a_0\in\Sigma$ such that
$\Ind_{(\MPt)^{(2)}N_P^\ast}^{\GLtt_{2q}}\varpi_\chi^{\psi^{a_0}}\otimes\delta_P^{1/4}$
has more than two quotients, say $\sigma_1$ and $\sigma_2$. For each
$i=1,2$ and $a\in\Sigma$, let ${^a\sigma_i}$ be the representation of
$\GLtt_{2q}$ obtained by twisting $\sigma_i$ by
$\s(\left(\begin{smallmatrix}1&\\ &a\end{smallmatrix}\right), I_2,\cdots,I_2)$. Then the
representation $\underset{a\in\Sigma}{\bigoplus}{^a\sigma_i}$ extends
to a representation of $\GLt_{2q}$ which can be seen as a quotient of
$\Ind_{\MPt N_P^\ast}^{\GLt_{2q}}\Pi_\chi\otimes\delta_P^{1/4}$. 
But this induced representation has a unique irreducible quotient,
namely $\vartheta_\chi$. Hence $\sigma_1=\sigma_2$.
\end{proof}

We call the representation $\vartheta_\chi^\psi$ constructed above
``the exceptional representation'' of $\GLtt_{2q}$.\\

 This exceptional representation also has the periodicity property.

\begin{Prop}[Local Periodicity]\label{P:local_periodicity3}
Assume $F$ is non-archimedean. Let ${(\vartheta_\chi^\psi)}_{N_P}$ be he
Jacquet module of $\vartheta_\chi$ along the parabolic $\Ptt$. Then 
\[
{(\vartheta_\chi^\psi)}_{N_p}={^{w_0}(\varpi_\chi^\psi)}\otimes\delta_P^{1/4}
=\varpi_\chi^\psi\otimes\delta_P^{1/4}.
\]
\end{Prop}
\begin{proof}
This follows from the above lemma together with the periodicity of
$\vartheta_\chi$ as in Proposition \ref{P:local_periodicity2}.
\end{proof}

\quad

\noindent\textbf{Global case:}

Let us consider the global case, and construct the (twisted)
exceptional representation
$\vartheta_\chi^\psi$ of $\GLtt_{2q}(\A)$. In a way, the construction
should be completely analogous to the $\GLt_{2q}(\A)$ case and one
would like to obtain $\vartheta_\chi^\psi$ as the representation
generated by the residues of the Eisenstein series. But a key ingredient missing for
this construction is the Langlands theory of Eisenstein series for the
group of the form $\GLtt_{2q}(\A)$. Probably, there is no danger to
assume that the theory of Eisenstein series for the
metaplectic group as developed in \cite{MW} can be carried over to
$\GLtt_{2q}(\A)$ at least to the extent necessary for the construction
of the exceptional representation. If one takes this for granted, the
exceptional representation $\vartheta_\chi^\psi$ of $\GLtt_{2q}(\A)$ can
be constructed in the same way as the exceptional representation
$\vartheta_\chi$ of $\GLt_{2q}(\A)$. However here we give an alternate
approach, in which we will show that the exceptional representation
$\vartheta_\chi^\psi$ of $\GLtt_{2q}(\A)$ is simply a constituent of the restriction to
$\GLtt_{2q}(\A)$ of the exceptional representation $\vartheta_\chi$ of
$\GLt_{2q}(\A)$. Here by restriction we mean the restriction of
automorphic forms as functions, not restriction of abstract
representation. Namely we have

\begin{Prop}\label{P:global_twisted_exceptional_even2}
For the global exceptional representation $\vartheta_\chi$ of
$\GLt_{2q}(\A)$, let $\vartheta_\chi^{(2)}$ be the representation of
$\GLtt_{2q}(\A)$ whose space is
$\{f|_{\GLtt_{2q}(\A)}:f\in\vartheta_\chi \}$, namely the space of
restrictions to $\GLtt_{2q}(\A)$ of automorphic forms in
$\vartheta_\chi$. Then we have
\[
    \vartheta_\chi^{(2)}=\bigoplus_{a\in\Sigma}\vartheta_\chi^{\psi^a},
\]
where $\vartheta_\chi^{\psi^a}$ is an irreducible quotient of the
global induced representation
$\Ind_{(\MPt)^{(2)}(\A)N_P(\A)}^{\GLtt_{2q}(\A)}\varpi_\chi^{\psi^a}\otimes\delta_P^{1/4}$,
and $\Sigma=(F^\times)^2\backslash F^\times$. 
\end{Prop}
\begin{proof}
Recall from our construction of $\vartheta_\chi$ in the previous
subsection that $\vartheta_\chi$ is constructed as the residual representation of the
Eisenstein series on $\Ind_{\Pt(\A)}^{\GLt_{2q}(\A)}\Pi_\chi^\nu$ at
$\nu=\rho_P/2$. For each $f^\nu$ in this
space, we have defined the Eisenstein series by
\[
    E(g,\Pi_\chi,f^\nu)=\sum_{\gamma\in P(F)\backslash \GL_{2q}(F)}f^\nu(\gamma g).
\]
But note that
\[
    P(F)\backslash \GL_{2q}(F)={(M_P)}^{(2)}(F)N_P(F)\backslash {\GL_{2q}^{(2)}}(F),
\]
and so one can write
\[
    E(g,\Pi_\chi,f^\nu)=\sum_{\gamma\in
      {(M_P)}^{(2)}(F)N_P(F)\backslash
      {\GL_{2q}^{(2)}}(F)}f^\nu(\gamma g).
\]
Hence we have
\[
 E(-,\Pi_\chi,f^\nu)|_{\GLtt_{2q}(\A)}= E(-,\Pi_\chi,f^\nu|_{\GLtt_{2q}(\A)}),
\]
where the latter may be called the ``Eisenstein series'' on $\GLtt_{2q}(\A)$.

To consider $f^\nu|_{\GLtt_{2q}(\A)}$, note that the induced representation
$\Ind_{(\MPt)^{(2)}(\A)N_P(\A)}^{\GLtt_{2q}(\A)}{\varpi_\chi^{\psi^a\;\nu}}$
is also realized in a space of automorphic forms on
$N_P(\A){M_P}^{(2)}(F)\backslash \GLtt_{2q}(\A)$. The map
$f^\nu\mapsto f^\nu|_{\GLtt_{2q}(\A)}$ gives a $\GLtt_{2q}(\A)$
surjection
\[
\Ind_{\MPt(\A)N_P(\A)}^{\GLt_{2q}(\A)}\Pi_\chi^\nu\longrightarrow
    \underset{a\in\Sigma}{\bigoplus}
\Ind_{(\MPt)^{(2)}(\A)N_P(\A)}^{\GLtt_{2q}(\A)}{\varpi_\chi^{\psi^a\;\nu}}.
\]
To see this, recall from (\ref{E:section}) and
(\ref{E:section_restriction}) that
$f^\nu=\phi\otimes\exp(\nu+\rho_P,H_P(\;))$ where $\phi$ is such that
the function $\phi_k$ is in $\Pi_\chi$. Then
$f^\nu|_{\GLtt_{2q}(\A)}=\phi|_{\GLtt_{2q}(\A)}\otimes
\exp(\nu+\rho_P,H_P(\;))$, where the Harish-Chandra map is also
restricted to $\GLtt_{2q}$. Hence the map $m\mapsto
\phi|_{\GLtt_{2q}(\A)}(mk)$ where $m\in\MPttt(\A)$ and
$k\in\widetilde{K}\cap\GLtt_{2q}(\A)$, is in $\Pi_\chi^{(2)}$ in the
notation of Proposition \ref{P:Weil_restriction2}. Then
use Proposition \ref{P:Weil_restriction2}.

Now if one chooses $f^\nu\in \Ind_{\MPt(\A)
  N_P(\A)}^{\GLt_{2q}(\A)}\Pi_\chi^\nu$ so that its image under the above
restriction map is in
$\Ind_{(\MPt)^{(2)}(\A)N_P(\A)}^{\GLtt_{2q}(\A)}{\varpi_\chi^{\psi^a\;\nu}}$
for a fixed $a\in\Sigma$, then the restriction of the Eisenstein series
$E(g,\Pi_\chi,f^\nu)$ is the Eisenstein series associated to
$\Ind_{(\MPt)^{(2)}(\A)N_P(\A)}^{\GLtt_{2q}(\A)}{\varpi_\chi^{\psi^a\;\nu}}$. 

Therefore at $\nu=\rho_P/2$, by taking the residues, we have the
commutative diagram of $\GLtt_{2q}(\A)$-intertwining maps
\[
\xymatrix{
\Ind_{\MPt(\A)
  N_P(\A)}^{\GLt_{2q}(\A)}\Pi_\chi^\nu\ar[r]\ar[d]_{E_{-1}}&
    \underset{a\in\Sigma}{\bigoplus}
\Ind_{(\MPt)^{(2)}(\A)N_P(\A)}^{\GLtt_{2q}(\A)}{\varpi_\chi^{\psi^a\;\nu}}
\ar[d]^{\oplus E_{-1}}\\\vartheta_\chi\ar[r]
&\underset{a\in\Sigma}{\bigoplus}V^a,
}
\]
where the vertical arrows are given by residue of Eisenstein series,
the horizontal arrows given by restriction of functions, and each $V^a$ is
the image of each
$\Ind_{(\MPt)^{(2)}(\A)N_P(\A)}^{\GLtt_{2q}(\A)}{\varpi_\chi^{\psi^a\;\nu}}$.

We need to show that each $V^a$ is $\vartheta_\chi^{\psi^a}$ as
claimed in the statement of the proposition, namely we need to show it is
irreducible. But the square integrability of
$\vartheta_\chi$ implies that of $V^a$, which implies complete
reducibility. Hence if $V^a$ is not irreducible, there is a place $v$
at which the $v$-component $V^a_v$ is a direct sum of more than two
irreducible representations of $\GLtt_{2q}(F_v)$. But since $V_v^a$ is
a quotient of the corresponding local induced representation, this
contradicts the uniqueness part of Proposition
\ref{P:local_twisted_exceptional_even2}. Hence $V^a$ is irreducible.
\end{proof}

Note that the representation $\vartheta_\chi^\psi$ is dependent on
$\psi$, and hence the notation.\\

Finally we need to prove the global periodicity property of
$\vartheta_\chi^\psi$.
\begin{Prop}[Global Periodicity]\label{P:global_periodicity3}
Let ${(\vartheta_\chi^\psi)}_{N_P}$ be the space generated by the constant
terms of automorphic forms in 
$\vartheta_\chi^\psi$ along the parabolic
$\Pt(\A)\cap\GLtt_{2q}(\A)=(\MPt)^{(2)}(\A)N_P(\A)$. Then we have 
\[
{(\vartheta_\chi^\psi)}_{N_P}={^{w_0}(\varpi_\chi^\psi)}\otimes\delta_P^{1/4}
=\varpi_\chi^\psi\otimes\delta_P^{1/4}.
\]
\end{Prop}
\begin{proof}
This follows from the above proposition and Proposition
\ref{P:global_periodicity2}, or strictly speaking from their
proofs. Namely from the proof of the above proposition, one knows that each element in
$\vartheta_\chi^\psi$ is written as
\[
\underset{\nu=\rho_P/2}{\Res}E(-,\Pi_\chi, f^\nu|_{\GLtt_{2q}(\A)})
\]
for some $f^\nu$ so that
$f^\nu|_{\GLtt_{2q}(\A)}\in
\Ind_{(\MPt)^{(2)}(\A)N_P(\A)}^{\GLtt_{2q}(\A)}\varpi_\chi^{\psi^a}\otimes\delta_P^\nu$,
and since the constant term is computed by integrating along
$N_P(F)\backslash N_P(\A)$, each element in
$(\vartheta_\chi^\psi)_{N_P}$ is generated by the elements of the form
\[
\underset{\nu=\rho_P/2}{\Res}E_P(-,\Pi_\chi, f^\nu|_{\GLtt_{2q}(\A)}),
\]
where the notation $E_P$ is as in the proof of Proposition
\ref{P:global_periodicity2}. But from the proof of Proposition
\ref{P:global_periodicity2}, we have
\[
\underset{\nu=\rho_P/2}{\Res}E_P(-,\Pi_\chi, f^\nu|_{\GLtt_{2q}(\A)})
=
\underset{\nu=\rho_P/2}{\Res}M(\nu, \Pi_\chi, w_0) f^\nu|_{\GLt_{2q}(\A)},
\]
where both sides are viewed as forms on $\MPttt(\A)$.
Here note that $M(\nu, \Pi_\chi, w_0) (f^\nu|_{\GLt_{2q}(\A)})=(M(\nu,
\Pi_\chi, w_0) f^\nu)|_{\GLt_{2q}(\A)}$ and that is why we can simply
write $M(\nu, \Pi_\chi, w_0) f^\nu|_{\GLt_{2q}(\A)}$. Viewed as forms on
$\MPttt(\A)$, we have 
\[
\underset{\nu=\rho_P/2}{\Res}M(\nu, \Pi_\chi, w_0)
f^\nu|_{\GLt_{2q}(\A)}\in
{^{w_0}(\varpi_\chi^\psi\otimes\delta_P^{1/4})}\otimes\delta_P^{1/2},
\]
and 
${^{w_0}(\varpi_\chi^\psi\otimes\delta_P^{1/4})}\otimes\delta_P^{1/2}=
{^{w_0}(\varpi_\chi^\psi)\otimes\delta_P^{1/4}}=\varpi_\chi^\psi\otimes\delta_P^{1/4}$.
\end{proof}

Finally, in this subsection let us mention that under
$\vartheta_\chi^\psi$, the center $\Zt$ of $\GLtt_{2q}$ acts by the
character
\begin{equation}\label{E:central_character2}
(1,\xi)\s(z)\mapsto\xi\chi(a)^q\mu_\psi(a)^q,\quad
z=\begin{pmatrix}a&&\\ &\ddots&\\ &&a\end{pmatrix}.
\end{equation}
This follows from (\ref{E:central_character}).

%%%%%%%%%%%%%%%%%%%%%%%%%%%%%%%%%%%%%%%%%%%%%%%%%%%%%%%%%%%%%%%%%%%

\subsection{\bf The semi-Whittaker functional on the exceptional
representation}

%%%%%%%%%%%%%%%%%%%%%%%%%%%%%%%%%%%%%%%%%%%%%%%%%%%%%%%%%%%%%%%%%%%

One of the key properties that we need for the exceptional
representations $\theta_\chi, \vartheta_\chi$ and
$\vartheta_\chi^\psi$ is that they do not possess Whittaker
functionals (unless $r=2$), but instead they possess what Bump and Ginzburg
call the semi-Whittaker functionals. This fact essentially follows
from the periodicity property for those exceptional representations.
To recall this notion, let us
define the character $\psi_N^{\e}$ on $N$ by
\[
\psi_N^{\e}
\begin{pmatrix}
1&x_{12}&\cdots&x_{1r}\\&1&&\vdots\\&&\ddots&x_{r-1,r}\\&&&1
\end{pmatrix}
=\psi(x_{r-1,r}+x_{r-3,r-2}+x_{r-5,r-4}+\cdots).
\] 
Then
\begin{Prop}\label{P:semi-Whittaker} 
Assume $F$ is non-archimedean, and $\theta$ is any of the exceptional
representations of $\GLt_r(F)$ (or $\GLtt_r(F)$). Then there is a
unique (up to scalar multiple) semi-Whittaker functional $L$ on
$\theta$, \ie a linear functional $L$ on $\theta$ such that
\[ 
L(\theta(\s(n))v)={\psi_N^{\e}}(n)L(v)
\] 
for $v\in V$ and $n\in N(F)$.
\end{Prop}
\begin{proof} 
For the non-twisted case $\theta=\theta_\chi$, this is
\cite[Theorem 1.4]{BG}. For the twisted case $\theta=\vartheta_\chi$
or $\vartheta_\chi^\psi$, this can be shown in the same way as
\cite{BG} by the periodicity property of $\vartheta_\chi$ and
$\vartheta_\chi^\psi$ (Proposition \ref{P:local_periodicity2} and
\ref{P:local_periodicity3}). For example, assume
$\theta=\theta_\chi^\psi$. Then the semi-Whittaker functional
is simply the composite of the surjection
$\vartheta_\chi^\psi\rightarrow\varpi_\chi^\psi\otimes\delta_P^{1/4}$ with the
Whittaker functional on the Weil representation
$\varpi_\chi^\psi$. The uniqueness follows from the uniqueness of the
Whittaker functional of the Weil representation.
\end{proof}

\begin{Rmk}\label{R:semi_Whittaker}
The important remark we have to make here is that the uniqueness of
the semi-Whittaker functional can be shown only for non-archimedean
$F$, because the proof requires
the periodicity of the Jacquet module of $\theta$, which is available
only for the non-archimedean $F$. Though this might hold for the
archimedean case as well, at this moment the author does not know if the
same technique can be applied to the archimedean case. (In \cite{BG}
it is simply stated without any proof or reference that the
periodicity holds for the archimedean case as well.) Because of this
lack of the uniqueness of the semi-Whittaker functional for the
archimedean places, it does not seem to be possible to prove the Euler
product of the Rankin-Selberg integral. But to get around this, we
obtain the ``almost Euler product", which is enough for our proof of
the main theorem of this paper.
\end{Rmk}

\quad

%%%%%%%%%%%%%%%%%%%%%%%%%%%%%%%%%%%%%%%%%%%%%%%%%%%%%%%%%%%%%%%%%%%

\subsection{\bf Exceptional representations of
  $\GLt_{r-1}\timest\GLt_1$}

%%%%%%%%%%%%%%%%%%%%%%%%%%%%%%%%%%%%%%%%%%%%%%%%%%%%%%%%%%%%%%%%%%%

Lastly in this section, we consider the exceptional representation of
$\GLt_{r-1}\timest\GLt_1$. Indeed the notion of the exceptional
representation can be generalized to the group
$\GLt_{r_1}\timest\cdots\timest\GLt_{r_k}$ both for twisted and
non-twisted cases following the method discussed in
\cite[p.142-143]{BG}. But here we specialize only to the case
$\GLt_{r-1}\timest\GLt_1$.

Let $Q$ be the $(r-1,1)$-parabolic subgroup of $\GL_r$ whose Levi part
is $M_Q=\GL_{r-1}\times\GL_1$.
Naively speaking, the exceptional representation of $\MQt$ is the
tensor product of the exceptional representation of $\GLt_{r-1}$ and a
character on $\GLt_1$. Things will slightly differ, depending on the parity of $r$.\\

\noindent\underline{Even case ($r=2q$)}
 
Assume $r$ is even, so $r=2q$. Fix $a\in F^\times$, where $F$ is either local or
global. Define a character $\omega_\chi^{\psi^a}$ on $\Tte$ by
\begin{equation}\label{E:exceptional_character2}
    {\omega_{\chi}^{\psi^a}}((1,\xi)\s(t))=\xi\chi(\det
    t)\mu_{\psi_a}(t_1)^{-1}\mu_\psi(t_3)^{-1}\mu_\psi(t_5)^{-1}\cdots \mu_\psi(t_{2q-1})^{-1}.
\end{equation}
Here unlike (\ref{E:exceptional_character}), we use
$\mu_\psi^{-1}$. Also for the first factor we use the character
$\psi_a$. This modification is needed for later purposes.

Let $B'$ be the Borel subgroup of $M_Q$, namely $M_Q\cap B$. For each
$\nu\in\Phi_{B'}$, we define
\[
\omega_{\chi}^{\psi^a\,\nu}:=
\omega_{\chi}^{\psi^a}\otimes\exp(\nu, H_{B'}(\;)).
\]
Then we have
\begin{Prop}\label{P:certain_exceptional}
The induced representation $\Ind_{\Tte
  N^\ast}^{\MQt}\omega_{\chi}^{\psi^a\,\nu}$ has a unique irreducible
quotient at $\nu=\rho_{B'}/2$, which we denote by
${\bar{\theta}_\chi}$. This is independent of $\psi$ and
$a$. If $F$ is global, it is a square integrable
automorphic representation in the residual spectrum of $\MQt(\A)$.
\end{Prop}
\begin{proof}
This is nothing but what Bump and Ginzburg call the exceptional
representation of $\GLt_{r-1}\timest\GLt_1$ in
\cite[p.142-143]{BG}. To show that it is independent of $\psi$ and
$a$, one can argue as we did for $\theta_\chi^\psi$ for the even case.
\end{proof}

For a (local or global) character $\eta$, define the character on
$\MQt$ by
\[
((g,a),\xi)\mapsto\eta(a)
\]
for $(g,a)\in M_Q$, namely the composition of $\eta$ with the projections
$\MQtt\rightarrow M_Q\rightarrow\GL_1$. We denote this character again
by $\eta$. We let
\[
\theta_{\chi, \eta}:=\bar{\theta}_\chi\otimes\eta,
\]
and call it the exceptional representation of $\MQt$ associated with
the characters $\chi$ and $\eta$.

Let $\MQttt$ be the metaplectic preimage of
\[
(M_Q)^{(2)}:=\{(g,a)\in\GL_{r-1}\times\GL_1: (\det g)a \text{ is a
  square}\}.
\]
As we did before for the exceptional representations of $\GLt_r$, the
restriction of $\theta_{\chi,\eta}$ to $\MQttt$ is described
as follows.
\begin{Prop}\label{certain_exceptional_restriction}
Assume $F$ is local. Then we have the decomposition
\[
\theta_{\chi,\eta}|_{\MQttt(F)}=\bigoplus_{a\in\Sigma}{\theta}_{\chi,\eta}^{\psi_a},
\]
where ${\theta}_{\chi,\eta}^{\psi_a}$ is a unique irreducible quotient of the
induced representation
$\Ind_{\Tte(F)N(F)^\ast}^{\MQttt(F)}\omega_\chi^{\psi^a}\otimes\delta_{B'}^{1/4}$,
and $\Sigma=(F^\times)^2\backslash F^\times$.

Assume $F$ is global and let $(\theta_{\chi,\eta})^{(2)}$ be
the space of the restrictions $\MQttt(\A)$ of the automorphic forms in
$\theta_{\chi,\eta}$. As representations of
$\MQttt(\A)$ we have the decomposition
\[
(\theta_{\chi,\eta})^{(2)}=\bigoplus_{a\in\Sigma}\theta_{\chi,\eta}^{\psi_a},
\]
where $\theta_{\chi,\eta}^{\psi_a}$ is an irreducible quotient of the
global induced representation
$\Ind_{\Tte(\A)N_A^\ast}^{\MQttt(\A)}\omega_\chi^{\psi^a}\otimes\delta_{B'}^{1/4}$,
and $\Sigma=(F^\times)^2\backslash F^\times$.
\end{Prop}
\begin{proof}
This can be proven in the same way as Proposition
\ref{P:local_twisted_exceptional_even2} and
\ref{P:global_twisted_exceptional_even2}.
\end{proof}

W call $\theta_{\chi,\eta}^{\psi^a}$ the exceptional representation
of $\MQttt$. We are mainly interested in $\theta_{\chi,\eta}^{\psi}$,
\ie $a=1$. Also let us note that $\Zt$ is in the center of $\MQttt$,
and each element $(1,\xi)\s(z)$ acts as
\begin{equation}\label{E:central_character3}
\theta_{\chi,\eta}^{\psi}((1,\xi)\s(z))=\xi\chi^{2q}(a)\eta(a)\mu_\psi(a)^{-q}, 
\quad z=\begin{pmatrix}a&&\\ &\ddots&\\ &&a\end{pmatrix}\in\GL_{2q}.
\end{equation}
As we mentioned in
(\ref{E:central_character1}), whether $q$ is even or odd, the map
$z\mapsto\mu_\psi(a)^{-q}$ is indeed a character.

\quad

\noindent\underline{Odd case ($r=2q+1$)}

Next we will consider the odd case. But this case is much simpler
because $\Zt$ is the center of $\GLt_{2q+1}$. First consider
$\cMQt$ as in Appendix \ref{S:tensor_product}. Note that $\Zt$ is in
the center of $\cMQt$, and $\GLt_{2q}$
naturally embeds into $\cMQt$ by $(g,\xi)\mapsto ((g,1),\xi)$. So we
have $\Zt\GLt_{2q}=\cMQt$.  Moreover inside $\cMQt$, we have
$\Zt\cap\GLt_{2q}=\{(1,\xi)\}$.

Let $\vartheta_{\chi}$ be the (local or global) exceptional
representation of $\GLt_{2q}$, where we include the case
$\chi^{1/2}$ exists, and $\eta$ a (local or global)
character. We define the representation
\[
    \vartheta_{\chi,\eta}
\]
of $\cMQt$ by extending $\vartheta_{\chi}$ to a representation of
$\Zt\GLt_{2q}=\cMQt$ by letting the center $\Zt\subseteq\cMQt$
act as
\begin{equation}\label{E:central_character4}
    (1,\xi)\s_Q(z)\mapsto \xi\chi(a)^q\eta(a)\mu_\psi(a)^{-q},
\quad z=\begin{pmatrix}a&&\\ &\ddots&\\ &&a\end{pmatrix},
\end{equation}
where $\s_Q:M_Q\rightarrow \cMQt$ is the map defined by $(g,a)\mapsto
((g,a), s_{2q}(g)^{-1}s_1(a)^{-1})$ for $g\in\GL_{2q}$ and
$a\in\GL_1$. (Strictly speaking it is a partial map if $F$ is global,
whose domain includes $B(\A)$ and hence $Z(\A)$.) As in the even case, the
map $z\mapsto\mu_\psi(a)^{-q}$ is a indeed a character for any $q$.

Now we identify  $\cMQt$ with the subgroup $\MQt$ of
$\GLt_{2q+1}$ via $\varphit_Q:\cMQt\rightarrow\MQt$. (See Appendix
\ref{S:tensor_product}.) Then we call the representation
$\vartheta_{\chi,\eta}\circ\varphit_Q^{-1}$ the exceptional
representation of $\MQt$, which we simply denote by
$\vartheta_{\chi,\eta}$ by abuse of notation. The central character
acts in the same way as in (\ref{E:central_character4}) with $\s_Q$
replaced by $\s$.

\quad

Finally for this subsection let us mention
\begin{Lem}\label{L:spherical_section2}
Let $r=2q$ or $2q+1$. Also assume $F$ is a non-archimedean local
field of odd residue characteristic. Further assume that
all of $\chi, \eta$ and $\psi$ are unramified. Consider the
intertwining operators
\begin{align*}
    &A(s,\theta^\psi_{\chi,\eta},w_0):
    \Ind_{\MQttt N_Q^\ast}^{\GLtt_{2q}}\theta_{\chi,\eta}^{\psi}\otimes\delta_Q^s\rightarrow 
\Ind_{(^{w_0}\MQttt) N^\ast_{(1.r-1)}}^{\GLtt_{2q}}\;^{w_0}(\theta_{\chi,\eta}^{\psi})\otimes\delta_Q^{-s},
\quad (r=2q)\\
&A(s,\vartheta_{\chi,\eta},w_0):
    \Ind_{\MQt N_Q^\ast}^{\GLt_{2q+1}}{\vartheta_{\chi,\eta}}\otimes\delta_Q^s\rightarrow 
\Ind_{(^{w_0}\MQt)N_{(1,r-1)}^\ast}^{\GLt_{2q+1}}\;^{w_0}({\vartheta_{\chi,\eta}})\otimes\delta_Q^{-s},
\quad (r=2q+1),
\end{align*}
 where $w_0=\begin{pmatrix}&1\\I_{r-1}&\end{pmatrix}$ and
 $N_{(1,r-1)}$ is the unipotent radical of the standard $(1,r-1)$-parabolic.

If $f_0^s\in\Ind_{\Qtt}^{\GLtt_{2q}}{\theta_{\chi,\eta}^{\psi}\otimes\delta_Q^s}$ (or
$\Ind_{\Qt}^{\GLt_{2q+1}}{\vartheta_{\chi,\eta}}\otimes\delta_Q^s $) is the
spherical section such that $f_0^s(1)=1$, then
\begin{align*}
   A(s,\theta^\psi_{\chi,\eta},w_0)f_0^s(1)&=
\frac{L(r(2s+\frac{1}{2})-r+1,
      \eta^{-2})}{L(r(2s+q+\frac{1}{2}), \eta^{-2})},\quad (r=2q)\\
A(s,\vartheta_{\chi,\eta},w_0) f_0^s(1)&=\frac{L(r(2s+\frac{1}{2})-r+1,
      \chi\eta^{-2})}{L(r(2s+q+\frac{1}{2}), \chi\eta^{-2})},\quad (r=2q+1).
\end{align*}
\end{Lem}
\begin{proof}
This is derived from the unramified computation by Kazhdan and
Patterson just as we did for Lemma \ref{L:spherical_section}. Since
the computation is straightforward, though quite tedious, we will omit
the details. Also this generalizes \cite[Proposition 5.6]{BG}. See the
proof there as well.
\end{proof}

\quad

%%%%%%%%%%%%%%%%%%%%%%%%%%%%%%%%%%%%%%%%%%%%%%%%%%%%%%%%%%%%%%%%%%%

\section{\bf The Rankin-Selberg integrals for the case $r=2q$}

%%%%%%%%%%%%%%%%%%%%%%%%%%%%%%%%%%%%%%%%%%%%%%%%%%%%%%%%%%%%%%%%%%%

In this section, we consider the Rankin-Selberg integral for the
cuspidal representation $\pi$ of $\GL_r(\A)$ when $r$ is even. So
throughout this section, we assume that
\[ 
r=2q=\text{even},
\] 
and $F$ is a number field, unless stated otherwise. We let
$\chi$ be a unitary Hecke character on $\A^\times$ and $\omega$ the
central character of $\pi$.

We let
\[ 
\theta:=\vartheta_{\chi\omega^{-2}}^\psi
\]
be the global twisted exceptional representation of $\GLtt_{r}(\A)$
associated with the character $\chi\omega^{-2}$. Also we let
\[ 
\theta':=\theta_{\omega,\;\omega^{-1}\chi^{-q}}
\] 
be the global exceptional representation of
$\MQt(\A)=\GLt_{r-1}(\A)\timest\GLt_1(\A)$ associated with the
characters $\omega$ and $\omega^{-1}\chi^{-q}$.

The global decomposition in Proposition \ref{P:certain_exceptional}
implies that we have $\GLtt_{2q}(\A)$ intertwining operator
\[
\ind_{\Qt(\A)}^{\GLt_{r}(\A)}\;\theta'\otimes\delta_{Q}^s\longrightarrow
\bigoplus_{a\in\Sigma}\ind_{\MQttt(\A)N_Q(\A)^\ast}^{\GLtt_{r}(\A)}
\;\theta_{\omega,\omega^{-1}\chi^{-q}}^{\psi^a}\otimes\delta_Q^s
\]
by restriction of functions. Here we assume the induction is NOT
normalized to be consistent with the convention in \cite{BG}.

In this section, unlike what we did in the previous section, we view each section
$f^s$ as a map 
\[
f^s:\GLt_r(\A)\longrightarrow \text{ space of }\theta'\otimes\delta_Q^s,
\]
and hence for each $\gt\in\GLt_r(\A)$, $f^s(\gt)$ is a function on
$\MQt(\A)$. We sometimes need to evaluate
$f^s(\gt)$ for each $\tilde{m}\in\MQt(\A)$. But we avoid
the notation $f^s(\gt)(\tilde{m})$, but rather use the notation
$f^s(\gt; \tilde{m})$. The advantage of this notation is that if we
have another $\tilde{m}_1\in\MQt(\A)$, then the translation of
$\tilde{m}_1$ from the first variable to the second is more naturally
written like $f^s(\tilde{m}_1\gt; \tilde{m})=f^s(\gt;
\tilde{m}\tilde{m}_1)$. 

Choose a section $f^s$ so that its image under the above surjection is
in
$\ind_{\Qtt(\A)}^{\GLtt_{r}(\A)}\;\theta_{\omega,\omega^{-1}\chi^{-q}}^{\psi}\otimes\delta_Q^s$,
\ie the component for $a=1$. Let
$E(-,s,f^s)$ be the Eisenstein series on $\GLt_{r}(\A)$ associated
with $f^s$. To
be precise,
\[ 
E(\gt,s,f^s)=\sum_{\gamma\in Q(F)\backslash
\GL_{r}(F)}f^s(\s(\gamma) \gt; e),
\] 
where $\gt\in\GLt_{r}(\A)$ and $e$ is the identity element in
$\GLt_{r}(\A)$. Note that the group $\GL_r(F)$ is viewed as a subgroup
of $\GLt_r(\A)$ via the splitting $\s$ and we simply write $\GL_r(F)$
for $\GL_r(F)^\ast$.
By an easy calculation, one sees that
\[ 
Q(F)\backslash \GL_{r}(F)=Q_{r-1}(F)\backslash \GL_{r}^{(2)}(F),
\] 
where
\[ 
Q_{r-1}(F):=Q(F)\cap \GL_{r}^{(2)}(F)=\{g\in Q:\det
g\in(F^\times)^2\}=\{\begin{pmatrix}h&\ast\\0&a\end{pmatrix}:(\det
h)a\in (F^\times)^2\},
\] 
and hence
\[ 
E(\gt,s,f^s)=\sum_{\gamma\in {Q_{r-1}(F)}\backslash
\GL_r^{(2)}(F)}f^s(\s(\gamma) \gt; e).
\] 
(The reason for the notation $Q_{r-1}$ will be clear in due
course.) Hence the restriction of the Eisenstein series $E(-,s,f^s)$ to
$\GLtt_r(\A)$ is the ``Eisenstein series'' on $\GLtt_r(\A)$ associated
with the induced representation
$\ind_{\Qtt(\A)}^{\GLtt_r(\A)}\theta_{\omega,\omega^{-1}\chi^{-q}}^{\psi}\otimes\delta_Q^s$.

Let $\Theta$ be an automorphic form in the space of $\theta$. Since
$\Theta(\gt)$ and $E(\gt,s,f^s)$ are genuine automorphic forms on
$\GLtt_{r}(\A)$, their product is a function on $\GLtwo_{r}(\A)$ in the
sense that
if $\gt\in\GLtt_{r}(\A)$ is any of the preimages of $g\in\GLtwo_{r}(\A)$,
then the function $g\mapsto \Theta(\gt)E(\gt,s,f^s)$ is independent of
the choice of $\gt$. 

Next let us consider how the center $\Zt(\A)$ acts. Let
$(1,\xi)\s(z)\in\Zt(\A)$ with $z=aI_r$ and $a\in\A^\times$. By
(\ref{E:central_character2}),
\[
\Theta((1,\xi)\s(z))=\xi\chi(a)^q\omega(a)^{-2q}\mu_\psi(a)^q\Theta(e),
\]
and by (\ref{E:central_character3}),
\[
E((1,\xi)\s(z), s,
f^s)=\xi\omega(a)^{2q-1}\chi^{-q}(a)\mu_\psi(a)^{-q}E(e,s,f^s),
\]
where $e$ is the identity element in $\GLtt_r(\A)$.
Hence on the
product $\Theta(-)E(-,s,f^s)$, the center acts as $\omega^{-1}$.\\

Now for a cusp form $\phi\in\pi$, the Rankin-Selberg integral we
consider is
\[
 Z(\phi,\Theta,f^s)=\underset{Z(\A)\GL_{r}^{(2)}(F)\backslash
\GL_{r}^{(2)}(\A)}{\int}\phi(g)\Theta(\kappa(g))E(\kappa(g),s,f^s) dg.
\] 
Note that since $\omega$ is the central character of $\pi$ and as we
have seen, on the product $\Theta(-)E(-,s,f^s)$ the center acts by
$\omega^{-1}$, this
integral is well-defined in the sense that the
center $Z(\A)$ acts trivially for the integrand. (Strictly speaking,
one needs to use the fact that the product $\Theta(\gt)E(\gt,s,f^s)$ is
independent of the choice of $\gt$ to check that the center acts by
$\omega^{-1}$.)
The reader should
notice that our integral differs from the
one in \cite[(3.4)]{BG}. (As is pointed out by Kable
(\cite[Appendix]{Kable}), the integral in \cite{BG} for the case $r=2q$ cannot
be well-defined.)  

However, if we define
$Z(\phi,\Theta,f^s)$ in this way, we cannot obtain the desired Euler
product simply by following the computation of \cite{BG}. Instead, we
have to take an alternate approach. But first note that
\begin{align*}
Z(\phi,\Theta,f^s)&=\underset{Z(\A)\GLtwo_{r}(F)\backslash
\GLtwo_{r}(\A)}{\int}\phi(g)\Theta(\kappa(g)) \sum_{\gamma\in
{Q_{r-1}(F)}\backslash {\GLtwo_{r}(F)}}f^s(\s(\gamma)\kappa(g); e) dg\\
&=\underset{Z(\A)\GLtwo_{r}(F)\backslash
\GLtwo_{r}(\A)}\int\sum_{\gamma\in {Q_{r-1}(F)}\backslash
{\GLtwo_{r}(F)}} \phi(\gamma g)\Theta(\s(\gamma)\kappa(g))f^s(\s(\gamma)
\kappa(g); e) dg.
\end{align*}

Now we would like to collapse the sum as usual using $\gamma$. To do
it, we would like to write
\[
\s(\gamma)\kappa(g)=\kappa(\gamma g).
\]
But
this equality does not hold in general. Yet, the fact that both
$\Theta$ and $f^s$ are genuine allows one to do this
manipulation. Let us explain this more in detail.
For each $\gamma$ and $g$, there is
$\xi=\xi(\gamma, g)\in\{\pm1\}$,
depending on $\gamma$ and $g$,  such that
\begin{equation}\label{E:n_and_g}
\s(\gamma)\kappa(g)=(1, \xi)\kappa(\gamma g).
\end{equation}
Since the induced representation is genuine, we have
\begin{align*}
 f^s(\s(\gamma)\kappa(g); e)= f^s((1, \xi)\kappa(\gamma g); e)
= \xi f^s(\kappa(\gamma g); e).
\end{align*}
The same consideration regarding $\s(\gamma)\kappa(g)$ and
$\kappa(\gamma g)$ applies
to $\Theta$. Then we have two $\xi$, one from $f^s$ and the other from
$\Theta$, and they get cancelled out when $f^s$ and $\Theta$ are
multiplied. Hence 
\begin{align}
\Theta(\s(\gamma)\kappa(g))f^s(\s(\gamma)
\kappa(g); e)
=\Theta(\kappa(\gamma g))f^s(
\kappa(\gamma g); e).
\end{align}
Namely the ``genuineness'' of $\Theta$ and $f^s$ takes care of the
discrepancy between $\s(\gamma)\kappa(g)$ and $\kappa(\gamma g)$.
This trick allows one to exchange $\s$ and $\kappa$ freely, as
long as one does the same to both $\Theta$ and $f^s$ at the same
time.  Since we need
this trick regularly, we call it the ``$\s-\kappa$ trick''.

Now we are allowed to collapse the sum and obtain

\begin{align*}
Z(\phi,\Theta, f^s)=\underset{Z(\A)Q_{r-1}(F)\backslash
\GLtwo_{r}(\A)}{\int}\phi(g)\Theta(\kappa(g))f^s(\kappa(g); e) dg.
\end{align*}

To get the desired Euler product, \cite{BG} used the well-know Fourier
expansion of the cusp form $\phi$ and collapsed the sum in the Fourier
expansion. But in our case, their method does not work because the
integral is over $\GLtwo_{r}$ instead of $\GL_r$, or
roughly put, the group $Q_{r-1}(F)$ is not large enough to collapse
the sum in the Fourier expansion of $\phi$. To get around this, we
carry out ``column-by-column" computations of the Fourier coefficients
of $\Theta$ and $f^s$ together with some of the properties of the
exceptional representations we developed in the previous section.

\quad

But before going into the computation, let us fix notations. Let $N$
be the unipotent radical of the Borel subgroup $B$ of $\GL_r$. For an
integer $1\leq m\leq r-1$, we define $N_m$ to be the subgroup of $N$
consisting of the elements whose only non-zero entries off the
diagonal are in the $m+1^{\text{st}}$ column, \ie
\[
N_m=\{\begin{pmatrix}I_{m}&\ast&0\\&1&0\\&&I_{r-m-1}\end{pmatrix}\}.
\] 
Note that
\[ 
N=N_{r-1}N_{r-2}\cdots N_1.
\] 
Also note that $N_m(F)\backslash N_m(\A)$ is a compact abelian
group, which is isomorphic to $(F\backslash\A)^m$. Since we use this
group so frequently, let us put
\[ 
[N_m]=N_m(F)\backslash N_m(\A).
\] 
Also each element in $[N_m]$ is often denoted by the symbol
$n_m$. Now for our fixed additive character $\psi$ on $F\backslash\A$
and $a\in F^\times$, we define the character $\psi_N^a$ on $N(\A)$ by
\[ 
\psi_N^a\begin{pmatrix}1&x_{12}&\cdots&
x_{1r}\\&1&&\vdots\\&&\ddots&x_{r-1\, r}\\&&&1\end{pmatrix}
=\psi(ax_{12}+\sum_{i=2}^{r-1}x_{i\,i+1}).
\] 
We write $\psi_N^1=\psi_N$, which is the one we usually use. We
often consider $\psi_N$ restricted to $N_m(\A)$, which we also denote
by $\psi_N$.

We let
\begin{align*}
H_{m}&=\{\begin{pmatrix}g_{m}&\\&aI_{r-m}\end{pmatrix}\in\GLtwo_r:g_{m}\in\GL_{m},
a\in\GL_1\},
\end{align*} 
so the product $(\det g_m)a^{r-m}$ is a square. Note
$H_{m-1}\subseteq H_m$. We let
\[
Q_{m-1}=H_{m-1}N_{m-1}=
\{\begin{pmatrix}g_{m-1}&n&\\&a&\\&&aI_{r-m}\end{pmatrix}\}
\subseteq H_{m},
\] 
where we assume $Q_0=\{a I_r:a\in\GL_1\}$. Also notice that
$Q_{r-1}=Q\cap\GLtwo_r$, and hence our previous notation for this
group.

For our cusp form $\phi$ we write the ``partial Whittaker coefficient"
by
\begin{equation}\label{E:partial_Whittaker}
 W_{m}(g):=\underset{[N_{m}]}{\int}\underset{[N_{m+1}]}{\int}\cdots
\underset{[N_{r-1}]}{\int}\phi(n_{r-1}n_{r-2}\cdots n_{m}g)
\psi_N(n_{r-1}n_{r-2}\cdots n_{m})dn_{r-1}dn_{r-2}\cdots dn_{m}.
\end{equation}
Strictly speaking $W_{m}$ depends on the choice of $\psi$ and
$\phi$ but we use this notation since it will not produce any
confusion. The following property of this partial Whittaker
coefficient will be necessary for our computation.
\begin{Lem}\label{L:partial_Whittaker} For any $h_{m-1}\in
H_{m-1}(F)$, one has
\[ W_{m}(h_{m-1}g)=W_{m}(g).
\]
\end{Lem}
\begin{proof} 
This follows from the automorphy of $\phi$ and the fact
that $h_{m-1}$ fixes $\psi_N$ in the sense that
$\psi_N(h_{m-1}n_{r-1}n_{r-2}\cdots
n_{m}h_{m-1}^{-1})=\psi_N(n_{r-1}n_{r-2}\cdots n_{m})$.
\end{proof}

For our unfolding argument we compute the Fourier expansions
of $\Theta$ and $f^s$ along $N_{m}$ and $N_{m-1}$
``alternatingly". Namely first we consider the Fourier expansion of
$\Theta$ along $N_{r-1}$, and then that of $f^s$ along $N_{r-2}$, and
then $\Theta$ along $N_{r-3}$ and then $f^s$ along $N_{r-4}$, etc. For
this computation, the following lemma plays a pivotal role.
\begin{Lem}\label{L:Fourier} 
Assume $m\geq 2$. The group $H_m(F)$ acts on
the dual space $\widehat{[N_m]}\cong F^m$ by conjugation as
\[ (h_m\cdot\psi)(n_m)=\psi(h_m^{-1}n_m{h_m}),\quad h_m\in H_m(F),
n_m\in N_m(\A)
\] with two orbits: the zero orbit and the orbit of $\psi_N$, where
$\psi_N$ is actually the restriction of $\psi_N$ to
$N_m(\A)$. Moreover the stabilizer of $\psi_N$ is $Q_{m-1}(F)$, and
hence the orbit of $\psi_N$ is indexed by $Q_{m-1}(F)\backslash
H_m(F)$.
\end{Lem}
\begin{proof} Straightforward computation.
\end{proof}

Now we are ready to work out our Rankin-Selberg integral to obtain the
desired Euler product. Recall we have obtained
\begin{equation}\label{E:zeta0}
Z(\phi,\Theta,f^s)=\underset{Z(\A)Q_{r-1}(F)\backslash
\GLtwo_r(\A)}{\int}\phi(g)\Theta(\kappa(g))f^s(\kappa(g); e) dg.
\end{equation}

In what follows, one should keep in mind
Lemma\ref{L:cocycle_generated} along with the fact that the partial
section $\s:\GL_r(\A)\rightarrow\GLt_r(\A)$ is not only defined but
also is a homomorphism on both of the groups $\GL_r(F)$ and
$N_B(\A)$.

\quad\\

\noindent\textbf{\underline{Unfolding Step 1}}\\

The first step starts with computing the Fourier expansion of $\Theta$ along
the ``last column" $N_{r-1}$. Consider $n_{r-1}\mapsto
\Theta(\s(n_{r-1})\kappa(g))$ as a function on $[N_{r-1}]$, and expand it. Then
one has
\begin{align*}
\Theta(\kappa(g))=\sum_{\psi\in\widehat{[N_{r-1}]}}\underset{[N_{r-1}]}{\int}
\Theta(\s(n_{r-1})\kappa(g))\psi(n_{r-1})^{-1}\,dn_{r-1}.
\end{align*}
By Lemma \ref{L:Fourier} we obtain
\begin{align*}
\Theta(\kappa(g))&=\underset{[N_{r-1}]}{\int}\Theta(\s(n_{r-1})\kappa(g))\,dn_{r-1}\\
&\qquad\qquad+\sum_{h_{r-1}\in Q_{r-2}(F)\backslash
H_{r-1}(F)}\underset{[N_{r-1}]}{\int}\Theta(\s(n_{r-1}h_{r-1})\kappa(g))
\psi_N(n_{r-1})^{-1}\,dn_{r-1}.
\end{align*}
By substituting this expression of $\Theta(\kappa(g))$ in (\ref{E:zeta0}),
one obtains
\begin{align*} 
&Z(\phi,\Theta,f^s)=\underset{Z(\A)Q_{r-1}(F)\backslash
\GLtwo_r(\A)}{\int}\phi(g)\\
&\quad\Bigg(\underset{[N_{r-1}]}{\int}\Theta(\s(n_{r-1})\kappa(g))\,dn_{r-1}\\
&\qquad\qquad
+\sum_{h_{r-1}\in Q_{r-2}(F)\backslash H_{r-1}(F)}\underset{[N_{r-1}]}{\int}
\Theta(\s(n_{r-1}h_{r-1})\kappa(g))\psi_N(n_{r-1})^{-1}\,dn_{r-1}\Bigg)
f^s(\kappa(g); e)\;dg.
\end{align*}

One of the key points in our computation is that the term coming from
the zero orbit (``zero orbit term'') vanishes because of the cuspidality
of $\phi$. To see it, we would like to multiply out the large
parentheses and write out the zero orbit term separately as
\begin{equation}\label{E:zero_orbit1}
\underset{Z(\A)Q_{r-1}(F)\backslash
\GLtwo_r(\A)}{\int}\phi(g)
\underset{[N_{r-1}]}{\int}\Theta(\s(n_{r-1})\kappa(g))\,dn_{r-1}
f^s(\kappa(g); e)\;dg.
\end{equation}
But we need justification for this process because we need to know
that the product 
\begin{equation}\label{E:zero_orbit2}
\phi(g)\underset{[N_{r-1}]}{\int}\Theta(\s(n_{r-1})\kappa(g))\,dn_{r-1}
f^s(\kappa(g); e),
\end{equation}
 viewed as a function on $g$, is indeed invariant on
$Z(\A)Q_{r-1}(F)$ so that we can carry out the integration for
$dg$. This is not immediately clear. To see it, let $h\in
Q_{r-1}(F)$. First of all, by the $\s-\kappa$ trick we introduced before,
we have
\begin{align}\label{E:zero_orbit3}
\notag&\underset{[N_{r-1}]}{\int}\Theta(\s(n_{r-1})\kappa(hg))\,dn_{r-1}
f^s(\kappa(hg); e)\\
=&\underset{[N_{r-1}]}{\int}\Theta(\s(n_{r-1})\s(h)\kappa(g))\,dn_{r-1}
f^s(\s(h)\kappa(g); e).
\end{align}

Now $f^s(\s(h)\kappa(g); e)=f^s(\kappa(g); \s(h))=f^s(\kappa(g);e)$ by
the automorphy of $f^s$. Also 
\allowdisplaybreaks\begin{align*}
&\underset{[N_{r-1}]}{\int}\Theta(\s(n_{r-1})\s(h)\kappa(g))\,dn_{r-1}\\
=&\underset{[N_{r-1}]}{\int}\Theta(\s(n_{r-1}h)\kappa(g))\,dn_{r-1}\quad
\text{by Lemma \ref{L:cocycle_generated}}\\
=&\underset{[N_{r-1}]}{\int}\Theta(\s(hh^{-1}n_{r-1}h)\kappa(g))\,dn_{r-1}\\
=&\underset{[N_{r-1}]}{\int}\Theta(\s(hn_{r-1})\kappa(g))\,dn_{r-1}\quad
\text{by change of variable for $n_{r-1}$}\\
=&\underset{[N_{r-1}]}{\int}\Theta(\s(h)\s(n_{r-1})\kappa(g))\,dn_{r-1}\quad
\text{by Lemma \ref{L:cocycle_generated}}\\
=&\underset{[N_{r-1}]}{\int}\Theta(\s(n_{r-1})\kappa(g))\,dn_{r-1}\quad
\text{by automorphy of $\Theta$}.
\end{align*}

Hence (\ref{E:zero_orbit3}) becomes
\[
\underset{[N_{r-1}]}{\int}\Theta(\s(n_{r-1})\kappa(g))\,dn_{r-1}
f^s(\kappa(g); e).
\]
Therefore indeed (\ref{E:zero_orbit1}) viewed as a function of $g$ is
left invariant on $Q_{r-1}(F)$. Similarly one can see that it is left
invariant on $Z(\A)$ by using the $\s-\kappa$ trick and the actions of
the center on $\Theta$ and $f^s$.

Thus the expression (\ref{E:zero_orbit1}) makes sense and we can work
on this integral. Indeed, we will show it is zero. For this,
note that we can write
\begin{align*} 
\int_{Z(\A)Q_{r-1}(F)\backslash\GLtwo_r(\A)}
&=\int_{Z(\A)N_{r-1}(F)H_{r-1}(F)\backslash \GLtwo_r(\A)}\\
&=\int_{Z(\A)H_{r-1}(F)N_{r-1}(\A)\backslash \GLtwo_r(\A)}
\int_{N_{r-1}(F)\backslash N_{r-1}(\A)}.
\end{align*} 
Then we can write the outer integral of (\ref{E:zero_orbit1}) as an integral over
those two sets $N_{r-1}(F)\backslash N_{r-1}(\A)$ and
$Z(\A)H_{r-1}(F)N_{r-1}(\A)\backslash \GLtwo_r(\A)$, whose
corresponding variables we denote by $n'_{r-1}$ and $g$ respectively,
so all the occurrences of $g$ in the integrant are replaced by
$n'_{r-1}g$, and $dg$ is replaced by $dn'_{r-1}dg$. Then the zero
orbit term (\ref{E:zero_orbit1}) is written as
\begin{align*}
&\underset{Z(\A)H_{r-1}(F)N_{r-1}(\A)\backslash
\GLtwo_r(\A)}{\int}\;\underset{[N_{r-1}]}{\int}\phi(n_{r-1}'g)\\
&\qquad\qquad\qquad\qquad
\underset{[N_{r-1}]}{\int}\Theta(\s(n_{r-1})\kappa(n_{r-1}'g))\,dn_{r-1}
f^s(\kappa(n_{r-1}'g);e)\;dn_{r-1}'\;dg.
\end{align*}

By using the $\s-\kappa$ trick, this is written as
\begin{align*}
&\underset{Z(\A)H_{r-1}(F)N_{r-1}(\A)\backslash
\GLtwo_r(\A)}{\int}\;\underset{[N_{r-1}]}{\int}\phi(n_{r-1}'g)\\
&\qquad\qquad\qquad\qquad
\underset{[N_{r-1}]}{\int}\Theta(\s(n_{r-1}n_{r-1}')\kappa(g))\,dn_{r-1}\,dn_{r-1}'
f^s(\kappa(g);e)\;dg,
\end{align*}
where we also used the fact that
$\s(N_B(\A))$ acts trivially on $f^s$, and 
$\s$ is a homomorphism on $N_B(\A)$.
By the invariance of the measure $dn_{r-1}$ for the integral for
$\Theta$, the two inner integrals are written as
\[ 
\underset{[N_{r-1}]}{\int}\phi(n'_{r-1}g)\,dn'_{r-1}
\underset{[N_{r-1}]}{\int}\Theta(\s(n_{r-1})\kappa(g))\,dn_{r-1}.
\] 
The cuspidality of $\phi$ makes this term vanish.

\quad

Therefore we obtain
\begin{align*} 
&Z(\phi,\Theta,f^s)\\
&=\underset{Z(\A)Q_{r-1}(F)\backslash \GLtwo_r(\A)}{\int}\phi(g)\\
&\quad \left(
\sum_{h_{r-1}\in Q_{r-2}(F)\backslash
H_{r-1}(F)}\underset{[N_{r-1}]}{\int}\Theta(\s(n_{r-1}h_{r-1})\kappa(g))\psi_N(n_{r-1})^{-1}
\,dn_{r-1}\right)
f^s(\kappa(g);e) dg.
\end{align*} 

For each $h_{r-1}\in H_{r-1}(F)$ one sees that
$f^s(\kappa(g); e) =f^s(\kappa(g);\s(h_{r-1})) =f^s(\s(h_{r-1})\kappa(g);e)$ by
the automorphy. We move around $\s$ and $\kappa$ by the
$\s-\kappa$ trick, and one can see that the above
integral is written as
\begin{align*} 
&\underset{Z(\A)Q_{r-1}(F)\backslash \GLtwo_r(\A)}{\int}\phi(h_{r-1}g)\\
&\quad \left(
\sum_{h_{r-1}\in Q_{r-2}(F)\backslash
H_{r-1}(F)}\underset{[N_{r-1}]}{\int}\Theta(\s(n_{r-1})\kappa(h_{r-1}g))\psi_N(n_{r-1})^{-1}
\,dn_{r-1}\right)
f^s(\kappa(h_{r-1}g);e) dg,
\end{align*} 
where we used the automorphy of $\phi$ as well as
$\s(n_{r-1}h_{r-1})=\s(n_{r-1})\s(h_{r-1})$ by Lemma \ref{L:cocycle_generated}.
Then we are allowed to
collapse the sum by using
$Q_{r-1}(F)$, and obtain
\begin{align*} &Z(\phi,\Theta,f^s)\\
&=\underset{Z(\A)N_{r-1}(F)Q_{r-2}(F)\backslash
\GLtwo_r(\A)}{\int}\phi(g)\Bigg(
\underset{[N_{r-1}]}{\int}\Theta(\s(n_{r-1})\kappa(g))\psi_N(n_{r-1})^{-1}\,dn_{r-1}\Bigg)
f^s(\kappa(g); e) dg.
\end{align*}

Note that
\[ 
\int_{Z(\A)N_{r-1}(F)Q_{r-2}(F)\backslash \GLtwo_r(\A)}
=\int_{Z(\A)Q_{r-2}(F)N_{r-1}(\A)\backslash \GLtwo_r(\A)}
\int_{N_{r-1}(F)\backslash N_{r-1}(\A)}.
\] 
So we have
\begin{align*} 
Z(\phi,\Theta,f^s)
&=\underset{Z(\A)Q_{r-2}(F)N_{r-1}(\A)\backslash
\GLtwo_r(\A)}{\int}\;\underset{[N_{r-1}]}{\int}\phi(n_{r-1}'g)\\
&\qquad\qquad\underset{[N_{r-1}]}{\int}
\Theta(\s(n_{r-1})\kappa(n_{r-1}'g))\psi_N(n_{r-1})^{-1}\,dn_{r-1}
f^s(\kappa(n_{r-1}'g); e) dn_{r-1}'dg.
\end{align*} 

By using the $\s-\kappa$ trick, this is written as
\begin{align*} 
Z(\phi,\Theta,f^s)
&=\underset{Z(\A)Q_{r-2}(F)N_{r-1}(\A)\backslash
\GLtwo_r(\A)}{\int}\;\underset{[N_{r-1}]}{\int}\phi(n_{r-1}'g)\\
&\qquad\qquad\underset{[N_{r-1}]}{\int}
\Theta(\s(n_{r-1}n_{r-1}')\kappa(g))\psi_N(n_{r-1})^{-1}\,dn_{r-1}
f^s(\kappa(g); e) dn_{r-1}'dg,
\end{align*} 
where we also used $f^s(\s(n_{r-1}')\kappa(g);e)=f^s(\kappa(g);e)$.
The change of variable $n_{r-1}n_{r-1}'\mapsto n_{r-1}$
in the inner most integral gives
\begin{align*} 
Z(\phi,\Theta,f^s)
&=\underset{Z(\A)Q_{r-2}(F)N_{r-1}(\A)\backslash
\GLtwo_r(\A)}{\int}\;\underset{[N_{r-1}]}{\int}\phi(n_{r-1}'g)\psi(n_{r-1}')dn_{r-1}'\\
&\qquad\qquad\underset{[N_{r-1}]}{\int}\Theta(\s(n_{r-1})\kappa(g))
\psi_N(n_{r-1})^{-1}\,dn_{r-1}f^s(\kappa(g); e) dg.
\end{align*} 
So we have
\begin{align*}
Z(\phi,\Theta,f^s)
&=\underset{Z(\A)Q_{r-2}(F)N_{r-1}(\A)\backslash \GLtwo_r(\A)}{\int}
W_{r-1}(g)\\
&\qquad\qquad\underset{[N_{r-1}]}{\int}\Theta(\s(n_{r-1})\kappa(g))
\psi_N(n_{r-1})^{-1}\,dn_{r-1}
f^s(\kappa(g);e) dg,
\end{align*} 
recalling the notation for $W_{r-1}(g)$ from (\ref{E:partial_Whittaker}).

\quad\\

\noindent\textbf{\underline{Unfolding Step 2}}\\

The second step starts with computing the Fourier expansion of $f^s(\kappa(g);-)$ along
the ``$r-1^\text{st}$ column" $N_{r-2}$. Recall that $f^s(\kappa(g);-)$ is an
automorphic form on $\GLt_{r-1}\timest\GLt_1$. By viewing the
function $n_{r-2}\mapsto f^s(\kappa(g);\s(n_{r-2}))$ as a function on
$[N_{r-2}]$, and expanding it by using Lemma \ref{L:Fourier}, we obtain
\begin{align*}
f^s(\kappa(g); e)=&\underset{[N_{r-2}]}{\int}f^s(\s(n_{r-2})\kappa(g);
e)\,dn_{r-2}\\
&+\sum_{h_{r-2}\in Q_{r-3}(F)\backslash
H_{r-2}(F)}\underset{[N_{r-2}]}{\int}f^s(\kappa(g);\s(n_{r-2}h_{r-2}))
\psi_N(n_{r-2})^{-1}\,dn_{r-2}.
\end{align*}

Now the first term (the zero orbit),  when integrated with
the cusp form $\phi$, vanishes by the cuspidality of
$\phi$ as we have seen in Step 1. The idea is
essentially the same but the computation is not completely
identical. Hence we will give the detailed computation here.

First if the above Fourier expansion of $f^s$ is plugged in to the
formula for the $Z(\phi, \Theta, f^s)$ we obtained at the end of
Step 1, the product of  $\underset{[N_{r-1}]}{\int}\Theta(\s(n_{r-1})\kappa(g))
\psi_N(n_{r-1})^{-1}\,dn_{r-1}$ with each term in the Fourier
expansion of $f^s$ viewed as a function on $g$ is invariant on
$Z(\A)Q_{r-2}(F)N_{r-1}(\A)$. (To see this, once again we need the $\s-\kappa$
trick.) Hence the expression for $Z(\phi,
\Theta, f^s)$ after the Fourier expansion of $f^s$ is plugged in can
be expanded and we can take out the zero orbit term separately as
\begin{align*}
&\underset{Z(\A)Q_{r-2}(F)N_{r-1}(\A)\backslash \GLtwo_r(\A)}{\int}
W_{r-1}(g)\\
&\qquad\qquad\underset{[N_{r-1}]}{\int}\Theta(\s(n_{r-1})\kappa(g))
\psi_N(n_{r-1})^{-1}\,dn_{r-1}
\underset{[N_{r-2}]}{\int}f^s(\s(n_{r-2})\kappa(g);
e)\,dn_{r-2} dg.
\end{align*}
Since
\begin{align*}
\int_{Z(\A)Q_{r-2}(F)N_{r-1}(\A)\backslash \GLtwo_r(\A)}
&=\int_{Z(\A)N_{r-2}(F)H_{r-2}(F)N_{r-1}(\A)\backslash
  \GLtwo_r(\A)}\\
&=\int_{Z(\A)H_{r-2}(F)N_{r-2}(\A)N_{r-1}(\A)\backslash
  \GLtwo_r(\A)}\int_{N_{r-2}(F)\backslash N_{r-2}(\A)}
\end{align*}
the outer integral can be written as integrals over those two sets
\[
Z(\A)N_{r-2}(F)N_{r-2}(\A)N_{r-1}(\A)\backslash\GLtwo_r(\A)\quad
\text{ and } \quad N_{r-2}(F)\backslash N_{r-2}(\A),
\] 
whose
corresponding variables we denote by $g$ and $n_{r-2}'$
respectively. Then all the occurrences of $g$ are replaced by
$n_{r-2}'g$. Namely the zero orbit term is written as
\begin{align*}
&{\int}\underset{[N_{r-2}]}{\int} W_{r-1}(n_{r-2}'g)\\
&\underset{[N_{r-1}]}{\int}\Theta(\s(n_{r-1})\kappa(n_{r-2}'g))
\psi_N(n_{r-1})^{-1}\,dn_{r-1}
\underset{[N_{r-2}]}{\int}f^s(\s(n_{r-2})\kappa(n_{r-2}'g);
e)\,dn_{r-2}\,dn_{r-2}'\,dg,
\end{align*}
where the outer integral is over the set 
$Z(\A)H_{r-2}(F)N_{r-2}(\A)N_{r-1}(\A)\backslash\GLtwo_r(\A)$. By
using the $\s-\kappa$ trick, we can
write $\s(n_{r-1})\kappa(n_{r-2}'g)=\s(n_{r-1}n_{r-2}')\kappa(g)$
inside $\Theta$ and $f^s$. So the integral becomes
\begin{align*}
&{\int}\underset{[N_{r-2}]}{\int} W_{r-1}(n_{r-2}'g)\\
&\underset{[N_{r-1}]}{\int}\Theta(\s(n_{r-1}n_{r-2}')\kappa(g))
\psi_N(n_{r-1})^{-1}\,dn_{r-1}
\underset{[N_{r-2}]}{\int}f^s(\s(n_{r-2}n_{r-2}')\kappa(g);
e)\,dn_{r-2}\,dn_{r-2}'\,dg.
\end{align*}
By the invariance of the measure $dn_{r-2}$, this is written as
\begin{align*}
&{\int}\underset{[N_{r-2}]}{\int} W_{r-1}(n_{r-2}'g)\\
&\underset{[N_{r-1}]}{\int}\Theta(\s(n_{r-1}n_{r-2}')\kappa(g))
\psi_N(n_{r-1})^{-1}\,dn_{r-1}
\,dn_{r-2}' \underset{[N_{r-2}]}{\int}f^s(\s(n_{r-2})\kappa(g);
e)\,dn_{r-2}\,dg.
\end{align*}
Hence to show the vanishing of the zero orbit term, it suffices to
show
\begin{equation}\label{E:vanishing2}
\underset{[N_{r-2}]}{\int} W_{r-1}(n_{r-2}'g)
\underset{[N_{r-1}]}{\int}\Theta(\s(n_{r-1}n_{r-2}')\kappa(g))
\psi_N(n_{r-1})^{-1}\,dn_{r-1}\,dn_{r-2}'=0.
\end{equation}

To proceed, we need following crucial property of the exceptional
representations.

\begin{Lem}\label{L:key_lemma} Let $\Theta$ be an automorphic form in
the space of the exceptional representation $\vartheta_\chi$,
$\vartheta_\chi^\psi$ or $\theta_\chi$ of $\GLt_{r}(\A)$,
$\GLtt_{r}(\A)$ or $\GLt_r(\A)$, respectively, where $r$ can be either
$2q$ or $2q+1$. For an integer $1\leq m\leq q$, the integral
\begin{align*}
&\underset{[N_{r-2m+1}]}{\int}\cdots\underset{[N_{r-2}]}{\int}\underset{[N_{r-1}]}{\int}
\Theta(\s(n_{r-1}n_{r-2}n_{r-3}\cdots
n_{r-2m})\kappa(g))\psi_N(n_{r-1}n_{r-3}n_{r-5}\cdots n_{r-2m+1})\\
&\qquad\qquad\qquad dn_{r-1}dn_{r-2}\cdots dn_{r-2m+1}
\end{align*} is independent of $n_{r-2m}\in N_{r-2m}(\A)$.

Consequently, by integrating over $[N_{r-2m}]$, this integral is equal
to
\begin{align*}
&\underset{[N_{r-2m}]}{\int}\cdots\underset{[N_{r-2}]}{\int}\underset{[N_{r-1}]}{\int}
\Theta(\s(n_{r-1}n_{r-2}n_{r-3}\cdots
n_{r-2m})\kappa(g))\psi_N(n_{r-1}n_{r-3}n_{r-5}\cdots n_{r-2m+1})\\
&\qquad\qquad\qquad dn_{r-1}dn_{r-2}\cdots dn_{r-2m},
\end{align*} provided the measure is so chosen that the volume of
$[N_{r-2m}]$ is $1$.
\end{Lem}
\begin{proof} The case for $\theta_\chi$ is Proposition 2.4 and 2.5 of
\cite{BG}. The case for $\vartheta_\chi$ can be proven
identically. The key ingredient for the proof is the non-existence of
the Whittaker functional for the exceptional representation, which
implies Proposition 2.1 of \cite{BG}. Once the case for
$\vartheta_\chi$ is taken care of, the case for $\vartheta_\chi^\psi$
trivially follows because any automorphic form in the space of
$\vartheta_\chi^\psi$ is simply the restriction of an automorphic form
in the space of $\vartheta_\chi$.
\end{proof}

By applying the first part of the lemma with $m=1$, the left hand side
of (\ref{E:vanishing2}) is written as
\[
\underset{[N_{r-2}]}{\int} W_{r-1}(n_{r-2}'g)\,dn_{r-2}'
\underset{[N_{r-1}]}{\int}\Theta(\s(n_{r-1})\kappa(g))
\psi_N(n_{r-1})^{-1}\,dn_{r-1}.
\]

By definition of $W_{r-1}$, which is given in
(\ref{E:partial_Whittaker}), we have
\begin{align*}
\underset{[N_{r-2}]}{\int} W_{r-1}(n_{r-2}'g)\,dn_{r-2}'
&=\underset{[N_{r-2}]}{\int}\underset{[N_{r-1}]}{\int}
 \phi (n_{r-1}n_{r-2}'g)\psi_N(n_{r-1})\,dn_{r-1}\,dn_{r-2}'\\
&=\underset{[A]}{\int}\left(\underset{[N_{(r-2, 2)}]}{\int}
 \phi (n_{(r-2,2)}ag)\,dn_{(r-2,2)}\right)\psi_N(a)\,da,
\end{align*}
where $N_{(r-2,2)}$ is the unipotent radical of the $(r-2,
2)$-parabolic, 
and $A$ is the set of the matrices of the form
\[
a=\begin{pmatrix}I_{r-2}&&\\ &1&\ast\\ &&1\end{pmatrix}.
\]
By the cuspidality of $\phi$, the inner integral is zero. This shows
that the zero orbit term vanishes.

\quad

Hence we obtain
\allowdisplaybreaks\begin{align*} &Z(\phi,\Theta,f^s)\\
=&\underset{Z(\A)Q_{r-2}(F)N_{r-1}(\A)\backslash
\GLtwo_r(\A)}{\int}W_{r-1}(g)
\underset{[N_{r-1}]}{\int}\Theta(\s(n_{r-1})\kappa(g))\psi_N(n_{r-1})^{-1}\,dn_{r-1}\\
&\qquad\sum_{h_{r-2}\in Q_{r-3}(F)\backslash
H_{r-2}(F)}\underset{[N_{r-2}]}{\int}
f^s(\s(n_{r-2}h_{r-2})\kappa(g);e)\psi_N(n_{r-2})^{-1}\,dn_{r-2}dg\\
=&\underset{Z(\A)Q_{r-2}(F)N_{r-1}(\A)\backslash
\GLtwo_r(\A)}{\int}\sum_{h_{r-2}\in Q_{r-3}(F)\backslash H_{r-2}(F)}
W_{r-1}(h_{r-2}g)\\ &\qquad\qquad
\underset{[N_{r-1}]}{\int}\Theta(\s(h_{r-2}n_{r-1}h_{r-2}^{-1}h_{r-2})\kappa(g))
\psi_N(n_{r-1})^{-1}\,dn_{r-1}\\
&\qquad\qquad
\underset{[N_{r-2}]}{\int}f^s(\s(n_{r-2}h_{r-2})\kappa(g);e)\psi_N(n_{r-2})^{-1}\,dn_{r-2}dg,
\end{align*} 
where for the second equality  we used
$W_{r-1}(h_{r-2}g)=W_{r-1}(g)$ by Lemma \ref{L:partial_Whittaker} and
the automorphy of $\Theta$.

Now by the change of variable $h_{r-2}n_{r-1}h_{r-2}^{-1}\mapsto
n_{r-1}$ for the integral for $\Theta$, the zeta integral becomes
\begin{align*} 
&\underset{Z(\A)Q_{r-2}(F)N_{r-1}(\A)\backslash
\GLtwo_r(\A)}{\int}\sum_{h_{r-2}\in Q_{r-3}(F)\backslash
H_{r-2}(F)}W_{r-2}(h_{r-2}g)\\
&\underset{[N_{r-1}]}{\int}\Theta(\s(n_{r-1}h_{r-2})\kappa(g))\psi_N(n_{r-1})^{-1}\,dn_{r-1}
\underset{[N_{r-2}]}{\int}f^s(\s(n_{r-2}h_{r-2})\kappa(g); e)\psi_N(n_{r-2})^{-1}\,dn_{r-2}dg
\end{align*} 
by using $\psi_N(h_{r-2}^{-1}n_{r-1}h_{r-2})=\psi_N(n_{r-1})$.

Then one can collapse the sum with the outer integral and obtain
\begin{align*} 
Z(\phi,\Theta,f^s)&=\int W_{r-1}(g)
\underset{[N_{r-1}]}{\int}\Theta(\s(n_{r-1})\kappa(g))\psi_N(n_{r-1})^{-1}\,dn_{r-1}\\
&\qquad\qquad\underset{[N_{r-2}]}{\int}f^s(\s(n_{r-2})
\kappa(g);e)\psi_N(n_{r-2})^{-1}\,dn_{r-2}dg,
\end{align*} 
where the outermost integral is over
\[ 
Z(\A)N_{r-2}(F)Q_{r-3}(F)N_{r-1}(\A)\backslash \GLtwo_r(\A).
\]

By applying the second part of  Lemma \ref{L:key_lemma} with $m=1$, one obtains
\begin{align*} 
Z(\phi,\Theta,f^s)
&={\int}W_{r-1}(g)\underset{[N_{r-2}]}{\int}\underset{[N_{r-1}]}{\int}
\Theta(\s(n_{r-1}n_{r-2})\kappa(g))\psi_N(n_{r-1})^{-1}\,dn_{r-1}dn_{r-2}\\
&\qquad\qquad\underset{[N_{r-2}]}{\int}f^s(\s(n_{r-2})\kappa(g);e)
\psi_N(n_{r-2})^{-1}\,dn_{r-2}dg.
\end{align*} 

By
\begin{align*} 
\int_{Z(\A)N_{r-2}(F)Q_{r-3}(F)N_{r-1}(\A)\backslash\GLtwo_r(\A)}
&=\int_{Z(\A)Q_{r-3}(F)N_{r-1}(\A)N_{r-2}(\A)\backslash\GLtwo_r(\A)}
\int_{N_{r-2}(F)\backslash N_{r-2}(\A)}
\end{align*} 
together with the $\s-\kappa$ trick, one obtains
\begin{align*} 
Z(\phi,\Theta,f^s)
&={\int}\underset{[N_{r-2}]}{\int}W_{r-1}(n_{r-2}'g)
\underset{[N_{r-2}]}{\int}\underset{[N_{r-1}]}{\int}\Theta(\s(n_{r-1}n_{r-2}n_{r-2}')\kappa(g))
\psi_N(n_{r-1})^{-1}\,dn_{r-1}dn_{r-2}\\
&\qquad\qquad\underset{[N_{r-2}]}{\int}f^s(\s(n_{r-2}n_{r-2}')\kappa(g);e)
\psi_N(n_{r-2})^{-1}\,dn_{r-2}dn_{r-2}'dg,
\end{align*} 
where the outermost integral is over
\[ 
Z(\A)Q_{r-3}(F)N_{r-1}(\A)N_{r-2}(\A)\backslash \GLtwo_r(\A).
\] 
The variable $n_{r-2}'$ inside $\Theta$ goes away by the invariance of the
measure $dn_{r-2}$. By the change of variable
$n_{r-2}n_{r-2}'\mapsto n_{r-2}$ inside $f^s$, the character
$\psi_N(n_{r-2}')$ comes out, and one obtains
\begin{align}\label{E:zeta2} 
Z(\phi,\Theta,f^s)\notag
&={\int}\underset{[N_{r-2}]}{\int}W_{r-1}(n_{r-2}'g)\psi(n_{r-2}')dn_{r-2}'\\
\notag
&\qquad\qquad\underset{[N_{r-2}]}{\int}\underset{[N_{r-1}]}{\int}
\Theta(\s(n_{r-1}n_{r-2})\kappa(g))\psi_N(n_{r-1})^{-1}\,dn_{r-1}dn_{r-2}\\
\notag
&\qquad\quad\qquad\underset{[N_{r-2}]}{\int}f^s(n_{r-2}g)(e)
\psi_N(n_{r-2})^{-1}\,dn_{r-2}dg\\
&={\int}W_{r-2}(g)\underset{[N_{r-2}]}{\int}\underset{[N_{r-1}]}{\int}
\Theta(\s(n_{r-1}n_{r-2})\kappa(g))\psi_N(n_{r-1})^{-1}\,dn_{r-1}dn_{r-2}\\
&\qquad\qquad\underset{[N_{r-2}]}{\int}f^s(\s(n_{r-2})\kappa(g);e)\psi_N(n_{r-2})^{-1}\,dn_{r-2}dg,\notag
\end{align} 
where the outermost integrals are over
\[ Z(\A)Q_{r-3}(F)N_{r-1}(\A)N_{r-2}(\A)\backslash \GLtwo_r(\A).
\] \quad\\

\noindent\textbf{\underline{Unfolding Step 3}}\\

The third step starts with computing the Fourier expansion of $\Theta$
along the ``$r-2^{\text{nd}}$ column" $N_{r-3}$, \ie consider the
Fourier expansion of the function
\[
n_{r-3}\mapsto\underset{[N_{r-2}]}{\int}\underset{[N_{r-1}]}{\int}
\Theta(\s(n_{r-1}n_{r-2}n_{r-3})\kappa(g))\psi_N(n_{r-1})^{-1}\,dn_{r-1}dn_{r-2}
\] 
on $[N_{r-3}]$. Again by Lemma \ref{L:Fourier} we have
\begin{align*}
&\underset{[N_{r-2}]}{\int}\underset{[N_{r-1}]}{\int}
\Theta(\s(n_{r-1}n_{r-2})\kappa(g))\psi_N(n_{r-1})^{-1}\,dn_{r-1}dn_{r-2}
=(\text{zero orbit})\\ 
&+\sum_{h_{r-3}\in Q_{r-4}(F)\backslash
H_{r-3}(F)} \underset{[N_{r-2}]}{\int}\underset{[N_{r-1}]}{\int}
\Theta(\s(n_{r-1}n_{r-2}n_{r-3}h_{r-3})\kappa(g))
\psi_N(n_{r-1}n_{r-3})^{-1}\,dn_{r-1}dn_{r-2}dn_{r-3}.
\end{align*} 
By substituting this in (\ref{E:zeta2}), one can show
that, first of all, the zero orbit vanishes thanks to the cuspidality
of $\phi$, and second of all, the sum can be collapsed with the
outermost integral by using the $\s-\kappa$ tick and the change of variable
$h_{r-3}n_{r-2}h_{r-3}^{-1}\mapsto n_{r-2}$ for the integral of $f^s$
along with $\psi_N(h_{r-3}^{-1}n_{r-2}h_{r-3})=\psi_N(n_{r-2})$ and
Lemma \ref{L:partial_Whittaker}. The computations are essentially the
same as the previous steps, and left to the reader.

By applying Lemma
\ref{L:key_lemma} to $f^s$, one obtains
\begin{align*} 
Z(\phi,\Theta,f^s)&={\int}W_{r-3}(g)\\
&\underset{[N_{r-3}]}{\int}\underset{[N_{r-2}]}{\int}
\underset{[N_{r-1}]}{\int}\Theta(\s(n_{r-1}n_{r-2}n_{r-3})\kappa(g))
\psi_N(n_{r-1}n_{r-3})^{-1}\,dn_{r-1}dn_{r-2}dn_{r-3}\\
&\underset{[N_{r-3}]}{\int}\underset{[N_{r-2}]}{\int}
f^s(\s(n_{r-2}n_{r-3})\kappa(g);e)\psi_N(n_{r-2})^{-1}\,dn_{r-2}dn_{r-3}\;dg,
\end{align*} 
where the outermost integral is over
\[ 
Z(\A)Q_{r-4}(F)N_{r-1}(\A)N_{r-2}(\A)N_{r-3}(\A)\backslash
\GLtwo_r(\A).
\]

\quad\\

\noindent\textbf{\underline{Unfolding Step 4 and further}}\\

We repeat this process. Namely the next step (Step 4) is to apply
the Fourier expansion formula (Lemma \ref{L:Fourier}) to the function
\[ 
n_{r-4}\mapsto
\underset{[N_{r-3}]}{\int}\underset{[N_{r-2}]}{\int}
f^s(\s(n_{r-2}n_{r-3}n_{r-4})\kappa(g);e)\psi_N(n_{r-2})^{-1}\,dn_{r-2}dn_{r-3}
\] 
on $[N_{r-4}]$ and one sees that the zero orbit vanishes by the
cuspidality of $\phi$, and collapse the sum by using Lemma
\ref{L:partial_Whittaker}. Then apply Lemma \ref{L:key_lemma} to
$\Theta$, which gives
\begin{align*} 
&Z(\phi,\Theta,f^s)\\
=&{\int}W_{r-4}(g)\\
&\underset{[N_{r-4}]}{\int}\underset{[N_{r-3}]}{\int}\underset{[N_{r-2}]}{\int}
\underset{[N_{r-1}]}{\int}\Theta(\s(n_{r-1}n_{r-2}n_{r-3}n_{r-4})\kappa(g))
\psi_N(n_{r-1}n_{r-3})^{-1}\,dn_{r-1}dn_{r-2}dn_{r-3}dn_{r-4}\\
&\underset{[N_{r-4}]}{\int}\underset{[N_{r-3}]}{\int}\underset{[N_{r-2}]}{\int}
f^s(\s(n_{r-2}n_{r-3}n_{r-4})\kappa(g);e)\psi_N(n_{r-2}n_{r-4})^{-1}\,dn_{r-2}dn_{r-3}dn_{r-4}\;dg,
\end{align*} 
where the outermost integral is over
\[
Z(\A)Q_{r-5}(F)N_{r-1}(\A)N_{r-2}(\A)N_{r-3}(\A)N_{r-4}(\A)\backslash
\GLtwo_r(\A).
\] 
For the next step (Step 5) one needs to compute the Fourier expansion for
$\Theta$ along $N_{r-5}$ using Lemma \ref{L:Fourier} (the zero orbit
goes away by the cuspidality of $\phi$), then collapse the
sum using Lemma \ref{L:partial_Whittaker}, and then apply Lemma
\ref{L:key_lemma} to $f^s$. Then the next step (Step 6) is to switch
the roles of $\Theta$ and $f^s$ and use those three lemmas, Lemma
\ref{L:Fourier}, \ref{L:partial_Whittaker} and \ref{L:key_lemma}, in
this order, and then proceed to the next step, and so on.

\quad\\

\noindent\textbf{\underline{Unfolding Final Step}}\\

After finishing step $r-2$, which is done by computing the Fourier
expansion of $f^s$, one obtains
\begin{align*} 
&Z(\phi,\Theta,f^s)\\=&{\int}W_2(g)\\
&\underset{[N_{2}]}{\int}\underset{[N_{3}]}{\int}\cdots\underset{[N_{r-1}]}{\int}
\Theta(\s(n_{r-1}n_{r-2}\cdots n_{2})\kappa(g))\psi_N(n_{r-1}n_{r-3}\cdots
n_{3})^{-1}\,dn_{r-1}dn_{r-2}\cdots dn_{2}\\
&\underset{[N_{2}]}{\int}\underset{[N_{3}]}{\int}\cdots\underset{[N_{r-2}]}{\int}
f^s(\s(n_{r-2}n_{r-3}\cdots n_{2})\kappa(g);e)\psi_N(n_{r-2}n_{r-4}\cdots
n_{2})^{-1}\,dn_{r-2}dn_{r-3}\cdots dn_{2}\;dg,
\end{align*} 
where the outermost integral is over
\[ 
Z(\A)Q_{1}(F)N_{r-1}(\A)N_{r-2}(\A)\cdots N_{3}(\A)N_{2}(\A)\backslash
\GLtwo_r(\A).
\] 

The final step (Step $r-1$) does not work out as before because
the key Lemma \ref{L:Fourier} does not hold for $m=1$. Namely, for
$m=1$, though $H_1$ acts on $\widehat{[N_1]}$ as before, the number of
orbits is not $2$ but rather the nonzero orbits are indexed by
$(F^\times)^2\backslash F^\times$, and indeed
\[ 
\widehat{[N_1]}=(\text{zero orbit})+\sum_{a\in
(F^\times)^2\backslash F^\times}H_1\psi_N^a,
\] 
where the stabilizer of each $\psi_N^a$ in $H_1$ is $Q_0$. But
everything else is the same as the previous steps, and we obtain
\begin{align*}
&Z(\phi,\Theta,f^s)={\int}W_1(g)\notag\\ 
&\sum_{a\in(F^\times)^2\backslash F^\times}\underset{[N_{1}]}{\int}
\underset{[N_{2}]}{\int}\cdots\underset{[N_{r-1}]}{\int}
\Theta(\s(n_{r-1}n_{r-2}\cdots n_{1})\kappa(g))\psi_N^a(n_{r-1}n_{r-3}\cdots
n_{2})^{-1}\,dn_{r-1}dn_{r-2}\cdots dn_{1}\\
&\notag\underset{[N_{1}]}{\int}\underset{[N_{2}]}{\int}\cdots\underset{[N_{r-2}]}{\int}
f^s(\s(n_{r-2}n_{r-3}\cdots n_{1})\kappa(g);e)\psi_N(n_{r-2}n_{r-4}\cdots
n_{2})^{-1}\,dn_{r-2}dn_{r-3}\cdots dn_{1}\;dg,
\end{align*} 
where the outermost integral is over
\[ 
Z(\A)N_{r-1}(\A)N_{r-2}(\A)\cdots N_{1}(\A)\backslash
\GLtwo_r(\A),
\] 
which is the same as
\[
Z(\A)N(\A)\backslash \GLtwo_r(\A)
\] 
because $N_{r-1}(\A)N_{r-2}(\A)\cdots N_{1}(\A)=N(\A)$.

\quad\\

\noindent\textbf{\underline{Almost Euler product}}\\

Now we are ready to obtain the (almost) Euler product from this last
expression. But as we have noted before, we are not able to
obtain the full Euler product. This is due to
the lack of the uniqueness result for the semi-Whittaker functional at
the archimedean places. Although such uniqueness result might hold at
the archimedean places, at this moment the author does not know how to
prove it. Hence the best we can do is to obtain the ``almost Euler
product", or the Euler product at the finite places.

First notice that $W_1(g)$ is the usual Whittaker coefficient with
respect to $\psi^{-1}$, so let us write
\[ 
W_1(g)=W_\phi^{\psi^{-1}}(g)=W(g),
\] 
where again we ignore the dependence of $W(g)$ on $\phi$ and
$\psi$. Also by following \cite{BG}, we define
\begin{align*}
Q^a(\kappa(g))&=\underset{[N_1]}{\int}\cdots\underset{[N_{r-1}]}{\int}
\Theta(\s(n_{r-1}n_{r-2}n_{r-3}\cdots n_1)\kappa( g))\\
&\qquad\qquad\qquad\qquad\qquad\qquad\qquad\qquad
\psi^a_N(n_{r-1}n_{r-3}n_{r-5}\cdots)dn_{r-1}\cdots dn_1\\
\end{align*} 
and
\begin{align*}
R^s(\kappa(g))&=\underset{[N_1]}{\int}\cdots\underset{[N_{r-2}]}{\int}
f^s(\s(n_{r-2}n_{r-3}\cdots n_1)\kappa( g);e)\\
&\qquad\qquad\qquad\qquad\qquad\qquad\qquad\qquad\qquad
\psi_N(n_{r-2}n_{r-4}\cdots)dn_{r-2}\cdots dn_1.
\end{align*} 
With this notation, the last formula we obtained for $Z(\phi,\Theta,f^s)$ is written as
\[ 
Z(\phi,\Theta,f^s)=\underset{Z(\A)N(\A)\backslash
\GLtwo(\A)}{\int}W(g)\left(\sum_{a\in (F^\times)^2\backslash
F^\times}Q^a(\kappa(g))\right)R^s(\kappa(g))dg.
\]

We need to take care of the sum $\sum_{a\in
(F^\times)^2\backslash F^\times}$. First for each
fixed $g\in\GL_r(\A)$ consider the map
$\vartheta_{\chi\omega^{-2}}^\psi\rightarrow\C$ defined by
\[ 
\Theta\mapsto Q^a(\kappa(g)).
\] 
This map is non-zero, because it is (a scalar multiple of) the
composite of the constant term map
$\vartheta_{\chi\omega^{-2}}^\psi\rightarrow \omega_{\chi}^\psi\otimes\delta_P^{1/4}$
along the $(2,\dots,2)$-parabolic $P$ 
with the $\psi_{(a,1,\dots,1)}$-Whittaker functional of
$\omega_{\chi}^\psi\otimes\delta_P^{1/4}$. (See Proposition
\ref{P:global_periodicity3} for the constant term and
(\ref{E:additive_a}) for the notation $\psi_{(a,1,\dots,1)}$.) With
this said, one can see that
Proposition \ref{P:tensor_product_Weil_global} implies
\begin{Lem}\label{L:Q(g)_non_vanishing}
The map $\Theta\mapsto  Q^a(\kappa(g))$ is not identically zero if and
only if $a\equiv 1\mod (F^\times)^2$.
\end{Lem}
This gives
\[
 Z(\phi,\Theta,f^s)=\underset{Z(\A)N(\A)\backslash
\GLtwo(\A)}{\int}W(g)Q(\kappa(g))R^s(\kappa(g))dg,
\] 
where we wrote $Q^1(\kappa(g))=Q(\kappa(g))$. This is precisely the analogue of
\cite[(3.5)]{BG}.

Now it would be ideal if we could show both $Q$ and
$R_s$ decompose into products of local components like
$Q(\kappa_v(g))=\prod_{v} Q_v(\kappa_v(g_v))$ and $R^s=\prod_{v}
R^s_v(\kappa_v(g_v))$
and hence by choosing $\phi$ so that the Whittaker function $W$
decomposes into a product $\prod_v W_v$, we could obtain the Euler
product. However, to achieve this, one needs the uniqueness of the semi-Whittaker
functional on the local exceptional representation for both
archimedean and non-archimedean cases. For the non-archimedean case,
the uniqueness of the semi-Whittaker functional follows from the
periodicity of the Jacquet module of the exceptional
representation (Proposition \ref{P:semi-Whittaker}), which seems to
be available only for the non-archimedean case, and the author
does not know if such uniqueness is available for the archimedean
case. (See Remark \ref{R:semi_Whittaker}). 
Because of this issue, we need to compromise with the almost Euler
product, which is, nonetheless, enough for proving our main
theorem.

First let
\[
 L:\vartheta_\chi^\psi\rightarrow\C,\quad\Theta\mapsto Q(\kappa(e)),
\] 
where $e\in\GLtwo_r(\A)$ is the identity element. This is a global
semi-Whittaker functional. Note that $L$ is not identically zero by
Lemma \ref{L:Q(g)_non_vanishing}.

Next let us define $\GLt_r(F_\infty)$ to be the image of the map
\[
\prod_{v|\infty}\GLt(F_v)\rightarrow\GLt_r(\A), \quad 
\prod_{v|\infty}(g_v,\xi_v)\mapsto (\prod_{v|\infty}g_v,\prod_{v|\infty}\xi_v),
\]
which may be called ``the archimedean component'' of
$\GLt_r(\A)$. Then we can write
$\theta=\theta_{\infty}\otimest\left(\otimest_{v<\infty}'\theta_v\right)$,
where
\[ 
\theta_{\infty}=\underset{v|\infty}{\otimest}\theta_v
\] 
is the metaplectic tensor product of $\theta_v$ for all archimedean
$v$, which is a representation of $\GLt_r(F_\infty)$.
And we write each simple tensor in
$\theta_{\infty}\otimest\left(\otimest_v'\theta_v\right)$ as
$x_\infty\otimes\left(\otimes'_vx_v\right)$, where
$x_\infty\in\theta_\infty$ and
$\otimes'_vx_v\in\otimest_v'\theta_v$. (Since the space of
restricted metaplectic tensor
product is the same as the usual restricted tensor product, we use
the notation $\otimes$ rather than $\otimest$ to denote each vector.) 
Let us fix the vector
$x_\infty^\circ\otimes\left(\otimes'_vx_v^\circ\right)\in\theta$ such
that at unramified $v$, $x_v^\circ$ is the spherical vector used to
define the tensor product $\otimest'_{v<\infty}\theta_v$ and
$L(x_\infty^\circ\otimes\left(\otimes'_vx_v^\circ\right))=1$. (Such
vector certainly exists.) We define
$L_\infty:\theta_\infty\rightarrow\C$ by
\[
L_\infty(x_\infty)=L(x_\infty\otimes\left(\otimes'_vx_v^\circ\right)).
\] 
One can show that
\begin{Prop} 
For $\Theta=x_\infty\otimes\left(\otimes'_v
x_v\right)\in\theta$, one has
\[ 
Q(\kappa(g))=Q_\infty(\kappa_\infty(g_\infty))\prod_{v<\infty} Q_v(\kappa_v(g_v)),
\] 
where
\[ 
Q_\infty(\kappa_\infty(g_\infty))=L_\infty(\theta_\infty(\kappa_\infty(g_\infty))x_\infty)
\] 
and for $v<\infty$
\[ 
Q_v(\kappa_v(g_v))=L_v(\theta_v(\kappa_v(g_v))x_v)
\] 
where $L_v$ is (a scalar multiple of) the semi-Whittaker functional
on $\theta_v$ such that $L_v(x_v^\circ)=1$ for almost all $v$. Also
note that for each $v$ (non-archimedean or not),
$\kappa_v:\GL_r(F_v)\rightarrow\GLt_r(F_v)$ is the set theoretic
section $g_v\mapsto (g_v,1)$, and
$\kappa_\infty:\prod_{v|\infty}\GL_r(F_v)\rightarrow\GLt_r(F_\infty)$
is given by
$\kappa_\infty(\prod_{v|\infty}(g_v))=(\prod_{v|\infty}(g_v), 1)$.
\end{Prop}
\begin{proof} 
The proof is the same as the usual proof that a
Whittaker-Fourier coefficient can be decomposed as the Euler
product. (See for example \cite[Theorem 3.5.2]{Bump}). We will repeat
the essential points here. First of all, we may assume $g=e$ because
that would simply replace $\Theta$ by $g\cdot\Theta$. Let $S$ be the
finite set of finite places at which $x_v\neq x_v^\circ$, so
$L_v(x_v)=Q_v(\kappa_v(e))=1$ if $v\notin S$. The proof is by induction on the
cardinality of $S$. Namely assume $S$ is empty. Then $x_v=x_v^\circ$ for
all finite $v$. Then $L(x_\infty(\otimes_v'
x_v)))=L_\infty(x_\infty)=Q_\infty(\kappa_\infty(e))=Q_\infty(e)\prod_{v<\infty}
Q_v(\kappa_v(e))$. This is the base step of induction. 

Assume the statement holds of all vectors $x_\infty\otimes(\otimes_v'
x_v)$ whose $S$ has cardinality for some $n$. Now assume $y_\infty\otimes(\otimes_v'
y_v)$ is such that the corresponding $S$ has cardinality $n+1$. Let $w$
be a place where $y_w\neq x_w^\circ$. Consider the map
\[
\theta_w\rightarrow\C,\quad
y'_w\mapsto L(y_\infty\otimes y'_w\otimes\left({\bigotimes_{v\neq w}}' y_v\right)).
\]
This is a semi-Whittaker functional for ${\vartheta_{\chi}^\psi}_w$.
By the uniqueness of the local semi-Whittaker functional
(Proposition \ref{P:semi-Whittaker}), this is equal
$cL_w(y'_w)$ for some scalar. Let $y'_w=x_w^\circ$, so that by the
induction hypothesis,
\[
cL_w(x_w^\circ)=L_\infty(y_\infty)L_w(x_w^\circ)\prod_{v\neq
  w}L_v(y_v).
\] 
But $L_w(x_w^\circ)=1$, which gives
$c=L_\infty(y_\infty)\prod_{v\neq w}L_v(y'_v)$. Thus we have 
\[
L(y_\infty\otimes y'_w\otimes\left({\bigotimes_{v\neq w}}' y_v\right))
=L_\infty(y_\infty)L_w(y'_w)\prod_{v\neq w}L_v(y_v)
\]
for any $y_w'$. By letting $y_w'=y_w$, the induction is complete.
\end{proof} 

Similarly we can obtain the decomposition
\[ 
R^s(\kappa(g))=R^s_\infty(\kappa_\infty(g_\infty))\prod_{v<\infty} R^s_v(\kappa_v(g_v))
\] 
for a decomposable
$f^s=f^s_\infty\otimes\left(\otimes'_v{f^s}_v\right)$ by defining
$L'_\infty:\theta'\rightarrow\C$ as before and setting
\[ 
R^s_\infty(\kappa_\infty(g_\infty))=L'_\infty(\theta_\infty'(\kappa_\infty(g_\infty))f^s_\infty)
\] 
and for $v<\infty$
\[ 
R^s_v(\kappa(g_v))=L'_v(\theta'_v(\kappa(g_v))f^s_v)
\] 
where $L'_v:\theta'_v\rightarrow\C$ is (a scalar multiple of) the
semi-Whittaker functional on $\theta_v'$. Hence we obtain the almost
Euler product of the zeta integral
\[ 
Z(W,Q,R^s)=Z_\infty(W_\infty,Q_\infty,R^s_\infty)\prod_{v<\infty}
Z_v(W_v,Q_v,R^s_v),
\] 
where
\[ 
Z_\infty(W_\infty,Q_\infty,R^s_\infty)=\underset{Z(F_\infty)
N(F_\infty)\backslash \GLtwo_r(F_\infty)}{\int}
W_\infty(g_\infty)Q_\infty(\kappa_\infty(g_\infty))R_\infty^s(\kappa_\infty(g_\infty))\;dg_\infty
\] 
and for $v<\infty$
\[ 
Z_v(W_v,Q_v,R^s_v)=\underset{Z(F_v) N(F_v)\backslash
\GLtwo_r(F_v)}{\int} W_v(g_v)Q_v(\kappa_v(g_v))R_v^s(\kappa_v(g_v))\;dg_v.
\] 
Here note that $F_\infty=\prod_{v|\infty}F_v$, which is a product
of copies of $\R$ and/or $\C$.

\quad\\

\noindent\textbf{\underline{Unramified factor}}\\

We will compute the unramified factor here. For this we need the
following ``Iwasawa decomposition'' of $\GL_r^{(2)}$.
\begin{Lem}\label{L:Iwasawa}
Assume $F$ is a non-archimedean local field and $P$ is any parabolic
subgroup of $\GL_r(F)$. We have the decomposition
\[
\GL_r^{(2)}(F)=P(F)^\#K^\#,
\]
where $P^\#(F)=P(F)\cap\GL_r^{(2)}(F)$ and $K^\#=\GL_r(\OF)\cap\GL_r^{(2)}(F)$.

For a measurable function $f$ on $G=\GL_r^{(2)}(F)$, we have
\[
\int_{G}f(g)\;dg=\int_{P^\#}\int_{K^\#}f(pk)\;dp\,dk
\]
where $dp$ is the left Haar measure on $P^\#$ and $dk$ is the Haar measure
on $K^\#$.
\end{Lem}
\begin{proof}
By the usual Iwasawa decomposition of $\GL_r$, each element $g\in
\GL_r^{(2)}(F)$ is written as $g=pk$ for $p\in P(F)$ and $k\in K$ such
that $\det(p)\det(k)\in(F^\times)^2$. We may assume $\det(p)=\varpi^n$ where
$\varpi$ is a uniformizer of $F$ and $n\in\Z$; For if $\det(p)=\varpi^n u$
for some $u\in\OF^\times$, let
$k_1$ be an element in $K\cap P$ with $\det(k_1)=u^{-1}$, for example
$k_1=\left(\begin{smallmatrix}u^{-1}&&&\\ &1&&\\ &&\ddots&\\
    &&&1\end{smallmatrix}\right)$. 
Then since $pk=(pk_1)(k_1^{-1}k)$ and $pk_1\in P$, we may simply
replace $p$ by $pk_1$.

Thus $\det(pk)=\varpi^n\det(k)\in (F^\times)^2$. But
$\det(k)\in\OF^\times$. Hence we must have $\det(p)\in
(F^\times)^2$. (Indeed this implies that $n$ is even.) 
Then $\det(k)\in (F^\times)^2$ as well

The decomposition of the measure is \cite[Proposition
2.1.5 (ii)]{Bump}. (The assumptions of \cite[Proposition
2.1.5 (ii)]{Bump} are satisfied by $P^\#$ and $K^\#$.)
\end{proof}

\begin{Prop} At each unramified place $v$,
\[ Z_v(W_v,Q_v,R^s_v)=L(2s-\frac{1}{2},\pi_v,
Sym^2\otimes\chi_v)L(r(2s-\frac{1}{2}),\chi_v^r\omega_v^2)^{-1}.
\]
\end{Prop}
\begin{proof} The computation is almost identical to \cite[Theorem
4.1]{BG}, and
hence we only give the key points. Also we omit the subscript $v$ in
our notation and simply write $\GL_{r}^{(2)}=\GL_{r}^{(2)}(F_v)$,
$N=N_B(F_v)$, $Z=Z(F_v)$, $T=T(F_v)$, $B=B(F_v)$ and $K=\GL_r({\OF}_v)$. 

We will work on the integral
\[
\underset{ZN\backslash \GL_{r}^{(2)}}{\int}W(g)Q(\kappa(g)){R^s}(\kappa(g)) dg,
\]
where all the data are unramified. By the above lemma, this is written
as
\[ 
\underset{ZN\backslash B^\#}{\int}\underset{K^\#}{\int}
W(bk)Q(\kappa(bk)){R^s}(\kappa(bk)) dkdb,
\] 
where $B^\#$ is as in the above lemma with $P=B$.
By the $\s-\kappa$ trick (or strictly speaking it should be called
``$\kappa-\kappa$'' trick in this case), this is written as
\[ 
\underset{ZN\backslash B^\#}{\int}\underset{K^\#}{\int}
W(bk)Q(\kappa(b)\kappa(k)){R^s}(\kappa(b)\kappa(k)) dkdb.
\] 
By the $K$ invariance of the integrand, we have
\[ 
\underset{ZN\backslash B^\#}{\int}W(b)Q(\kappa(b)){R^s}(\kappa(b)) db.
\] 
Since the integrand is left $N$ invariant, this is
written as
\[ 
\underset{Z\backslash T^\#}{\int}W(t)
Q(\kappa(t))R^s(\kappa(t))\delta_B(t)^{-1} dt,
\] 
where $T^\#=T\cap B^\#$ and $\delta_B$ is the modular character of
the Borel subgroup $B$. (Once again, one also need the $\s-\kappa$
trick for this formulation.) For each
\[ 
\lambda=(\lambda_1,\dots,\lambda_{r})\in \Z^r,
\] 
we write
\[
t_\lambda=\begin{pmatrix}
\varpi^{\lambda_1}&&\\&\ddots&\\&&\varpi^{\lambda_{r}}\end{pmatrix}.
\] 
Then the integral is equal to
\[ 
\sum_{\lambda\in\Z^{r}, t_\lambda\in
T^\#}W(t_\lambda)Q(\kappa(t_\lambda))R^s(\kappa(t_\lambda))\delta_B(t_\lambda)^{-1},
\] 
where $\lambda$ runs through the elements of the form
$(\lambda_1,\dots,\lambda_{r-1},0)$ with
$\sum_{i=1}^{r-1}\lambda_i=\text{ even }$. (Since we mod out by $Z$,
we always have $\lambda_{r}=0$.)

We have
\[
R^s(\kappa(t_\lambda))=Q'(\kappa(t_\lambda))\delta_Q(t_\lambda),
\]
where $Q'$ is the semi-Whittaker functional for the inducing
representation $\vartheta_{\omega,\omega^{-1}\chi^{-q}}^\psi$, and
$\delta_Q$ is the modular character for the parabolic
$Q$. (Unfortunately we have two different $Q$ here, but we assume it should
not create any confusion.)

Let $\alpha=\{\alpha_1,\dots,\alpha_r\}$ be the Satake parameter of
the $v$ component $\pi_v$ of our cuspidal representation $\pi$. 
By Shintani's formula (\cite{Shintani}),
\[
W(t_\lambda)=\begin{cases}\delta_B(t_\lambda)^{1/2}s_\lambda(\alpha),\text{
if $\lambda_1\geq\cdots\geq\lambda_{r-1}\geq 0$};\\ 0,\text{
otherwise,}\end{cases}
\] 
where $s_\lambda$ is the symmetric function of $r$ variables as defined in
\cite[Section I.3]{Mac}, and $s_\lambda(\alpha)$ is the value of the
function evaluated at the Satake parameter.

Following \cite{BG}, we call $\lambda$ even if all the components
$\lambda_i$ are even. Since $\chi^{1/2}$ certainly exists in the
unramified situation, which we fix, one can see that
\[
Q(\kappa(t_\lambda))=\begin{cases}\delta_B^{1/4}(t_\lambda)\chi^{1/2}\omega^{-1}
(\det(t_\lambda)),\quad\text{$\lambda$ is even};\\ 0,\quad\text{ otherwise}\end{cases}.
\] 
One can also see
\[
Q'(\kappa(t_\lambda))=\begin{cases}\delta_{B'}^{1/4}(t_\lambda)\omega(\det(t_\lambda)),
\quad\text{$\lambda$ is even};\\ 
0,\quad \text{ otherwise,}\end{cases}
\] 
where $\delta_{B'}$ is the modulus character of the Borel subgroup
$B'$ of $\GL_{r-1}$ viewed as a subgroup of $\GL_{r}$ with the
embedding
$h\mapsto\left(\begin{smallmatrix}h&\\&1\end{smallmatrix}\right)$. By
multiplying all those, one can see that the local zeta integral is
equal to
\[ 
\sum_{\substack{ \text{even } \lambda\in\Z^{r}\\
\lambda_1\geq\lambda_2\geq\dots\geq\lambda_{r-1}\geq0,\; \lambda_r=0}}
s_{\lambda}(\alpha)\delta_P(t_{\lambda})^{s-\frac{1}{4}}\chi^{1/2}(\det(t_{\lambda})),
\] 
which is precisely the twisted analogue of \cite[(4.7)]{BG}. Hence
the computation in the proof of \cite[Theorem 4.1]{BG} can be directly
applied to our integral, which yields the proposition. Namely, as in
p.171 of \cite{BG}, we have
\[ \prod_{1\leq i\leq j\leq r}(1-\alpha_i\alpha_jX)^{-1}
=\left[\sum_{\substack{ \text{even } \lambda\in\Z^{r-1}\\
\lambda_1\geq\lambda_2\geq\dots\geq\lambda_{r-1}\geq0}}
s_{\lambda}(\alpha)X^{\left(\sum_{i=1}^{r-1}\lambda_i\right)/2}\right]
\times (1-\omega(\varpi)^2 X^r)^{-1}.
\] (Here notice that in the exponent of $X$ in the first factor of the
corresponding formula in \cite{BG}, there is a typo.) By taking
$X=\chi(\varpi)q^{-2s+1/2}$, we obtain our proposition.
\end{proof}

\quad

%%%%%%%%%%%%%%%%%%%%%%%%%%%%%%%%%%%%%%%%%%%%%%%%%%%%%%%%%%%%%%%%%%%

\section{\bf The Rankin-Selberg integrals for the case $r=2q+1$}

%%%%%%%%%%%%%%%%%%%%%%%%%%%%%%%%%%%%%%%%%%%%%%%%%%%%%%%%%%%%%%%%%%%

We consider
\[ 
r=2q+1=\text{odd}.
\] 
Note that for the $r=2q+1$ case there is no issue raised by Kable
(\cite{Kable}) for the Rankin-Selberg integral of Bump and Ginzburg. But
in order to incorporate the character twist into the Bump-Ginzburg
integral, we need to choose
\[ 
\theta=\theta_{\omega^{-1}},
\] 
where $\theta_{\omega^{-1}}$ is the global non-twisted exceptional
representation of $\GLt_r(\A)$ with determinantal character
$\omega^{-1}$, and
\[ 
\theta'=\vartheta'_{\chi\omega^2, \chi^{-q}}
\] 
for the exceptional representation of
$\GLt_{r-1}(\A)\timest\GLt_1(\A)\subseteq\GLt_r(\A)$ associated with
$\chi\omega^2$ and $\chi^{-q}$.  Notice that the
central character of $\theta$ is
\[ 
(1,\xi)\s(z)\mapsto \xi\omega^{-2q-1}(a)\mu_\psi(a)^q
\] 
by (\ref{E:central_character1}), and the central character of $\theta'$ is
\[ 
(1,\xi)\s(z)\mapsto \xi\omega^{2q}(a)\mu_\psi(a)^{-q}
\]
by (\ref{E:central_character4}), where $z=aI_{2q+1}$.

Then for $\Theta\in\theta$ and
$f^s\in\Ind_{\Pt(\A)}^{\GLt_{r}(\A)}\theta'\otimes\delta_{P}^s$, we
define the global zeta integral as
\[ 
Z(\phi,\Theta,
f^s)=\underset{Z(\A)\GL_r(F)\backslash\GL_r(\A)}{\int}\phi(g)
\Theta(\kappa(g))E(\kappa(g),s,f^s)\;dg,
\] 
where $\phi\in\pi$ is a cusp form and $E(-,s,f^s)$ is the
Eisenstein series as before. Note that the product
$\Theta(\kappa(g))E(\kappa(g),s,f^s)$ is not genuine, on which 
the center $z\in Z(\A)$ acts as the character $\omega^{-1}$, and hence 
the integral is well-defined. By following the computation of \cite{BG}, the global
integral decomposes into the almost Euler product
\[ 
Z(\phi,\Theta,
f^s)=Z(W,Q,R^s)=Z_\infty(W_\infty,Q_\infty,R^s_\infty)\prod_{v<\infty}Z_v(W_v,Q_v,R^s_v),
\] 
where $W, Q$ and $R^s$ and their local components are just as the
$r=2q$ case. Note again that just like the case $r=2q$ because of the
issue on the uniqueness of the semi-Whittaker functional at the
archimedean places, we need to compromise with this almost Euler
product instead of the full Euler product.

We can compute the unramified factor as follows.
\begin{Prop} At each unramified place $v$,
\[ Z_v(W_v,Q_v,R^s_v)=L(2s-\frac{1}{2},\pi_v,
Sym^2\otimes\chi_v)L(r(2s-\frac{1}{2}),\chi_v^r\omega_v^2)^{-1}.
\]
\end{Prop}
\begin{proof} This is even more straightforward modification of
\cite{BG} than the $r=2q$ case. Also see \cite[Theorem 7]{Banks2} for
the case $r=3$.
\end{proof}

\quad

%%%%%%%%%%%%%%%%%%%%%%%%%%%%%%%%%%%%%%%%%%%%%%%%%%%%%%%%%%%%%%%%%%%

\section{\bf The poles of $L^S(s,\pi,Sym^2\otimes\chi)$ }

%%%%%%%%%%%%%%%%%%%%%%%%%%%%%%%%%%%%%%%%%%%%%%%%%%%%%%%%%%%%%%%%%%%

Now we are ready to prove the following main theorem of this paper.
\begin{Thm}\label{T:main} Let $\pi$ be a cuspidal automorphic
representation of $\GL_r(\A)$ with central character $\omega$ and
$\chi$ a unitary Hecke character. Then for each archimedean $v$, there
exists an integer $N_v\geq 0$ such that the product
\[ 
L^S(s,\pi,Sym^2\otimes\chi)\prod_{v|\infty}L_v(rs-r+1,
\chi_v^r\omega_v^2)^{-N_v}
\] 
is holomorphic everywhere except at $s=0$ and $s=1$. Moreover there
is no pole if $\chi^r\omega^2\neq 1$.
\end{Thm}
\begin{proof} 
The proof is a modification of the one given by Bump and
Ginzburg \cite[Theorem 7.5]{BG}. Since the essential points are
already in \cite{BG}, we only give a sketch of the proof for most of
the time. Our Rankin-Selberg integral gives
\begin{align}\label{E:integral_representation} 
&L(r(2s-\frac{1}{2}),
\chi^r\omega^2)Z(\phi,\Theta,f^s)\notag\\
&=L^S(2s-\frac{1}{2},\pi,Sym^2\otimes\chi) L_\infty(r(2s-\frac{1}{2}),
\chi_\infty^r\omega_\infty^2)Z_\infty(W_\infty,Q_\infty,R^s_\infty)\\
&\qquad\prod_{v\in S,\; v<\infty}L_v(r(2s-\frac{1}{2}),
\chi_v^r\omega_v^2)Z_v(W_v,Q_v,R^s_v)\notag
\end{align} 
for the factorizable
$f^s=f_\infty^s\otimes\left(\otimes'f^s_v\right)$. Recall that
\[ 
Z(\phi,\Theta,f^s)={\int}\phi(g)\Theta(\kappa(g))E(\kappa(g),s,f^s)dg,
\] 
where the integral is over $Z(\A)\GLtwo_r(F)\backslash
\GLtwo_r(\A)$ if $r$ is even and $Z(\A)\GL_r(F)\backslash \GL_r(\A)$
if $r$ is odd. Let us define the
normalized Eisenstein series by
\[ 
E^\ast(g,s,f^s):=L^S(r(2s-\frac{1}{2}), \chi^r\omega^2)E(g,s,f^s).
\] 
Let us note the following proposition, whose proof will be given after the
proof of this main theorem.
\begin{Prop}\label{P:poles_Eisenstein} 
Let $f^s$ be a flat
section. Then for each archimedean $v$, there exists an integer
$N_v\geq 0$ such that the product
\[ 
E^\ast(g,s,f^s)\prod_{v|\infty}L_v(r(2s-\frac{1}{2})-r+1,
\chi_v^r\omega_v^2)^{-N_v}
\] 
is entire except that, if $\chi^r\omega^2=1$, it has simple poles
at $s=1/4$ and $s=3/4$.
\end{Prop}

\begin{Rmk} 
Let us note that we are not able to show that the normalized
Eisenstein series $E^\ast(g,s,f^s)$ has the desired analytic
properties, but we need to multiply a kind of compensation factor
$L_v(r(2s-\frac{1}{2})-r+1,\chi_v^r\omega_v^2)^{-N_v}$ at each
archimedean place. This is because of a subtle issue about asymptotic 
expansions of matrix coefficients to be explained later.
\end{Rmk}

We also have

\begin{Prop} 
The local zeta integral $Z_v(W_v,Q_v,R^s_v)$ (resp. the
archimedean $Z_\infty(W_\infty,Q_\infty,R^s_\infty)$) has meromorphic
continuation as a function in $s\in\C$. Moreover, for each fixed
$s=s_0$, one may choose the local data so that $Z_v(W_v,Q_v,R^s_v)$
(resp. $Z_\infty(W_\infty,Q_\infty,R^s_\infty)$) does not have a zero
at $s=s_0$.
\end{Prop}
\begin{proof} 
For the non-archimedean zeta integral, the first part is proven in the
same way as
\cite[Proposition 5.2]{BG} and the second part is as \cite[Theorem
7.2]{BG}. For the archimedean zeta integral, we can apply their arguments to the product
of copies of $\GL_r(\R)$ and/or $\GL_r(\C)$ instead of just one copy
of each.

But since our zeta integrals are not identical to those treated in
\cite{BG}, we repeat the essential points by making clear how the
proofs have to be modified. First for the meromorphic continuation
(\cite[Proposition 5.2]{BG}), there are two key ingredients. One is the
asymptotic expansion of the Whittaker functions (\cite[\S4]{JS}) 
for $W_v$ and the product
$Q_vR_v^s$. (Since $Q_v$ and $R_v^s$ are semi-Whittaker
functionals, in which the alternating entries one above the diagonal
come out via the additive character, the
product $Q_vR_v^s$ is a Whittaker functional.) Note that for
$r=2q$, all the data $W_v, Q_v$ and $R_v^s$ are restrictions to
$\GLtwo_r(F_v)$ of those defined over $\GL_r(F_v)$, and hence there is no
issue for applying this theory. The second ingredient is the Iwasawa
decomposition. For the non-archimedean case, this is Lemma
\ref{L:Iwasawa}. For the archimedean case, we also have the analogous
decomposition. Namely if
$F_v=\C$, then $\GLtwo_r(\C)=\GL_r(\C)$, so there is no issue
here. If $F_v=\R$, then $\GLtwo_r(\R)=\GL_r(\R)^+=\{g\in\GL_r(\R):
\det(g)>0\}$, and we have the Iwasawa decomposition with
$K^\#=\SO(n)$. Using those two ingredients, one can reduce the problem
to meromorphic continuation of a torus integral of a finite sum of a
product of a Schwartz function and a finite function (see \cite[p.178]{BG}),
where by torus we mean $T\cap B^\#$ when $r=2q$. The rest of the
computation is identical.

For the second part of the proposition, which corresponds to
\cite[Theorem 7.2]{BG}, again the key ingredient is the Iwasawa
decomposition. With it, one can reduce the problem to a problem on a
integral over $GL_{r-1}(F_v)$, where $\GL_{r-1}(F_v)$ sits in the
Levi part of the $(r-1, 1)$-parabolic of $\GL_r(F_v)$, and show the
non-vanishing of the integral by induction. For the case $r=2q$, one
can argue in the same way using $\GLtwo_{r-1}(F_v)$.
\end{proof}

Hence by taking all those into account, we know that the poles of
$L^S(2s-\frac{1}{2},\pi,Sym^2\otimes\chi)$ are among the poles of the
normalized Eisenstein series $E^\ast(g,s,f^s)$ because by canceling
the local factors $L_v(r(2s-\frac{1}{2}),
\chi_v^r\omega_v^2)$ for all $v\in S$ in
(\ref{E:integral_representation}), we ave
\[
Z^\ast(\phi,\Theta,f^s)
=L^S(2s-\frac{1}{2},\pi,Sym^2\otimes\chi)Z_\infty(W_\infty,Q_\infty,R^s_\infty)
\prod_{v\in S,\; v<\infty}Z_v(W_v,Q_v,R^s_v)\notag
\]
where
\[ 
Z^\ast(\phi,\Theta,f^s):={\int}\phi(g)\Theta(\kappa(g))E^\ast(\kappa(g),s,f^s)dg.
\] 
Thus the theorem follows.
\end{proof}

\quad\\

We give a proof of Proposition \ref{P:poles_Eisenstein}.

\begin{proof}[Proof of Proposition \ref{P:poles_Eisenstein}] 
The proof is almost identical to the one given by \cite[Theorem
7.4]{BG} except a subtle issue about asymptotic expansions of matrix
coefficients at the archimedean places. 
Since the proof is essentially the same as in \cite{BG} for most of
the part,
we will reproduce only the main points. Moreover since the case of our
main interest is the case for $r=2q$, the case $r=2q+1$ being more
similar to \cite{BG}, we only consider $r=2q$. 
(Also the twisted case for $r=3$ is treated by
\cite{Banks2}.)

First let us note that as we explained at the beginning of the previous
section, the Eisenstein series $E(-,s, f^s)$ on $\GLt_{2q}^{(2)}(\A)$
is simply the restriction of the Eisenstein series on
$\GLt_{2q}(\A)$. Hence one can apply the theory of Eisenstein series
(\cite{MW}) to this case. 

As the
proof in \cite[p.195-196]{BG}, the computation of the poles boils down
to determining the poles of the  intertwining operator
\[
M(s):\ind_{\Qt(\A)}^{\GLt_{2q}(\A)}\theta_{\omega, \omega^{-1}\chi^{-q}}\otimes\delta_Q^s\rightarrow
\ind_{(^{w_0}\MQt(\A))N_{(1,r-1)}(\A)^\ast}^{\GLt_{2q}(\A)}\;^{w_0}(\theta_{\omega,\omega^{-1}\chi^{-q}}
\otimes\delta_Q^s),
\] 
where the induction is NOT normalized and
$w_0=\left(\begin{smallmatrix}&1\\ I_{2q-1}&\end{smallmatrix}\right)$.

For each factorizable section $f^s=\otimes'f_v^s$, we know from
Lemma \ref{L:spherical_section2} that
\[ 
M(s)f^s=\frac{L(r(2s-\frac{1}{2})-r+1,
\chi^r\omega^2)}{L(r(2s-\frac{1}{2}), \chi^r\omega^2)}
\left(\underset{v}{\otimes'}\frac{L_v(r(2s-\frac{1}{2}),
\chi_v^r\omega_v^2)}{L_v(r(2s-\frac{1}{2})-r+1,
\chi_v^r\omega_v^2)}M_v(s)f_v^s\right),
\] 
where $M_v(s)$ is the corresponding local intertwining operator. (Note
that in Lemma \ref{L:spherical_section2}, the induction is normalized,
and hence we need to shift $s$ by $1/2$.)

Hence the poles of the
normalized Eisenstein series $E^\ast(g,s,f^s)$ are the poles of
\begin{align*}
 &L^S(r(2s-\frac{1}{2}), \chi^r\omega^2)M(s)f^s\\
=&L(r(2s-\frac{1}{2})-r+1, \chi^r\omega^2)\underset{v\notin
S}{\otimes'} \frac{L_v(r(2s-\frac{1}{2}),
\chi_v^r\omega_v^2)}{L_v(r(2s-\frac{1}{2})-r+1,\chi_v^r\omega_v^2)}M_v(s)f_v^s\\
&\qquad\qquad\underset{v\in
S}{\otimes'}\frac{1}{L_v(r(2s-\frac{1}{2})-r+1,
\chi_v^r\omega_v^2)}M_v(s)f_v^s.
\end{align*} 
The Hecke $L$-function $L(r(2s-\frac{1}{2})-r+1, \chi^r\omega^2)$ has no
pole unless $\chi^r\omega^2=1$. (Note that if $\chi^r\omega^2= 1$,
then this $L$-function has
poles at $r(2s-\frac{1}{2})-r+1=1$ \ie $s=3/4$ and
$r(2s-\frac{1}{2})-r+1=0$ \ie $s=1/4$. This is why the normalized
Eisenstein series could have a pole at $s=1/4$ and $s=3/4$ for this
case.) Also for almost all $v$, we know from Lemma
\ref{L:spherical_section2} that 
\[ 
\frac{L_v(r(2s-\frac{1}{2}),
\chi_v^r\omega_v^2)}{L_v(r(2s-\frac{1}{2})-r+1,
\chi_v^r\omega_v^2)}M_v(s)f_v^s
\] 
is the spherical section in
$\ind_{^{w_0}\MQt(F_v)N_{(1,r-1)}(F_v)^\ast}^{\GLt_{2q}(F_v)}\;^{w_0}({\theta_{\chi_v\omega_v^2}}\otimes\delta_Q^s)$,
which means that this has no pole.

By taking all those into account, what we have to prove is that for
non-archimedean $v$ the local intertwining operator
\begin{align}\label{E:local_normalized_operator}
&\frac{1}{L_v(r(2s-\frac{1}{2})-r+1, \chi_v^r\omega_v^2)}M_v(s):\\
&\qquad\qquad\ind_{\Qt(F_v)}^{\GLt_{2q+1}(F_v)}
{\theta_{\omega,\omega^{-1}\chi^{-q}}}_v\otimes\delta_Q^s\rightarrow
\ind_{^{w_0}\MQt(F_v)N_{(1,r-1)}(F_v)^\ast}^{\GLt_{2q+1}(F_v)}\;^{w_0}({\theta_{\omega,\omega^{-1}\chi^{-q}}}_v
\otimes\delta_Q^s)\notag
\end{align} 
has no pole, and for archimedean $v$ it has no pole except
those which are canceled by the compensation factor
$L_v(r(2s-\frac{1}{2})-r+1, \chi_v^r\omega_v^2)^{-N_v}$.

First assume $v$ is non-archimedean. But our situation is identical to
\cite{BG}, because our $\theta_{\omega,\omega^{-1}\chi^{-q}}$ is the
same as theirs, except that ours has the twist $\omega^{-1}\chi^{-q}$
by the $\GLt_1$ factor of the parabolic $\Qt$, which does no harm when
one applies the method of \cite{BG}. To show the above
intertwining operator (\ref{E:local_normalized_operator}) has no pole, Bump and
Ginzburg considered the inner product of the induced representations
$\ind_{\Qt(F_v)}^{\GLt_{2q+1}(F_v)}{\theta_{\omega,\omega^{-1}\chi^{-q}}}_v\otimes\delta_Q^s$ 
and
$\ind_{\Qt(F_v)}^{\GLt_{2q+1}(F_v)}{\theta_{\omega,\omega^{-1}\chi^{-q}}}_v\delta_Q^{1-s}$
and reduced the problem to a computation of the asymptotic behavior of the
matrix coefficients given by
\begin{align}\label{E:matrix_coefficient} &\int_{F_v^{r-2}}\int_{F_v}
\left\langle
\theta_v\left(\s\begin{pmatrix}y&&\\&I_{r-2}&\\&&y^{-1}\end{pmatrix}\right)[u_1],
\theta_v\left(\s\begin{pmatrix}Z&-1&\\I_{r-1}&&\\&&1\end{pmatrix}\right)[u_2]
\right\rangle_{\theta'}\\
&\qquad\qquad\qquad\qquad\qquad\qquad\qquad
\qquad\qquad\qquad\qquad\qquad\qquad\qquad
\phi(y,Z)|y|^{r(s-1)}\;dy\,dZ,\notag
\end{align} 
where $\theta_v={\theta_{\omega,\omega^{-1}\chi^{-q}}}_v$, 
$u_i$ is a vector in the space of $\theta'_v$, and
$\phi(y,Z)$ is a Schwartz function on $F_v^{r-1}$. (This is equation
(7.14) of \cite{BG}, and so the details can be found there.) So it
suffices to show that this integral has no pole. It is shown by
\cite{BJ} (Casselman's theorem applied to the metaplectic group) that
the asymptotic of the matrix coefficients as $|y|\rightarrow 0$ is
determined by the Jacquet module of $\theta_v$ along the Borel
subgroup. But since the representation $\theta_v$ is the exceptional
representation for which we know the exact expression for the Jacquet
module by Proposition \ref{P:local_periodicity1}, one can explicitly
compute the asymptotic expansion, which is carried out in \cite[p.200]{BG}.

Assume $v$ is real. (Let us mention that what follows is explained to the author by
N. Wallach, and the author would like to thank him for it.) 
Unlike the non-archimedean case, we do not have such description of
the Jacquet module. But instead, we (and Bump-Ginzburg) use the theory
of Harrish-Chandra (\cite[p.200-201]{BG}). 
The basic idea is essentially
analogous to the non-archimedean case in that one needs to consider the
analogous integral of the matrix coefficient, and instead of
Casselman's theorem, one needs to use the asymptotic expansion of the matrix
coefficient due to Harish-Chandra. Then one obtains
\begin{align*} &\left\langle
\theta_v\left(\s\begin{pmatrix}y&&\\&I_{r-2}&\\&&y^{-1}\end{pmatrix}\right)[u_1],
\theta_v\left(\s\begin{pmatrix}Z&-1&\\I_{r-1}&&\\&&1\end{pmatrix}\right)[u_2]
\right\rangle_{\theta'}\\
&\qquad\qquad\qquad\qquad\qquad\qquad\qquad\qquad
\sim\sum_{n=0}^{\infty}a_n(Z)|y|^{n+(r-2)/4}\omega_v\chi_v^q(y)P(\log|y|)
\quad\text{ as $|y|\rightarrow 0$},
\end{align*} 
where $P(\log|y|)$ is some polynomial in $\log|y|$. (The
reader is advised to compare it with the formula in
\cite[p.201]{BG}. In \cite{BG}, the factor $P(\log|y|)$ is missing.)
Then as in \cite[p.201]{BG} if one carries out the integration, one obtains
the Mellin transform of a function in $y$, which vanishes for $|y|$
large, and the possible poles are determined by the asymptotic as
$|y|\rightarrow 0$. Indeed, for example if $\chi^{1/2}_v$ exists, the possible poles are at
\[ 
s=\frac{3}{4}-\frac{1}{2r}-\frac{\rho}{r}-\frac{n}{r}\quad\text{
for $n\geq0$}
\] 
where $\rho$ is the purely imaginary number so that
$\chi_v^{1/2}\omega_v(y)=(y/|y|)^\epsilon|y|^\rho$ where $\epsilon=0$
or $1$ as in \cite{BG}. (This computation is done by
integration by parts.) 
And those are precisely where the local
archimedean factor $L_v(r(2s-\frac{1}{2})-r+1, \chi_v^r\omega_v^2)$
has poles. However, this theory only tells us the locations of the
possible poles, but the orders of the possible poles cannot be shown to be
simple. Indeed, this theory only tells that the order of each possible
pole is at most
\[ 
(\text{the degree of the polynomial $P$})+1.
\] 
All those issues are explained quite in detail in
\cite[p.361-362]{Wallach}. Hence unless one can show that the degree
of $P$ is $0$, one can not conclude that the possible poles are
canceled with the poles of $L_v(r(2s-\frac{1}{2})-r+1,
\chi_v^r\omega_v^2)$. Although it might be still possible that the polynomial $P$
indeed has degree $0$, at least the author does not know how to show
it. Hence it should be considered that even after the factor
$L_v(r(2s-\frac{1}{2})-r+1, \chi_v^r\omega_v^2)^{-1}$ is multiplied to
the intertwining operator $M_v(s)$, we still have the possible poles
at the above locations. Hence the best we have is the product
\[ L_v(r(2s-\frac{1}{2})-r+1, \chi_v^r\omega_v^2)^{-N_v-1}M_v(s)
\] is holomorphic, where $N_v$ is the degree of the polynomial $P$.

Hence by taking all
those into account, we can show the holomorphy of the product
\[
E^\ast(g,s,f^s)\prod_{v|\infty}L_v(r(2s-\frac{1}{2})-r+1,
\chi_v^r\omega_v^2)^{-N_v}
\]
as in the proposition.
\end{proof}

\quad

From the main theorem (Theorem \ref{T:main}), it is immediate that the
possible poles other than at $s=0$ and $s=1$ come from the poles of
the local archimedean factors $L_v(r(2s-\frac{1}{2})-r+1,
\chi_v^r\omega_v^2)$, which are just gamma functions. Hence we have

\begin{Cor}\label{C:main} 
The (incomplete) twisted symmetric square
$L$-function $L^S(s,\pi,Sym^2\otimes\chi)$ is holomorphic everywhere
in the region $\Re(s)>1-\frac{1}{2r}$ except at $s=1$. Moreover there
is no pole at $s=1$ if $\chi^r\omega^2\neq 1$.
\end{Cor}

The reason we can have our result only for $\Re(s)>1-\frac{1}{2r}$ is
the archimedean issue pointed out above. But we believe this issue can
be resolved and hope to prove

\begin{Con}\label{Con:main}
The (incomplete) twisted
symmetric square $L$-function $L^S(s,\pi,Sym^2\otimes\chi)$ is
holomorphic everywhere except at
$s=0$ and $s=1$. Moreover there is no pole if $\chi^r\omega^2\neq 1$.
\end{Con}

We hope this can be done in our forthcoming paper \cite{Takeda2}.\\

Finally let us note that this corollary does NOT tell us that the $L$-function
$L^S(s,\pi,Sym^2\otimes\chi)$ has a pole at $s=1$ if $\chi^r\omega^2=
1$. We only know it might have a pole at $s=1$, but it might not. We
are not able to determine this. However if $r$ is odd, the following
theorem due to Jacquet-Shalika and Shahidi allows one to tell exactly when
$L^S(s,\pi,Sym^2\otimes\chi)$ has a pole at $s=1$.

\begin{Thm} 
Assume $r$ is odd. Then the (complete) twisted exterior
square $L$-function $L(s,\pi,\wedge^2\otimes\chi)$ is non-zero
holomorphic at $s=1$.
\end{Thm}
\begin{proof} 
The non-vanishing part is the main theorem of \cite{Sh97}, and
the holomorphy is \cite[Theorem 9.6.2]{JS}.
\end{proof}

This theorem implies
\begin{Cor}\label{C:main2} Assume $r$ is odd. Then the $L$-function
$L^S(s,\pi,Sym^2\otimes\chi)$ has a pole at $s=1$ if and only if
$\check{\pi}=\pi\otimes\chi$, where $\check{\pi}$ is the
contragredient of $\pi$.
\end{Cor}
\begin{proof} Recall
\[
L^S(s,\pi\times\pi\otimes\chi)=L^S(s,\pi,\wedge^2\otimes\chi)L^S(s,\pi,Sym^2\otimes\chi),
\] and the Rankin-Selberg $L$-function
$L^S(s,\pi\times\pi\otimes\chi)$ has a pole at $s=1$ if and only if
$\check{\pi}=\pi\otimes\chi$. Hence the corollary follows from the
above theorem.
\end{proof}

\quad

%%%%%%%%%%%%%%%%%%%%%%%%%%%%%%%%%%%%%%%%%%%%%%%%%%%%%%%%%%%%%%%%%%%

\appendix\section{\bf Metaplectic tensor product }\label{S:tensor_product}

%%%%%%%%%%%%%%%%%%%%%%%%%%%%%%%%%%%%%%%%%%%%%%%%%%%%%%%%%%%%%%%%%%%

In this appendix, we will recall the notion of metaplectic tensor
product for $\GLtt_r$ both locally and globally. For the local case,
if one uses the block-compatible cocycle $\sigma_r$, the formulation of
metaplectic tensor product is done in several places. (See \cite{Banks, Kable2,
  Mezo}.) But since we use our $\tau_r$, which works both for the
local and global cases, we need another formulation. Let us mention
that this appendix is a portion of 
\cite{Takeda1} in which we developed the theory of metaplectic tensor
products for automorphic representations of the $n$-fold cover of
$\GL_r(\A)$, and in the interest of space, we only recall the basic
facts necessary for our purposes and we will occasionally omit the
detailes of the proofs, all of which are
available in \cite{Takeda1}.

Let $P$ be a parabolic subgroup of $\GL_r$ whose Levi is
\[
M_P=\GL_{r_1}\times\cdots\times\GL_{r_k}.
\]
Of course we assume $M_P$ sits in $\GL_r$ diagonally. Let us denote by $\MPt$
the metaplectic preimage of $M_P$, and write
\[
\MPt=\GLt_{r_1}\timest\cdots\timest\GLt_{r_k},
\]
where the group structure of $\MPt$ is defined via the restriction of
the cocycle $\tau_r$.

%%%%%%%%%%%%%%%%%%%%%%%%%%%%%%%%%%%%%%%%%%%%%%%%%%%%%%%%%%%%%%%%%%%

\subsection{\bf The group $\cMPt$}\label{S:cMPt}

%%%%%%%%%%%%%%%%%%%%%%%%%%%%%%%%%%%%%%%%%%%%%%%%%%%%%%%%%%%%%%%%%%%

One difficulty to work with $\tau_r$ is that it is not known that it
is block-compatible unless $r=2$. To get around it, let us define a cocycle
\[
\tau_P: M_P\times M_P\rightarrow\{\pm 1\},
\]
both locally and globally, by
\[
\tau_P(\begin{pmatrix}g_1&&\\ &\ddots&\\ &&g_k\end{pmatrix},
\begin{pmatrix}g'_1&&\\ &\ddots&\\ &&g'_k\end{pmatrix})
=\prod_{i=1}^k\tau_{r_i}(g_i,g_i')\prod_{1\leq i<j\leq k}(\det(g_i), \det(g_j')),
 \]
where $(-,-)$ is the local or global Hilbert symbol.
Note that the definition makes sense both locally and
globally. Moreover the global $\tau_P$ is the product of the local
ones.

We define the group $\cMPt$ to be 
\[
\cMPt=M_P\times\{\pm 1\}
\]
as a set and the group structure is
given by $\tau_P$. The superscript $^c$ is for
``compatible''. One advantage to work with $\cMPt$ is that each
$\GLt_{r_i}$ embeds into $\cMPt$ via the natural map
\[
(g_i,\xi)\mapsto(\begin{pmatrix}I_{r_1+\cdots+r_{i-1}}&&\\ &g_i&\\
    &&I_{r_{i+1}+\cdots+r_k}\end{pmatrix}, \xi).
\]
Or rather, the cocycle $\tau_P$ is so chosen that we have this
embedding.

Also recall our notation
\[
M_P^{(2)}=\GL_{r_1}^{(2)}\times\cdots\times\GL_{r_k}^{(2)},
\]
and
\[
\MPtt=\GLtt_{r_1}\timest\cdots\timest\GLtt_{r_k}.
\]
We define $\cMPtt$ analogous to $\cMPt$, namely the group structure of
$\cMPtt$ is defined via the cocycle $\tau_P$. Of course, $\cMPtt$ is a
subgroup of $\cMPt$. Note that each
$\GLtt_{r_i}$ naturally embeds into $\cMPtt$ as above.

\begin{Lem}
The subgroups $\GLtt_{r_i}$ and $\GLtt_{r_j}$ in $\cMPtt$ commute
pointwise for $i\neq j$.
\end{Lem}
\begin{proof}
Locally or globally, it suffices to show
$\tau_P(\gt_i,\gt_j)=\tau_P(\gt_j,\gt_i)$. But since the global
$\tau_r$ is the product of local ones, it suffices to show the local
case. So assume our groups are over a local field. By the relation
between $\tau_P$ and $\sigma_r$, it suffices to show
$\sigma_r(\gt_i,\gt_j)=\sigma_r(\gt_j,\gt_i)$. But this follows from 
the block-compatibility of the 2-cocycle $\sigma_r$ as in
(\ref{E:compatibility}). (See also \cite[p.141]{BG}.)
\end{proof}

\begin{Lem}
There is a surjection
\[
\GLtt_{r_1}\times\cdots\times\GLtt_{r_k}\rightarrow\; \cMPtt
\]
given by the map
\[
((g_1,\xi_1),\dots,(g_k,\xi_k))\mapsto
(\begin{pmatrix}g_1&&\\ &\ddots&\\ &&g_k\end{pmatrix},
\xi_1\cdots\xi_k),
\]
whose kernel is
\[
\mathcal{K}_P:=\{((1,\xi_1),\dots,(1,\xi_k)):\xi_1\cdots\xi_k=1\},
\]
so that $\cMPtt\cong
\GLtt_{r_1}\times\cdots\times\GLtt_{r_k}/\mathcal{K}_P$.
\end{Lem}
\begin{proof}
The above lemma together with the block-compatibility of $\tau_P$
guarantees that the map is indeed a group homomorphism. The
description of the kernel is immediate.
\end{proof}

Note that for the group $\MPt$, the group structure
is defined by the restriction of $\tau_r$ to $M_P\times M_P$, and hence
each $\GLt_{r_i}$ might not embed into $\GLt_r$ in
the natural way because of the possible failure of the
block-compatibility of $\tau_r$ unless $r=2$.
To make explicit the relation between $\cMPt$ and $\MPt$, the discrepancy between
$\tau_r|_{M_P\times M_P}$ (which we denote simply by $\tau_r$) and
$\tau_P$ has to be clarified.

\quad

\noindent{\bf Local case:}

Assume $F$ is local. Then we have
\begin{align*}
&\tau_P(\begin{pmatrix}g_1&&\\ &\ddots&\\ &&g_k\end{pmatrix},
\begin{pmatrix}g'_1&&\\ &\ddots&\\ &&g'_k\end{pmatrix})
\\
=&\sigma_r(\begin{pmatrix}g_1&&\\ &\ddots&\\ &&g_k\end{pmatrix},
\begin{pmatrix}g'_1&&\\ &\ddots&\\ &&g'_k\end{pmatrix})
\prod_{i=1}^ks_{r_i}(g_i)s_{r_i}(g_i')/s_{r_i}(g_ig_i'),
\end{align*}
so $\tau_P$ and $\sigma_r|_{M_P\times M_P}$ are cohomologous via the function
$\prod_{i=1}^ks_{r_i}$. Here recall from Section \ref{S:group} that the map
$s_{r_i}:\GL_{r_i}\rightarrow\{\pm 1\}$ relates $\tau_{r_i}$ with
$\sigma_{r_i}$ by
\[
\sigma_{r_i}(g_i,g_i')=\tau_{r_i}(g_i,g_i')\cdot\frac{s_{r_i}(g_i,g_i')}{s_{r_i}(g_i)s_{r_i}(g_i')},
\]
for $g_i,g_i'\in\GL_{r_i}$. Moreover if the residue characteristic is odd, $s_{r_i}$ is chosen to
be ``canonical'' in the sense that (\ref{E:canonical_section}) is satisfied.

The block-compatibility of $\sigma_r$ implies
\[
\tau_r(m, m')\cdot\frac{s_r(mm')}{s_r(m)s_r(m')}
=\tau_P(m,m')\prod_{i=1}^k\frac{s_{r_i}(g_1g_2)}{s_{r_i}(g_1)s_{r_i}(g_2)},
\]
for $m=\begin{pmatrix}g_1&&\\ &\ddots&\\ &&g_k\end{pmatrix}$ and
$m'=\begin{pmatrix}g_1'&&\\ &\ddots&\\ &&g_k'\end{pmatrix}$. Hence if
we define $\hat{s}_P:M_P\rightarrow\{\pm1\}$ by
\[
\hat{s}_P(m)=\frac{\prod_{i=1}^ks_{r_i}(g_i)}{s_r(m)}, 
\]
we have
\begin{equation}\label{E:s_hat}
\tau_P(m, m')=\tau_r(m, m')\cdot
\frac{\hat{s}_P(m)\hat{s}_P(m')}{\hat{s}_P(mm')}.
\end{equation}

Therefore we have the isomorphism
\[
\varphit_P:\cMPt\rightarrow\MPt,\quad (m,\xi)\mapsto (m, \hat{s}_P(m)\xi).
\]

\quad

An important fact about the map $\hat{s}_P$ is
\begin{Lem}\label{L:s_hat=1}
Assume $F$ is non-archimedean of
odd residual characteristic. Then for
all $k\in M_P(\OF)$, we have $\hat{s}_P(k)=1$.
\end{Lem}
\begin{proof}
This is \cite[Lemma 3.5]{Takeda1} and essentially follows from the
``canonicality'' of $s_r$ and $s_{r_i}$, so that $s_r$ has been chosen to satisfy
$s_r={\sss_r}|_{l(\GL_r(\OF))}$, where $\sss_r$ is the map on
$\G_{r}(F)$ that makes the diagram (\ref{E:canonical_diagram})
commute, and from the fact that the
cocycle for $\G_r$ is block-compatible for a very strong sense as in
\cite[Lemma 5, Theorem 7 \S 2]{BLS}. See \cite{Takeda1} for the detail.
\end{proof}

\quad

\noindent{\bf Global case:}

Assume $F$ is global. Define $\hat{s}_P:
M_P(\A)\rightarrow\{\pm1\}$ by
\[
\hat{s}_P(\prod_vm_v):=\prod_v{\hat{s}_{P_v}}(m_v)
\]
for $\prod_vm_v\in M_P(\A)$. The product is finite thanks to Lemma
\ref{L:s_hat=1}. Since both of the cocycles $\tau_r$ and $\tau_P$ are
the products of the corresponding local ones, one can
see that the relation (\ref{E:s_hat}) holds globally as well.

Thus analogously to the local case, we have
the isomorphism
\[
\varphit_P:\; \cMPt(\A)\rightarrow\MPt(\A),\quad
(m,\xi)\mapsto (m, \hat{s}_P(m)\xi).
\]

\quad

\begin{Lem}\label{L:splitting_cMPt}
The splitting of $M_P(F)$ into $\cMPt(\A)$ is given by
\[
\s_P:M_P(F)\rightarrow\;\cMPt(\A),\quad
\begin{pmatrix}g_1&&\\ &\ddots&\\ &&g_k\end{pmatrix}\mapsto 
(\begin{pmatrix}g_1&&\\ &\ddots&\\ &&g_k\end{pmatrix},\;
\prod_{i=1}^k s_i(g_i)^{-1}).
\]
\end{Lem}
\begin{proof}
For each $i$ the splitting
$\s_{r_i}:\GL_{r_i}(F)\rightarrow\GLt_{r_i}(\A)$ is given by
$g_i\mapsto(g_i,\;s_{r_i}(g_i)^{-1})$, where $\GLt_{r_i}(\A)$ is
defined via the cocycle $\tau_{r_i}$. Then the lemma follows by the
block-compatibility of $\tau_P$ and the product formula for the
Hilbert symbol.
\end{proof}

This splitting is related to the splitting
$\s:\GL_r(F)\rightarrow\GL_r(\A)$ by
\begin{Prop}\label{P:diagram}
We have the following commutative diagram:
\[
\xymatrix{\cMPt(\A)\; \ar@{^{(}->}[r]^{\varphit_P}&\GLt_r(\A)\\
M_P(F)\;\ar@{^{(}->}[r]\; \ar[u]^{\s_P}&\GL_r(F)\ar[u]_{\s_r}.
}
\]
\end{Prop}
\begin{proof}
Note that for the elements in $\GL_r(F)$, all of $s_{r_i}$ and $s_r$
are defined globally, and then the proposition follows from the
definition of $\s_P$ and $\s_r$.
\end{proof}
\begin{comment}
\begin{proof}
For $m=\begin{pmatrix}g_1&&\\ &\ddots&\\ &&g_k\end{pmatrix}
\in M_P(F)$, we have
\[
\varphit_P(\s_P(m))=\varphit_P(m, \prod_{i=1}^ks_{r_i}(g_i)^{-1})
=(m, \hat{s}_P(m)\prod_{i=1}^ks_{r_i}(g_i)^{-1})
=(m, s_r(m)^{-1})=\s_r(m).
\]
Here note that for the elements in $\GL_r(F)$, all of $s_{r_i}$ and $s_r$
are defined globally. 
\end{proof}
\end{comment}

This proposition implies
\begin{Cor} \label{C:diagram}
Assume $\pi$ is an automorphic representation of
$\cMPt(\A)$. The representation of
$\MPt(\A)$ defined by $\pi\circ\varphit_P^{-1}$ is
also automorphic.
\end{Cor}
\begin{proof}
If $\pi$ is realized in a space $V$ of automorphic forms on
$\cMPt(\A)$, then $\pi\circ\varphit_P^{-1}$ is realized in the
space of functions of the form $f\circ\varphit_P^{-1}$ for $f\in
V$. Then the automorphy follows from the commutativity of the diagram
in the above lemma.
\end{proof}

%%%%%%%%%%%%%%%%%%%%%%%%%%%%%%%%%%%%%%%%%%%%%%%%%%%%%%%%%%%%%%%%%%%

\subsection{\bf Metaplectic tensor product}\label{S:cMPt}

%%%%%%%%%%%%%%%%%%%%%%%%%%%%%%%%%%%%%%%%%%%%%%%%%%%%%%%%%%%%%%%%%%%

We are ready to define the notion of metaplectic tensor
product. We treat both local and global cases at the same time. 
Let $\pi_1,\dots,\pi_k$ be irreducible admissible representations of
$\GLtt_{r_1},\dots,\GLtt_{r_k}$, respectively, where each $\pi_i$ is
realized in the space $V_i$. Further assume each $\pi_i$ is genuine.
Consider the usual tensor product
representation $\pi_1\otimes\cdots\otimes\pi_k$ of the direct product
$\GLtt_{r_1}\times\cdots\times\GLtt_{r_k}$ realized in the space
$V_1\otimes\cdots\otimes V_k$. Since each $\pi_i$ is genuine, the
kernel $\mathcal{K}_P$ of the above lemma acts trivially on
$\pi_1\otimes\cdots\otimes\pi_k$. Thus this tensor product
representation descends to a representation of
$\cMPtt$, which we denote by
\[
\pi_1\otimest\cdots\otimest\pi_k,
\]
and we call it the {\it metaplectic tensor product}
representation of $\cMPtt$. Let us emphasize that the space of the metaplectic tensor
product representation is the same as that of the tensor product.

Of course one can pullback the metaplectic tensor product
$\pi_1\otimest\cdots\otimest\pi_k$ of $\cMPtt$ to a representation of $\MPtt$ via
the map $\varphit_P^{-1}$, which we often denote by the same symbol
$\pi_1\otimest\cdots\otimest\pi_k$, when there is no danger of
confusion, and we call it the metaplectic tensor product representation of $\MPtt$.

\begin{Prop}\label{P:tensor_product_automorphy}
Assume $F$ is global, and $\pi_1,\dots,\pi_k$ are genuine irreducible automorphic
representations of $\GLtt_{r_1}(\A)$, $\dots$,
$\GLtt_{r_k}(\A)$, respectively. Then the metaplectic tensor
product representation $\pi_1\otimest\cdots\otimest\pi_k$ of
$\cMPtt(\A)$ is aslo automorphic.
\end{Prop}
\begin{proof}
This is \cite[Proposition 5.2]{Takeda1}. The proof is quite
straightforward by viewing each function $f_1\otimes\cdots\otimes f_k\in
\pi_1\otimes\cdots\otimes \pi_k$ naturally as a function on $\cMPtt(\A)$ by
\[
(f_1\otimes\cdots\otimes f_k)(
\begin{pmatrix}g_1&&\\ &\ddots&\\ &&g_k\end{pmatrix},\xi)
=\xi f_1(g_1,1)\cdots f_k(g_k,1).
\]
The automorphy follows from the definition of $\s_P$ and $\s_{r_i}$
along with the block compatibility of $\tau_P$.
\end{proof}

\end{document}